\documentclass[ejs,preprint]{imsart}

\usepackage{multicol}
\usepackage{multirow}

\usepackage[latin1]{inputenc}
\usepackage{amsfonts}
\usepackage{amssymb}
\usepackage{amsmath}
\usepackage{amsthm}
\usepackage{latexsym}
\usepackage{verbatim}
\usepackage{epsfig}
\usepackage{color}
										
\usepackage{graphicx}
\usepackage{amsbsy}
\usepackage{amscd}
\usepackage{bm}

\newcommand{\R}{{\mathbb R}}

\DeclareMathOperator*{\E}{\mathbf{E}}
\DeclareMathOperator*{\argmin}{argmin}

\DeclareMathOperator*{\p}{\mathbf{P}}
\providecommand{\wt}[1]{\widetilde{#1}}
\providecommand{\wh}[1]{\widehat{#1}}

\providecommand{\norm}[1]{\left \lVert#1 \right \rVert}
\providecommand{\nnorm}[1]{ \lVert#1 \rVert}
\newcommand{\scp}[2]{\left\langle#1, #2\right\rangle}
\newcommand{\nscp}[2]{\langle#1, #2\rangle}

\providecommand{\mc}[1]{\mathcal#1}
\providecommand{\T}{\top}

\newcommand{\blanco}[1]{  }

\newcommand{\deriv}[3]{%
\ifthenelse{#1 = 1}{\frac{d\,#2}{d\,#3}}{\frac{d^{{#1}} #2}{d{#3}^{{#1}}}}
}

\newcommand{\partials}[3]{%
\ifthenelse{#1 = 1}{\frac{\partial\,#2}{\partial\,#3}}{\frac{\partial^{#1}
    #2}{\partial#3^{#1}}}
} 

\def\su{\sum_{i=1}^n}

\def \invcoloneq{=\mathrel{\mathop:}}
\def \eps{\varepsilon} 

\def \gec{\succeq}
\def \lec{\preceq} 
   
\def \fon{\text{\scriptsize{$\frac{1}{n}$}}}
\def \ffon{\text{\footnotesize{$\frac{1}{n}$}}}

\newtheorem{theo}{Theorem}
\newtheorem{propo}{Theorem}
\newtheorem{definitio}{Theorem}
\newtheorem{conditio}{Theorem}

\newtheorem{lemmachen}{Theorem}

\newtheorem{remark}{Rem}

\newtheorem{theorem}[theo]{Theorem}
  \newtheorem{defn}[definitio]{Definition}

  \newtheorem{prop}[propo]{Proposition}
  
  \newtheorem{lemma}[lemmachen]{Lemma}

\newtheorem{cond}[conditio]{Condition}

\newtheorem{remarque}[remark]{Remark}

\newenvironment{bew}{\begin{proof}[Proof]}{\end{proof}}

\def\plasso{\wh{\beta}^{\ell_{1, \gec}}}

\RequirePackage[OT1]{fontenc}
\RequirePackage{amsthm,amsmath}
\RequirePackage[numbers]{natbib}
\RequirePackage{hypernat}

\providecommand{\se}[1]{{\footnotesize$\pm$#1}}

\usepackage{bbm}
\def\one{\mathbbm{1}}
\startlocaldefs
\numberwithin{equation}{section}
\theoremstyle{plain}

\endlocaldefs

\begin{document}

\begin{frontmatter}
  
\title{Non-negative least squares for high-dimensional linear
  models: consistency and sparse recovery without regularization}

\runtitle{NNLS in high-dimensions}

\begin{aug}
\author{\fnms{Martin}
  \snm{Slawski}\thanksref{t1}\ead[label=e1]{ms@cs.uni-saarland.de}}
\and
\author{\fnms{Matthias} \snm{Hein} \ead[label=e2]{hein@cs.uni-saarland.de}}


\thankstext{t1}{Supported by the Cluster of Excellence 'Multimodal Computing
        and Interaction' (MMCI) of Deutsche Forschungsgemeinschaft (DFG)}
\runauthor{Slawski and Hein}

\affiliation{Saarland University}

\address{Saarland University\\
Campus E 1.1\\
Saarbr\"ucken, Germany\\
e-mail: \printead*{e1} ; \printead*{e2}}

\end{aug}

\begin{abstract}
Least squares fitting is in general not useful for high-dimensional linear
models, in which the number of predictors is of the same or even larger
order of magnitude than the number of samples. Theory developed in recent
years has coined a paradigm according to which sparsity-promoting regularization is
regarded as a necessity in such setting. Deviating from this paradigm, we show that
non-negativity constraints on the regression coefficients may be similarly
effective as explicit regularization if the design matrix has additional
properties, which are met in several applications of non-negative
least squares (NNLS). We show that for these designs, the 
performance of NNLS with regard to prediction and estimation is comparable
to that of the lasso. We argue further that in specific cases, NNLS may have a better $\ell_{\infty}$-rate
in estimation and hence also advantages with respect to support recovery when 
combined with thresholding. From a practical point of view, NNLS does not
depend on a regularization parameter and is hence easier to use.  
\end{abstract}



\end{frontmatter}

\section{Introduction}

Consider the linear regression model
\begin{equation}\label{eq:linearmodel}
y = X \beta^{\ast} + \eps,
\end{equation}
where $y$ is a vector of observations, $X \in \R^{n \times p}$ a design
matrix, $\eps$ a vector of noise and $\beta^{\ast}$ a vector of
coefficients to be estimated. Throughout this paper, we are concerned with a
high-dimensional setting in which the number of unknowns $p$ is at least of the same
order of magnitude as the number of observations $n$, i.e.~$p=\Theta(n)$ or even $p
\gg n$, in which case one cannot hope to recover the target $\beta^{\ast}$ if
it does not satisfy one of various kinds of sparsity constraints, the simplest
being that $\beta^{\ast}$ is supported on $S = \{j: \beta_j^{\ast} \neq 0
\}$, $|S| = s < n$. In this paper, we additionally assume that $\beta^{\ast}$
is non-negative, i.e.~$\beta^{\ast} \in \R_+^p$. This constraint is
particularly relevant when modelling non-negative data, which emerge e.g.~in the
form of pixel intensity values of an image, time measurements, histograms or
count data, economical quantities such as prices, incomes and growth
rates. Non-negativity constraints occur naturally in numerous deconvolution and unmixing problems in
diverse fields such as acoustics \cite{LinLeeSaul2004}, astronomical imaging
\cite{Bardsley2006}, hyperspectral imaging \cite{Szlam2010}, genomics \cite{Speed2000}, proteomics
\cite{SlawskiHein2010}, spectroscopy \cite{Donoho1992} and network tomography \cite{Meinshausen2013}; see
\cite{ChenPlemmons2007} for a survey. As reported in these references,
non-negative least squares (NNLS) yields at least reasonably good, sometimes
even excellent results in practice, which may seem surprising in view of the
simplicity of that approach. The NNLS problem is given by the quadratic program
\begin{equation}\label{eq:nnls}
\min_{\beta \gec 0} \frac{1}{n} \norm{y - X \beta}_2^2,
\end{equation}
which is a convex optimization problem that can be solved efficiently
\cite{Kim2010}. A minimizer $\wh{\beta}$ of \eqref{eq:nnls} will be referred
to as an NNLS estimator. Solid theoretical support for the empirical success of NNLS from a statistical
perspective scarcely appears in the literature. An early reference is \cite{Donoho1992} dating back already two decades. The authors show
that, depending on $X$ and the sparsity of $\beta^{\ast}$, NNLS may
have a 'super-resolution'-property that permits reliable estimation of
$\beta^{\ast}$. Rather recently, sparse recovery of non-negative signals in a
noiseless setting ($\eps = 0$) has been studied in \cite{Bruckstein2008, DonohoTanner2010, WangTang2009, WangXuTang2011}. One
important finding of this body of work is that non-negativity constraints alone may
suffice for sparse recovery, without the need to employ sparsity-promoting
$\ell_1$-regularization as usually. On the other hand, it remains unclear
whether similar results continue to hold in a more realistic noisy setup. At first
glance, the following considerations suggest a negative answer. The fact that 
NNLS, apart from sign constraints, only consists of a fitting
term, fosters the intuition that NNLS is prone to overfitting, specifically in
a high-dimensional setting. Usefulness of NNLS in such setting appears to be
in conflict with the well-established paradigm that a regularizer is necessary
to prevent over-adaptation to noise and to enforce desired structural
properties of the solution, like sparsity of the vector of coefficients. As
one of the main contributions of the present paper, we characterize the
\emph{self-regularizing property} which NNLS exhibits for a certain class of designs that
turn out to be tailored to the non-negativity constraint, thereby
disentangling the apparent conflict above and improving the understanding of
the empirical success of NNLS. More precisely, we show that for these designs,
NNLS is rather closely related to the non-negative lasso problem
\begin{equation}\label{eq:nnlasso}
\min_{\beta \gec 0} \frac{1}{n} \norm{y - X \beta}_2^2 + \lambda \bm{1}^{\T}
\beta, \quad \lambda > 0, 
\end{equation}
the sign-constrained counterpart to the popular lasso problem
\cite{Tib1996}, which is given by   
\begin{equation}\label{eq:lasso}
\min_{\beta \gec 0} \frac{1}{n} \norm{y - X \beta}_2^2 + \lambda \nnorm{\beta}_1, \quad \lambda > 0, 
\end{equation}
where $\lambda > 0$ is a tuning parameter. Elaborating on the relation between
NNLS and the non-negative lasso, we establish that for the aforementioned class
of designs, NNLS achieves comparable performance with regard to prediction,
estimation of $\beta^*$ and support recovery when combined with subsequent
thresholding. We here refer to the monograph \cite{Buhlmann2011}, the
survey \cite{Geer2009} and the retrospective \cite{Tib2011} for an overview on
the appealing performance guarantees established for the lasso in the last decade, which have, along 
with favourable computational properties, contributed to its enormous
popularity in data analysis. In this paper, we argue
that, in view of both theoretical considerations and empirical evidence, improvements of NNLS over the
(non-negative) lasso are possible, even though they are limited to a
comparatively small set of designs. Differences in performance arise from the
bias of the $\ell_1$-regularizer in \eqref{eq:nnlasso} and \eqref{eq:lasso}
that is responsible for a in general sub-optimal rate for estimation of
$\beta^*$ in the $\ell_{\infty}$-norm \cite{Zhang2009b}. Unlike for NNLS, a
tuning parameter needs to be specified for the (non-negative) lasso,
as it is necessary for all regularization based-methods. Selection of the
tuning parameter by means of cross-validation increases the computational burden and may be error-prone if done by a grid
search, since the grid could have an unfavourable range or a too small number of
grid points. Theoretical results on how to set the regularization parameter
are often available, but require a sufficient degree of acquaintance with existing
literature and possibly also knowledge of quantities such as the noise
level. The last issue has been withstanding problem of the lasso until
recently. In \cite{Belloni2011} and \cite{SunZhang2012} two related variants
of the lasso are proposed that have similar theoretical guarantees, while
the tuning parameter can be set without knowledge of the noise level. On the
other hand, NNLS is directly applicable, since it is free of tuning
parameters. 


%

\emph{Outline and contributions of the paper.}
The paper significantly extends a previous conference publication
\cite{SlawskiHein2011nips}, which contains the first systematic analysis of
NNLS in a high-dimensional setting. Recently, Meinshausen
\cite{Meinshausen2013} has
independently studied the performance of NNLS in such a
setting. His work is related to ours in Section
\ref{subsec:oracleprediction}. The paper is organized as follows. In Section  
\ref{sec:nnlsdoesnotoverfit} we work out the self-regularizing property 
NNLS may have in conjunction with certain design matrices. Equipped with that
property, a bound on the prediction error is stated that resembles 
a corresponding 'slow rate' bound available for the lasso. Developing further the connection
to the lasso, we use techniques pioneered in Bickel \emph{et al.} \cite{Bickel2009} to prove bounds
on the estimation error $\nnorm{\wh{\beta} - \beta^*}_q$ in
$\ell_q$-norm, $q \in [1,2]$, and an improved bound on the prediction error of
NNLS in Section \ref{subsec:oracleprediction}. In Section
\ref{sec:sparserecovery}, we finally provide bounds on the sup-norm
error $\nnorm{\wh{\beta} - \beta^*}_{\infty}$ of NNLS. Hard
thresholding of $\,\wh{\beta}$ is proposed for sparse recovery, and a data-driven 
procedure for selecting the threshold is devised. Section
\ref{sec:sparserecovery} also contains a discussion of advantages and
disadvantages of NNLS relative to the non-negative lasso. In Section
\ref{sec:designs}, we have a closer look at several designs, which satisfy the
conditions required throughout the paper. In Section \ref{sec:empiricalresults}, we discuss the  
empirical performance of NNLS in deconvolution and sparse recovery in
comparison to standard methods, in particular to the non-negative lasso. The
appendix contains most of the proofs, apart from those that have been placed
into a supplement.  

\emph{Notation.}
We denote the usual $\ell_q$-norm by
$\norm{\cdot}_q$. The cardinality of a set is denoted by $|\cdot|$. Let $J,K \subseteq \{1,\ldots,m \}$ be index
sets. For a matrix $A \in \R^{n \times m}$, $A_J$ denotes the matrix one obtains
by extracting the columns corresponding to $J$. For $j=1,\ldots,m$, $A_j$
denotes the $j$-th column of $A$. The
matrix $A_{JK}$ is the sub-matrix of $A$ by extracting rows in $J$ and columns
in $K$. Likewise, for $v \in \R^m$, $v_J$ is the
sub-vector corresponding to $J$. Identity matrices and vectors of ones of
varying dimensions are denoted by $I$ respectively $\bm{1}$. The symbols $\lec, \gec$ and $\prec,\succ$ denote
componentwise inequalities and componentwise strict inequalities,
respectively. In addition, for some matrix $A$, $A \gec a$ means that all
entries of $A$ are at least equal to a scalar $a$. The
non-negative orthant $\{x \in \R^m: x \gec 0 \}$ is denoted by
$\R_+^m$. The standard simplex in $\R^m$, that is the set $\{x \in \R_+^m: \; \sum_{j = 1}^m x_j
= 1 \}$ is denoted by $T^{m - 1}$. Lower and uppercase c's like $c,c', c_1$ and $C,C',C_1$ etc.~denote positive constants not depending on the
sample size $n$ whose values may differ from line to line. In general, the
positive integers $p = p_n$ and $s = s_{n}$ depend on $n$. Landau's symbols
are denoted by $o(\cdot),O(\cdot), \Theta(\cdot), \Omega(\cdot)$. Asymptotics is to be understood w.r.t.
a triangular array of observations $\{(X_n, y_n), X_n \in \R^{n \times p_n}
\}, \; \; n=1,2,\ldots$. 
\\ 
\emph{Normalization.} If not stated otherwise, the design matrix $X$ is considered as
deterministic, having its columns normalized such that $\norm{X_j}_2^2 = n$,
$j=1,\ldots,p$.\\ 
\emph{General linear position.}
We say that the columns of $X$ are in general linear position in $\R^n$ if the following condition (GLP) holds
\begin{equation}\label{eq:glp}
(\text{GLP}):  \quad  \forall J \subset \{1,\ldots, p\}, \; |J| = \min\{n,p\} \; \; \forall \lambda \in \R^{|J|}   \; \;
X_J \lambda = 0 \; \; \Longrightarrow \; \lambda = 0,  
\end{equation}
where for $J \subseteq \{1,\ldots,p\}$, $X_J$ denotes the submatrix of $X$
consisting of the columns corresponding to $J$.

\section{Prediction error: a bound for 'self-regularizing' designs}\label{sec:nnlsdoesnotoverfit}

The main result of the following section is a bound on the excess error of
NNLS that resembles the so-called slow rate bound of the lasso
\cite{GreenshteinRitov2004, Bartlett2012, Buhlmann2011}. In contrast to the
latter, the corresponding bound for NNLS only holds for a certain class of designs. We
first show that extra conditions on the design are in fact necessary to obtain
such bound. We then introduce a condition on $X$ that we refer to as 'self-regularizing
property' which is sufficient to establish a slow rate bound for NNLS. The term 'self-regularization'
is motivated from a resulting decomposition of the least squares objective into a
modified fitting term and a quadratic term that plays a role similar to an
$\ell_1$-penalty on the coefficients. This finding provides an intuitive  
explanation for the fact that NNLS may achieve similar performance than the
lasso, albeit no explicit regularization is employed. 




\subsection{A minimum requirement on the design for non-negativity
 being an actual constraint}

In general, the non-negativity constraints in \eqref{eq:nnls} may not be
meaningful at all, given the fact that any least squares problem can be
reformulated as a non-negative least squares problem with an augmented design
matrix $[X \; -X]$. More generally, NNLS can be as ill-posed as least squares if 
the following condition $(\mc{H})$ does not hold.   
\begin{equation}\label{eq:halfspaceconstraint}
(\mc{H}):\;\; \exists w \in \R^n \; \text{such that} \; X^{\T} w \succ 0.
\end{equation}
Condition $(\mc{H})$ requires the columns of $X$ be contained in the interior
of a half-space containing the origin. If $(\mc{H})$ fails to hold, $0 \in
\text{conv}\{X_j \}_{j = 1}^p$, so that there are infinitely many minimizers of the NNLS problem \eqref{eq:nnls}. 
If additionally $p > n$ and the columns of $X$ are in general linear
position (condition $(\text{GLP})$ in \eqref{eq:glp}), $0$ must be in the
interior of $\text{conv}\{X_j \}_{j = 1}^p$. It then follows that $\mc{C} = \R^n$, where $\mc{C} =
\text{cone} \{X_j \}_{j = 1}^p$ denotes the polyhedral cone generated by the
columns of $X$. Consequently, the non-negativity constraints become vacuous
and NNLS yields perfect fit for any observation vector $y$. In light of this,
$(\mc{H})$ constitutes a necessary condition for a possible improvement 
of NNLS over ordinary least squares.

\subsection{Overfitting of NNLS for orthonormal design}\label{sec:orthonormaldesign} 

Since NNLS is a pure fitting approach, over-adaptation to noise is a natural
concern. Resistance to overfitting can be quantified in terms of $\frac{1}{n}
\nnorm{X \wh{\beta}}_2^2$ when $\beta^* = 0$ in \eqref{eq:linearmodel}. It
turns out that condition $(\mc{H})$ is not sufficient to ensure that $\frac{1}{n} \nnorm{X
  \wh{\beta}}_2^2 = o(1)$ with high probability, as can be seen from studying
orthonormal design, i.e. $X^{\T} X = n I$. Let $y = \eps$ be a standard Gaussian
random vector. The NNLS estimator has the closed form expression
\begin{equation*}
\wh{\beta}_j = \max \{X_j^{\T} \eps, 0\} / n, \; \; j=1,\ldots,p, 
\end{equation*}
so that the distribution of each component of $\wh{\beta}$ is given by a
mixture of a point mass $0.5$ at zero and a half-normal distribution. We conclude that $\frac{1}{n} \nnorm{X
  \wh{\beta}}_2^2 = \frac{1}{n} \nnorm{\wh{\beta}}_2^2$ is of the order $\Omega(p/n)$
with high probability. The fact that $X$ is orthonormal is much stronger than the obviously necessary
half-space constraint $(\mc{H})$. In fact, as rendered more precisely
in Section \ref{sec:designs}, orthonormal design turns out to be at the edge of the set of designs still
leading to overfitting. 

\subsection{A sufficient condition on the design preventing NNLS from overfitting}\label{sec:designsselfreg}

We now present a sufficient condition $X$ has to satisfy so that overfitting
is prevented. That condition arises as a direct strengthening of
($\mc{H}$). In order to quantify the separation required in $(\mc{H})$, we define
\begin{equation}\label{eq:tau0}
\tau_0 =  \left \{ \max \tau: \;\exists w \in \R^n,\; \nnorm{w}_2 \leq 1 \;
  \; \;
\text{s.t.} \; \; \; \frac{X^{\T} w}{\sqrt{n}} \gec \tau \bm{1} \right \}.   
\end{equation}
Note that $\tau_0 > 0$ if and only if $(\mc{H})$ is fulfilled. Also note that with $\nnorm{X_j}_2^2 = n \; \forall j$, it holds that $\tau_0
\leq 1$. Introducing the Gram matrix $\Sigma = \frac{1}{n} X^{\T} X$, we have
by convex duality that
\begin{equation}\label{eq:tau0sq}
\tau_0^2 = \min_{\lambda \in T^{p-1}} \frac{1}{n} \nnorm{X \lambda}_2^2 =
\min_{\lambda \in T^{p-1}} \lambda^{\T} \Sigma \lambda, \quad \;
\text{where} \; \; T^{p-1} = \{\lambda\in \R_+^p: \; \bm{1}^{\T} \lambda  = 1
\}, 
\end{equation}
i.e.~in geometrical terms, $\tau_0$ equals the distance of the convex hull of
the columns of $X$ (scaled by $1/\sqrt{n}$) to the origin. Using terminology
from support vector machine classification (e.g.~\cite{ScholkopfSmola2002}, Sec. 7.2), $\tau_0$ can be interpreted as margin of a maximum margin
hyperplane with normal vector $w$ separating the columns of $X$ from the
origin. As argued below, in case that $\tau_0$ scales as a constant, 
overfitting is curbed. This is e.g.~\emph{not} fulfilled for orthonormal design, where
$\tau_0 = 1/\sqrt{p}$ (cf.~Section \ref{sec:designs}). 
\begin{cond}\label{cond:1}
A design $X$ is said to have a \textbf{self-regularizing property} if there exists
a constant $\tau > 0$ so that with $\tau_0$ as defined in \eqref{eq:tau0}, it
holds that $\tau_0 \geq \tau > 0$.  
\end{cond}
The term 'self-regularization' expresses the fact that the design itself
automatically generates a regularizing term, as emphasized in the next
proposition and the comments that follow. We point out that Proposition
\ref{prop:selfreg} is a qualitative statement preliminary to the main result
of the section (Theorem \ref{theo:prediction}) and mainly serves an
illustrative purpose.  
\begin{prop}\label{prop:selfreg} Consider the linear model
  \eqref{eq:linearmodel} with $\beta^* = 0$ and $y = \eps$ having entries that are
  independent random variables with zero mean and finite variance.
Suppose that $X$ satisfies Condition \ref{cond:1}. We then have
\begin{equation}\label{eq:selfregobj}
\min_{\beta \gec 0} \frac{1}{n} \nnorm{\eps - X \beta}_2^2 = \min_{\beta \gec 0}
 \frac{1}{n} \nnorm{\eps - \wt{X} \beta}_2^2 + \tau^2
    (\bm{1}^{\T} \beta)^2  + O_{\p}\left(\frac{1}{\sqrt{n}} \right),
\end{equation}
with $\wt{X} = (\Pi X) D$, where $\Pi$ is a projection onto an
$(n-1)$-dimensional subspace of $\R^n$ and $D$ is a diagonal matrix, the
diagonal entries being contained in $[\tau, 1]$.
Moreover, if $\frac{1}{n} \nnorm{X^{\T} \eps}_{\infty} = o_{\p}(1)$, then any
minimizer $\wh{\beta}$ of \eqref{eq:selfregobj} obeys $\frac{1}{n}
\nnorm{X \wh{\beta}}_2^2 = o_{\p}(1)$. 
\end{prop} 
In Proposition \ref{prop:selfreg}, the pure noise fitting problem is decomposed into
a fitting term with modified design matrix $\wt{X}$, a second term that can be interpreted as \emph{squared} non-negative lasso
penalty $\tau^2 (\bm{1}^{\T} \beta)^2$  (cf.~\eqref{eq:nnlasso}) and an additional stochastic term of
lower order. As made precise in the proof, the lower bound on $\tau$ implies that 
the $\ell_1$-norm of any minimizer is upper bounded by a constant. Prevention
of overfitting is then an immediate consequence under
the further assumption that the term  $\frac{1}{n} \nnorm{X^{\T}
  \eps}_{\infty} = o_{\p}(1)$ tends to zero. This holds under rather mild additional conditions
on $X$ \cite{Lounici2008} or more stringent conditions on the tails of the
noise distribution. As a last comment, let us make the connection of
the r.h.s.~of \eqref{eq:selfregobj} to a non-negative lasso problem more
explicit. Due to the correspondence of the level sets of the mappings $\beta \mapsto \bm{1}^{\T}
\beta$ and $\beta \mapsto (\bm{1}^{\T} \beta)^2$ on $\R_+^p$, we have
\begin{equation}\label{eq:squarednnlasso}
\min_{\beta \gec 0} \frac{1}{n} \nnorm{\eps - \wt{X} \beta}_2^2 + \tau^2
(\bm{1}^{\T} \beta)^2 = \min_{\beta \gec 0} \frac{1}{n} \nnorm{\eps - \wt{X} \beta}_2^2 + \gamma(\tau) \bm{1}^{\T} \beta,
\end{equation}
where $\gamma$ is a non-negative, monotonically increasing function of
$\tau$. Proposition \ref{prop:selfreg} in conjunction with
\eqref{eq:squarednnlasso} provides a high-level understanding of what will be shown in the sequel, namely that NNLS may inherit desirable properties of the
(non-negative) lasso with regard to prediction, estimation and sparsity of the solution.

\subsection{Slow rate bound}\label{sec:designsselfreg}

Condition \ref{cond:1} gives rise to the following general bound on the
$\ell_2$-prediction error of NNLS. Note that in Theorem \ref{theo:prediction}
below, it is not assumed that the linear model is specified
correctly. Instead,  we only assume that there is a
fixed target $f = (f_1,\ldots,f_n)^{\T}$ to be approximated by a non-negative
combination of the columns of $X$.
\begin{theo}\label{theo:prediction} 
Let $y = f + \eps$, where $f \in \R^n$ is fixed and $\eps$ has
i.i.d.~zero-mean sub-Gaussian entries with parameter $\sigma$ \footnote{See
  Appendix \ref{app:subgaussian} for background on sub-Gaussian random variables}. 
Define 
\begin{equation*}
\mc{E}^* = \min_{\beta \gec 0} \frac{1}{n} \nnorm{X \beta - f}_2^2, \qquad
\wh{\mc{E}} = \frac{1}{n} \nnorm{X \wh{\beta} - f}_2^2.
\end{equation*}
Suppose that $X$ satisfies Condition \ref{cond:1}. Then, for any minimizer 
$\wh{\beta}$ of the NNLS problem \eqref{eq:nnls} and any $M \geq 0$, it holds
with probability no less than $1 - 2 p^{-M^2}$ that
\begin{equation}\label{eq:theoprediction}
\wh{\mc{E}} \leq \mc{E}^* +  \left( \frac{6 \nnorm{\beta^{\ast}}_1
     + 8 \sqrt{\mc{E^{*}}}}{\tau^2}   \right) (1 + M) \sigma \sqrt{\frac{2 \log
    p}{n}} + \frac{16 (1 + M)^2 \sigma^2}{\tau^2} \frac{ \log p}{n},
\end{equation}
for all $\beta^{\ast} \in \argmin_{\beta \gec 0} \frac{1}{n} \nnorm{X \beta - f}_2^2$. 
\end{theo}

\paragraph{Comparison with the slow rate bound of the lasso.}

Theorem \ref{theo:prediction} bounds the excess error by a term of order
$O(\norm{\beta^*}_1 \sqrt{\log(p)/n})$, which implies that NNLS can be
consistent in a regime in which the number of predictors $p$ is nearly
exponential in the number of observations $n$. That is, NNLS constitutes a 'persistent procedure' in the 
spirit of Greenshtein and Ritov \cite{GreenshteinRitov2004} who coined the
notion of 'persistence' as distinguished from classical consistency with a
fixed number of predictors. The excess error bound of Theorem
\ref{theo:prediction} is of the same order of magnitude as the corresponding
bound of the lasso \cite{GreenshteinRitov2004, Bartlett2012, Buhlmann2011}  
that is typically referred to as slow rate bound. Since the bound \eqref{eq:theoprediction} depends on $\tau$,
it is recommended to solve the quadratic program in \eqref{eq:tau0sq} before
applying NNLS, which is roughly of the same computational cost. Unlike Theorem
\ref{theo:prediction}, the slow rate bound of the lasso does not require any conditions on the
design and is more favourable than \eqref{eq:theoprediction} regarding the constants. 
In \cite{Lederer2012, HebiriLederer2012}, improvements of the slow rate bound
are derived. On the other hand, the results for the lasso require the regularization parameter to
be chosen appropriately.      

\begin{remarque}
NNLS has been introduced as a tool for 'non-negative data'. In this
context, the assumption of zero-mean noise in Theorem \ref{theo:prediction} is questionable. In case
that the entries of $\eps$ have a positive mean, one can decompose
$\eps$ into a constant term, which can be absorbed into the linear
model, and a second term which has mean zero, so that Theorem
\ref{theo:prediction} continues to be applicable. 
\end{remarque}

\section{Fast rate bound for prediction and bounds on the $\ell_q$-error for estimation, $1
  \leq q  \leq 2$}\label{subsec:oracleprediction}

Within this section, we further elaborate on the similarity in performance of 
$\ell_1$-regularization and NNLS for designs with a self-regularizing
property. We show that the latter admits a reduction to the scheme pioneered 
in \cite{Bickel2009} to establish near-optimal performance guarantees of the
lasso and the related Dantzig selector \cite{CandesTao2007} with respect to estimation of
$\beta^*$ and and prediction under a sparsity scenario. Similar results are
shown e.g.~in \cite{Bunea2007, CandesTao2007, Geer2008, ZhangHuang2008, MeinshausenYu2009,
  CandesPlan2009, Zhang2009b}, and we shall state results of that flavour
for NNLS below. Throughout the rest of the paper, the data-generating model
\eqref{eq:linearmodel} is considered for a sparse target $\beta^*$
with support $S = \{j:\beta_j^{*} > 0 \}, \; 1 \leq |S| = s < n$. 



\paragraph{Reduction to the scheme used for the lasso.}

As stated in the next lemma, if the design satisfies Condition
\ref{cond:1}, the NNLS estimator $\wh{\beta}$ has, with high
probability, the property that $\wh{\delta} = \wh{\beta} - \beta^*$
has small $\ell_1$-norm, or that $\wh{\delta}$ is contained in the convex cone
\begin{equation}\label{eq:lassocone}
\{\delta \in \R^p:\, \nnorm{\delta_{S^c}}_1 \leq c_0 \nnorm{\delta_S}_1 \}, \; \quad
\text{where} \; \; c_0 = \frac{3}{\tau^2}. 
\end{equation}
The latter property is shared by the lasso and Dantzig selector with different
values of the constant $c_0$ \cite{CandesTao2007, Bickel2009}. 
\begin{lemma}\label{lem:lassocone}  Assume that $y = X \beta^{\ast} + \eps$, where $\beta^{\ast} \gec
0$ has support $S$, $\eps$ has i.i.d.~zero-mean sub-Gaussian entries with
parameter $\sigma$. Further assume that $X$ satisfies Condition
\ref{cond:1}. Denote $\wh{\delta} = \wh{\beta} -
\beta^*$. Then, for any $M \geq 0$, at least one of the following two events occurs with probability
no less than $1 - 2 p^{-M^2}$:  
\begin{align*}
\left\{\nnorm{\wh{\delta}_{S^c}}_1 \leq \frac{3}{\tau^2} \nnorm{\wh{\delta}_S}_1
 \right \}, \quad \text{and} \quad \left\{ \nnorm{\wh{\delta}}_1 \leq 4  (1 +
   M) \left(1 + \frac{3}{\tau^2} \right) \sigma \sqrt{\frac{2 \log p}{n}} \; \; \right \}.
\end{align*}
\end{lemma}
The rightmost event is most favourable, since it immediately yields the assertion of
Theorem \ref{theo:oracle} below.  On the other hand, under the leftmost event,
one is in position to carry over techniques used for analyzing the lasso. When
combined with the following Condition \ref{cond:2}, near-optimal rates with
regard to estimation and prediction can be obtained. 
\begin{cond}\label{cond:2}
Let $\mc{J}(s) = \{ J \subseteq \{1,\ldots,p\}: 1 \leq |J| \leq s \}$ and for $J
\in \mc{J}(s)$ and $\alpha \geq 1$,
\begin{equation*}
  \mc{R}(J, \alpha) = \{\delta \in \R^p:\, \nnorm{\delta_{J^c}}_1 \leq  \alpha
\nnorm{\delta_J}_1 \}. 
\end{equation*}
We say that the design satisfies the $(\alpha, s)$-\textbf{restricted
  eigenvalue condition} if there exists a constant $\phi(\alpha, s)$ so
that 
\begin{equation}\label{eq:recondition}
\min_{J \in \mc{J}(s)} \, \min_{\delta \in \mc{R}(J, \alpha) \setminus
  0} \, \, \frac{\delta^{\T} \Sigma \delta}{\nnorm{\delta_J}_2^2} \geq
\phi(\alpha, s) > 0. 
\end{equation}   
\end{cond}
Condition \ref{cond:2} is introduced in \cite{Bickel2009}. It is weaker than several
other conditions such as those in \cite{CandesTao2007, MeinshausenYu2009}
employed in the same context; for a comprehensive comparison of these and
related conditions, we refer to \cite{Geer2008}. Moreover, we note that
Condition \ref{cond:2} is satisfied with overwhelming probability if $X$
belongs to a rather broad class of random sub-Gaussian matrices with
independent rows as long as $n$ scales as $\Omega(s \log p)$ \cite{RudelsonZhou2012}.\\ 
Using Lemma \ref{lem:lassocone}, the next statement follows along the lines of the analysis in \cite{Bickel2009}.
\begin{theo} \label{theo:oracle} 
In addition to the conditions of Lemma \ref{lem:lassocone},
assume further that $X$ satisfies the $(3/\tau^2, s)$-restricted eigenvalue condition. It then
holds for any $q \in [1,2]$ and any $M \geq 0$ that 
\begin{align*}\label{eq:optimalrate}   
& \nnorm{\wh{\beta} -
  \beta^{\ast}}_q^q \leq \frac{2^{3q - 2}}{\left\{ \phi(3/\tau^2, s)
  \right\}^q} \left(1 + \frac{3}{\tau^2} \right)^{2q} \, s \; \left((1 + M)
  \sigma \sqrt{\frac{2 \log p}{n}} \right)^q\\
&\frac{1}{n} \nnorm{X \wh{\beta}
  - X \beta^{\ast}}_2^2 \leq \frac{8 (1 + M)^2 \sigma^2}{\phi(3/\tau^2, s)} \left(1 +
  \frac{3}{\tau^2} \right)^2  \; \frac{s \; \log p}{n} ,
\end{align*}
with probability no less than $1 - 2 p^{-M^2}$. 
\end{theo}
Theorem \ref{theo:oracle} parallels Theorem 7.2 in \cite{Bickel2009},
establishing that the lasso adapts to the underlying sparsity, since its 
performance attains (apart from factors logarithmic in $p$) what could be achieved if the support $\beta^*$ were
known in advance. The rates of Theorem \ref{theo:oracle} for NNLS are the same
as those for the lasso, modulo (less favourable) multiplicative constants.\\ 
The required condition on the design is a combination of the self-regularizing
property and the restricted eigenvalue condition. At first glance, these two
conditions might appear to be contradicting each other, since the first one is
not satisfied if the off-diagonal entries of $\Sigma$ are too small, while
for $\alpha \geq 1$, we have $\phi(\alpha, s) \leq 2 (1 - \max_{j,k, \, j \neq k}
\scp{X_j}{X_k}/n)$. We resolve this apparent contradiction in Section \ref{sec:designs} by
providing designs satisfying both conditions simultaneously. The use 
of Condition \ref{cond:2} in place of more restrictive conditions like
restricted isometry properties (RIP, \cite{CandesTao2007}) used earlier
in the literature turns out to be crucial here, since these conditions
impose much stronger constraints on the magnitude of the off-diagonals entries of
$\Sigma$ as discussed in detail in \cite{Raskutti2010}.   
\\
A result of the same spirit as Theorem \ref{theo:oracle} is shown in the recent paper
\cite{Meinshausen2013} by Meinshausen who has independently studied the
performance of NNLS for high-dimensional linear models. That paper provides
an $\ell_1$-bound for estimation of $\beta^*$ and a fast rate bound for
prediction with better constants than those in above theorem, even though the
required conditions are partly more restrictive. The ingredients leading to
those bounds are the self-regularizing property, which is termed 'positive
eigenvalue condition' there, and the 'compatibility condition' \cite{Geer2008, Geer2009}  
which is used in place of Condition \ref{cond:2}. We prefer the latter here,
because the 'compatibility condition' is not sufficient to establish $\ell_q$-bounds for
estimation for $q > 1$. As distinguished from our Theorem \ref{theo:oracle}, a lower bound
on the minimum non-zero coefficient of $\beta^{\ast}$ is additionally required
in the corresponding result in \cite{Meinshausen2013}.

\section{Estimation error with respect to the $\ell_{\infty}$ norm and support
recovery by thresholding}\label{sec:sparserecovery} 

In the present section, we directly derive bounds on the
$\ell_{\infty}$-estimation error of NNLS without resorting to 
techniques and assumptions used in the analysis of the lasso. Instead, we
build on the geometry underlying the analysis of NNLS for sparse recovery in
the noiseless case
\cite{DonohoTanner2005,DonohoTanner2005b,WangTang2009,DonohoTanner2010}. In
light of the stated bounds, we subsequently study the performance of a thresholded 
NNLS estimator with regard to support recovery.  

\subsection{Main components of the analysis}

In the sequel, we provide the main steps towards the results stated in this section. Basic to our approach is a decomposition of the NNLS
problem into two sub-problems corresponding to the support $S$ and the off-support $S^c$. For this purpose, we need to introduce the following
quantities. For a given support $S$, let $\Pi_S$ and $\Pi_S^{\perp}$ denote
the projections on the subspace spanned by $\{ X_j \}_{j \in S}$ and its
orthogonal complement, respectively. In the context of the linear model
\eqref{eq:linearmodel}, we then set 
\begin{equation}\label{eq:offsupportquantities}  
Z = \Pi_S^{\perp} X_{S^c} \quad \text{and} \quad \xi = \Pi_S^{\perp} \eps. 
\end{equation}
These quantities appear in the following key lemma that contains the aforementioned decomposition of the 
NNLS problem. 
\begin{lemma}\label{lem:decoupling} Let $Z$ and $\xi$ be as defined in
  \eqref{eq:offsupportquantities}. Consider the two non-negative least squares problems
\begin{align*}
&(P1): \min_{\beta^{(P1)} \gec 0} \frac{1}{n} \nnorm{\xi  - Z
  \beta^{(P1)}}_2^2,  \\
&(P2): \; \min_{\beta^{(P2)} \gec 0} \frac{1}{n} \nnorm{\Pi_S \eps + X_S \beta_S^{*} - X_S \beta^{(P2)}  - \Pi_S X_{S^c} \wh{\beta}^{(P1)}}_2^2
\end{align*}
with minimizers $\wh{\beta}^{(P1)}$ of $(P1)$ and
$\wh{\beta}^{(P2)}$ of $(P2)$, respectively. If $\wh{\beta}^{(P2)} \succ 0$,
then setting $\wh{\beta}_S = \wh{\beta}^{(P2)}$ and $\wh{\beta}_{S^c} =
\wh{\beta}^{(P1)}$ yields a minimizer $\wh{\beta}$ of the
non-negative least squares problem \eqref{eq:nnls}.
\end{lemma}
Lemma \ref{lem:decoupling} is used in the proof of Theorem
\ref{theo:robustsparserecovery} below in the following way.
We first study the off-support problem $(P1)$ separately, establishing an upper bound 
on the $\ell_1$-norm of its minimizer $\wh{\beta}^{(P1)}$ in dependence of
the \emph{separating hyperplane constant} introduced in the next paragraph. Substituting that bound
into $(P2)$, we conclude an upper bound on $\nnorm{\beta_S^* -
  \wh{\beta}^{(P2)}}_{\infty}$ and in turn, by the lemma, on
$\nnorm{\wh{\beta} - \beta^*}_{\infty}$. In this second step, we exploit the
fact that if $\wh{\beta}^{(P2)} \succ 0$, it equals the corresponding
unconstrained least squares estimator. 
\paragraph{Separating hyperplane constant.} To establish an upper bound on
$\nnorm{\wh{\beta}^{(P1)}}_1$, we require a positive lower bound on the following
quantity to which we refer as separating hyperplane constant, which is
nothing else than the constant \eqref{eq:tau0} introduced in the context of
self-regularization designs in Section \ref{sec:nnlsdoesnotoverfit} with respect to
the matrix $Z$ in \eqref{eq:offsupportquantities}. The term 'separating
hyperplane constant' follows the geometric interpretation as margin of a
hyperplane that contains the origin and that separates $\{X_{j} \}_{j \in S}$
from $\{ X_j \}_{j \in S^c}$. Accordingly, for given $S$, we define
\begin{equation}\label{eq:tau}
\tau(S) = \left \{ \max \tau: \;\exists w \in \R^n,\; \nnorm{w}_2 \leq 1 \; \; \;
\text{s.t.} \; \; 
\frac{1}{\sqrt{n}} X_S^{\T}w = 0 \; \; \text{and} \; \; \frac{1}{\sqrt{n}} X_{S^c}^{\T}{w} \gec \tau \bm{1} \right \}.       
\end{equation}
By convex duality, we have 
\begin{align}\label{eq:taudual}
\begin{split}
  \tau^2(S) &= \hspace*{-0.3cm}\min_{\substack{\theta \in \R^s \\ \;\;\;\;\;\;\;\lambda \in T^{p-s-1}}} \frac{1}{n} \norm{X_S \theta
  - X_{S^c} \lambda}_2^2, \; \, \text{where} \; \, T^{p-s-1} = \{\lambda \in
\R_{+}^{p-s}: \; \lambda^{\T} \bm{1} = 1 \} \\
&=  \min_{\lambda \in T^{p-s-1}}  \lambda^{\T} \, \frac{1}{n} X_{S^c}^{\T}
\Pi_S^{\perp} X_{S^c} \, \lambda =  \min_{\lambda \in T^{p-s-1}} \lambda^{\T} \,
\frac{1}{n} Z^{\T} Z \, \lambda. 
\end{split}
\end{align}
The last line highlights the connection to \eqref{eq:tau0sq} in Section
\ref{sec:nnlsdoesnotoverfit}. Expanding $\frac{1}{n} Z^{\T} Z$ under the
assumption that the submatrix $\Sigma_{SS}$ is invertible, $\tau^2(S)$
can also be written as
\begin{equation}\label{eq:taucondcovariance}
\tau^2(S) = \min_{\lambda \in T^{p-s-1}} \lambda^{\T} \left(\Sigma_{S^cS^c}
  - \Sigma_{S^cS} \Sigma_{SS}^{-1} \Sigma_{SS^c}\right)\lambda
\end{equation}
It is shown in \cite{Slawski2012} that having $\tau(S) > 0$ is a necessary and
sufficient condition for recovery of $\beta^*$ by NNLS in the absence of noise
($\eps = 0$). Thus, the appearance of $\tau(S)$ in the present context is
natural.  

\subsection{Bounds on the $\ell_{\infty}$-error}

The upper bound of the next theorem additionally depends on the
quantities below, which also appear in the upper bound on the
$\ell_{\infty}$-error of the lasso \cite{Wain2009}.  
\begin{equation}\label{eq:constantsp1}
\beta_{\min}(S) = \min_{j \in S} \beta_j^{\ast}, \; \; \;  K(S) = \max_{\norm{v}_{\infty} = 1} \nnorm{\Sigma_{SS}^{-1} v}_{\infty},
\; \; \;  \phi_{\min}(S) = \min_{\norm{v}_2 = 1} \norm{\Sigma_{SS}  v}_2.  
\end{equation}
\begin{theo}\label{theo:robustsparserecovery}
Assume that $y = X \beta^{\ast} + \eps$, where $\beta^{\ast} \gec
0$ and $\eps$ has i.i.d. zero-mean sub-Gaussian entries with parameter $\sigma$. For $M \geq 0$, set 
\begin{equation}\label{eq:bounds_robustsparserecovery} 
b = \frac{2(1 + M) \sigma
  \sqrt{\frac{2 \log p}{n}}}{\tau^2(S)}, \quad \text{and} \quad \wt{b} = b \cdot K(S) +
\frac{(1 + M)
  \sigma}{\sqrt{\phi_{\min}(S)}} \sqrt{\frac{2 \log p}{n}}.   
\end{equation}
If $\beta_{\min}(S) > \wt{b}$, then the NNLS estimator $\wh{\beta}$ has the following properties with probability no less
than $1 - 4 p^{-M^2}$:
\begin{equation*}
\nnorm{\wh{\beta}_{S^c}}_1 \leq b \quad \text{and} \quad \nnorm{\wh{\beta}_S -
  \beta_S^*}_{\infty} \leq \wt{b}. 
\end{equation*}
\end{theo}
\paragraph{Discussion.}
Theorem \ref{theo:robustsparserecovery} can be summarized as follows. Given a
sufficient amount of separation between $\{X_j \}_{j \in S}$ and $\{X_j \}_{j
\in S^c}$ as quantified by $\tau^2(S)$, the $\ell_1$-norm of the off-support
coefficients is upper bounded by the effective noise level proportional to
$\sqrt{\log(p)/n}$  divided by $\tau^2(S)$, provided that the entries of
$\beta_S^*$ are all large enough. The upper bound $\wt{b}$ depends in particular on the ratio $K(S)/\tau^2(S)$. In Section 
\ref{sec:designs}, we discuss a rather special design for which $\tau^2(S) =
\Omega(1)$; for a broader class of designs that is shown to satisfy the
conditions of Theorem \ref{theo:oracle} as well, $\tau^2(S)$ roughly scales as
$\Omega(s^{-1})$. Moreover, we have $\{\phi_{\min}(S) \}^{-1} \leq K(S) \leq
\sqrt{s} \{\phi_{\min}(S) \}^{-1}$. In total, the $\ell_{\infty}$-bound can
hence be as large as $O(s^{3/2} \sqrt{\log(p)/n})$ even if $\tau^2(S)$
scales favourably, a bound that may already be implied by
the $\ell_2$-bound in Theorem \ref{theo:oracle}. On the positive side, Theorem
\ref{theo:robustsparserecovery} may yield a satisfactory result for $s$ constant or
growing only slowly with $n$, without requiring the restricted eigenvalue
condition of Theorem \ref{theo:oracle}.
\paragraph{Towards a possible improvement of Theorem \ref{theo:robustsparserecovery}.}  
The potentially sub-optimal dependence on the sparsity level
$s$ in the bounds of Theorem \ref{theo:robustsparserecovery} is too pessimistic relative
to the empirical behaviour of NNLS as discussed in Section
\ref{sec:empiricalresults}. The performance reported there can be better understood in light
of Theorem \ref{theo:ellinfsmin} below and the comments that follow. Our reasoning is based 
on the fact that any NNLS solution can be obtained from an ordinary least squares solution restricted to the variables
in the active set $F = \{j:\,\wh{\beta}_j > 0 \}$, cf.~Lemma
\ref{lem:NNLSopt} in Appendix \ref{app:lem:decoupling}. For the subsequent discussion to be meaningful, it is necessary that the NNLS solution
and thus its active set are unique, for which a sufficient condition is
thus established along the way.
\begin{theorem}\label{theo:ellinfsmin} Let the data-generating model be as in Theorem \ref{theo:robustsparserecovery} and
let $M \geq 0$ be arbitrary. If the columns of $X$ are in general 
linear position \eqref{eq:glp} and if 
\begin{equation}\label{eq:objectivepositivecond}
\frac{32 (1 + M)^2 \sigma^2}{\E[\eps_1^2]}
\; \, \frac{\log p}{\tau^2(S) \, n} \leq \left(1 - \frac{s}{n} \right),
\end{equation}
then, with probability at least $1 - \exp(-c (n-s)/\sigma^4) - 2 p^{-M^2}$,
the NNLS solution is unique and its active set $F = \{j:\,\wh{\beta}_j > 0 \}$ 
satisfies $|F| \leq \min\{n-1, p \}$. Conditional on that event, if
furthermore $\beta_ {\min}(S) > \wt{b}$ as defined in
\eqref{eq:bounds_robustsparserecovery}, then $S \subseteq F$ and 
\begin{equation}
\nnorm{\wh{\beta} - \beta^*}_{\infty} \leq \frac{(1 + M)
  \sigma}{\sqrt{\phi_{\min}(F)}} \sqrt{\frac{2 \log p}{n}},
\end{equation}
in with probability at least $1 - 6 p^{-M^2}$. 
\end{theorem}   
We first note that for $s / n$ bounded away from $1$, condition \eqref{eq:objectivepositivecond} is fulfilled if $n$ scales as
$\Omega(\log(p) / \tau^2(S))$. Second, the condition on
$\beta_{\min}(S)$ is the same as in the previous Theorem \ref{theo:robustsparserecovery}, so that the scope of 
application of the above theorem remains limited to designs with an
appropriate lower bound on $\tau^2(S)$. At the same time, Theorem
\ref{theo:ellinfsmin} may yield a significantly improved bound on
$\nnorm{\wh{\beta}- \beta^*}_{\infty}$  as compared to Theorem
\ref{theo:robustsparserecovery} if $\left\{ \phi_{\min}(F) \right \}^{-1/2}$, the smallest singular value of $X_F/\sqrt{n} \in \R^{n \times |F|}$,
scales more favourably than $K(S)/\tau^2(S)$, noting that as long as
$S \subseteq F$, $\left\{\phi_{\min}(S) \right \}^{-1/2} \leq \left\{ \phi_{\min}(F) \right
\}^{-1/2}$. In the first place, control of $\left\{ \phi_{\min}(F) \right
\}^{-1/2}$ requires control over the cardinality of the
set $F$. In a regime with $|F|$ scaling as a constant multiple of $s$ with $s = \alpha n$,
$0 \leq \alpha \ll 1$, it is not restrictive to assume that $\left\{
  \phi_{\min}(F) \right \}^{1/2}$ as the smallest singular value of a tall
submatrix of $X$ is lower bounded by a positive constant, as
it has been done in the literature on $\ell_1$-regularization
\cite{CandesTao2007, ZhangHuang2008, MeinshausenYu2009}. That assumption is
strongly supported by results in random matrix theory
\cite{Litvak2005,RudelsonVershynin2010}. In Section \ref{sec:designs}
the hypothesis of having $|F| \ll n$ is discussed in more detail for the class
of so-called equi-correlation-like designs. For equi-correlated
design, it is even possible to derive the distribution of $|F|$ conditional on
having $S \subseteq F$ (Proposition \ref{prop:sparsity_equicor} in Section \ref{sec:designs}).

\subsection{Support recovery by thresholding}

The bounds on the estimation error presented in the preceding two sections
imply that hard thresholding of the NNLS estimator may be an effective means
for recovery of the support $S$. Formally, for a threshold $t \geq 0$,  
the hard-thresholded NNLS estimator is defined by  
\begin{equation}\label{eq:nnlshard}
\wh{\beta}_j(t) = \begin{cases}
                        \wh{\beta}_j, \quad &\wh{\beta}_j > t, \\
                         0, \quad & \text{otherwise}, \; j=1,\ldots,p,
                       \end{cases}
\end{equation}
and we consider $\wh{S}(t) = \{j: \; \wh{\beta}_j > 0 \}$ as an estimator for
$S$. In principle, the threshold might be chosen according to Theorem
\ref{theo:robustsparserecovery} or \ref{theo:ellinfsmin}: if $t > b$ and
$\beta_{\min}(S) > b + \wt{b}$, where $b$ and $\wt{b}$ denote upper bounds on
$\nnorm{\wh{\beta}_{S^c}}_{\infty}$ and $\nnorm{\wh{\beta}_S -
  \beta_S^*}_{\infty}$, respectively, one has that $S = \wh{S}(t)$ with the
stated probabilities. This approach, however, is not practical, since the
bounds $b$ and $\wt{b}$ depend on constants that are not accessible. In the
sequel, we propose a data-driven approach as devised in
\cite{Genovese2012} for support recovery on the basis of marginal regression. A
central observation in \cite{Genovese2012} is that direct specification of the
threshold can be avoided if the purpose of thresholding is support
recovery. In fact, given a ranking $(r_j)_{j = 1}^p$ of the predictors $\{ X_j \}_{j = 1}^p$ so that $r_j
\leq s$ for all $j \in S$, it suffices to estimate $s$. In light of Theorems
\ref{theo:oracle} to \ref{theo:ellinfsmin}, NNLS may give
rise to such ranking by setting   
\begin{align}\label{eq:nnlsranks}
r_j = k \; \, \Longleftrightarrow \; \, \wh{\beta}_j = \wh{\beta}_{(k)}, \;
j=1,\ldots,p,  
\end{align}
where $\wh{\beta}_{(1)} \geq \wh{\beta}_{(2)} \geq \ldots \geq
\wh{\beta}_{(p)}$ is the sequence of coefficients arranged in decreasing
order. Theorem \ref{theo:Sdatadriven} below asserts that conditional on having
an ordering in which the first $s$ variables are those in $S$, support
recovery can be achieved by using the procedure in \cite{Genovese2012}. Unlike
the corresponding result in \cite{Genovese2012}, our statement is
non-asymptotic and comes with a condition that is easier to verify. We point
out that Theorem \ref{theo:Sdatadriven} is of independent interest, since it
is actually not specific to NNLS, but would equally apply to any estimator yielding the correct ordering of the
variables.  
\begin{theo}\label{theo:Sdatadriven} Consider the data-generating model of Theorem
  \ref{theo:robustsparserecovery} and suppose that the NNLS estimator has the
  property that according to \eqref{eq:nnlsranks}, it holds that $r_j \leq s$
  for all $j \in S$. For any $M \geq 0$, set 
\begin{align}\label{eq:shat}
\begin{split}
&\wh{s} = \max \left \{ 0 \leq k \leq (p-1): \delta(k) \geq (1 + M) \sigma \sqrt{2  
    \;  \log p}\right\} + 1,\\  
\text{where} \; \; &\delta(k) = \norm{(\Pi(k + 1) - \Pi(k)) y}_2, \; k=0,\ldots,(p-1), 
\end{split}
\end{align}  
with $\Pi(k)$ denoting the projection on the linear space spanned by the
variables whose ranks are no larger than $k$ (using $\Pi(0) = 0$). 
Let $\wh{S} = \{j:\,r_j \leq \wh{s} \}$.\\ If $\beta_{\min}(S) \geq 2 (1 + M)
\sigma \left\{ \phi_{\min}(S) \right \}^{-1/2}\sqrt{2 \log(p)/n}$, then
$\wh{S} = S$ with probability no less than $1 - 4 p^{-M^2}$. 
\end{theo}
We note that the required lower bound on $\beta_{\min}(S)$ is rather moderate. Similar
or even more stringent lower bounds are required throughout the literature on
support recovery in a noisy setup \cite{Wain2009, Lounici2008, CandesPlan2009,
  Zhang2009, Zhang2009b, Bunea2008}, and are typically already needed to ensure
that the variables in $S$ are ranked at the top (cf.~also Theorems \ref{theo:oracle} to \ref{theo:ellinfsmin}).\\   
Strictly speaking, the estimate $\wh{s}$ in Theorem \ref{theo:Sdatadriven} is not
operational, since knowledge of the noise level $\sigma$ is assumed. In
practice, $\sigma$ has to be replaced by a suitable estimator. Variance
estimation in high-dimensional linear regression with Gaussian errors continues to be a topic of
active research, with several significant advances made very recently
\cite{Huet2012}. In our experiments, this issue appears to be minor, because
even naive plug-in estimation of the form $\wh{\sigma}^2 = \frac{1}{n}
\nnorm{y - X \wh{\beta}}_2^2$ yields satisfactory results \footnote{We note that the denominator $n$ could
be replaced by $n - \nu$, with $\nu$ denoting the degrees of freedom of NNLS
(which, to the best of our knowledge, is not known).}(cf.~Section
\ref{sec:empiricalresults}).
A nice property of the approach is its computational
simplicity. Repeated evaluation of $\delta(k)$ in \eqref{eq:shat} can be
implemented efficiently by updating QR decompositions.
Finally, we note that subsequent to thresholding, it is beneficial to
re-compute the NNLS solution using data $(y, X_{\wh{S}})$ only, because the
removal of superfluous variables leads to a more accurate estimation of the
support coefficients.

\subsection{Comparison of NNLS and the non-negative lasso}\label{subsec:plasso} 

Let us recall the non-negative lasso problem \eqref{eq:nnlasso} given by
\begin{equation*}
\min_{\beta \gec 0} \frac{1}{n} \norm{y - X \beta}_2^2 + \lambda \bm{1}^{\T}
\beta, \quad \lambda > 0, 
\end{equation*}
with minimizer denoted by $\plasso$. In the present subsection, we elaborate
on similarities and differences of NNLS and the non-negative lasso regarding 
the $\ell_{\infty}$-error in estimating $\beta^*$. We first state a result
according to which the non-negative lasso succeeds in support recovery. We
then argue that in general, the non-negative lasso does not attain the optimal
rate $O(\sqrt{\log(p)/n})$ in estimation with respect to the
$\ell_{\infty}$-norm by providing a specific design as counterexample. We finally
conclude by summarizing benefits and drawbacks of the non-negative lasso and 
NNLS in an informal way.     

\paragraph{Non-negative irrepresentable condition.} We study the performance
of the non-negative lasso estimator $\plasso$ under a version of the 'irrepresentable condition' 
that takes into account non-negativity of the regression coefficients. The
irrepresentable condition has been employed e.g.~in \cite{Mei2006, ZhaoYu2006,
Wain2009, Geer2009} to study the (ordinary) lasso from the point of view of
support recovery. For given $S \subset \{1,\ldots, p\}$, the \emph{non-negative irrepresentable constant} is defined as 
\begin{equation}\label{eq:iotaS}
\iota(S) = \max_{j \in S^c} \Sigma_{jS} \Sigma_{SS}^{-1} \bm{1} = \max_{j \in
  S^c} X_j^{\T} X_S (X_S^{\T} X_S)^{-1} \bm{1}. 
\end{equation}
It follows from the analysis in \cite{Wain2009} that the non-negative
irrepresentable condition $\iota(S) < 1$ is necessary for the non-negative
lasso to recover the support $S$, cf.~Theorem \ref{theo:plasso} below.
\begin{remarque} We here point out that the condition $\iota(S) < 1$ can be regarded
as a strengthening of the condition 
\begin{equation}\label{eq:sephyperplane_plasso}
\exists w \in \R^n \; \, \text{s.t.} \; \, X_S^{\T} w = \bm{1} \; \;
\text{and} \; \, X_{S^c}^{\T} w \prec \bm{1},  
\end{equation}
which is a necessary condition for recovering a non-negative solution with
support $S$ as minimum $\ell_1$-norm solution of an underdetermined linear system
of equations \cite{Zhao2012}, which corresponds to a non-negative lasso
problem in the absence of noise. The non-negative
irrepresentable condition $\iota(S) < 1$ results from
\eqref{eq:sephyperplane_plasso} by additionally requiring the vector $w$ to be contained
in the column space of $X_S$. Condition \eqref{eq:sephyperplane_plasso}
highlights a conceptual connection to the separating hyperplane constant as
defined in \eqref{eq:tau}. Note that as distinguished from \eqref{eq:tau}, the
separating hyperplane underlying \eqref{eq:sephyperplane_plasso} does \textbf{not}
contain the origin.    
\end{remarque}     

\paragraph{Support recovery with the non-negative lasso.} Using the scheme developed
in \cite{Wain2009}, one can show that under the non-negative irrepresentable
condition and a suitable choice of the regularization parameter $\lambda$, the
non-negative lasso has the property that $\plasso_{S^c} = 0$, i.e.~no
false positive variables are selected, and support recovery can be deduced
from having an appropriate lower bound on the minimum support coefficient
$\beta_{\min}(S)$. Along with a bound on $\nnorm{\plasso - \beta^*}_{\infty}$,
this is stated in the next theorem.
\begin{theo}\label{theo:plasso}
Assume that $y = X \beta^{\ast} + \eps$, where $\beta^{\ast} \gec
0$ has support $S$ and $\eps$ has i.i.d.~zero-mean sub-Gaussian entries with
parameter $\sigma$. Suppose further that the non-negative
irrepresentable condition $\iota(S) < 1$ according to \eqref{eq:iotaS}
holds. For any $M \geq 0$, if
\begin{align}\label{eq:nnlassocondition}
\begin{split}
&\lambda > \frac{2 \lambda_M}{1 - \iota(S)}, \; \; \; \text{where} \; \; \; \lambda_M = (1 + M) \, \sigma \sqrt{\frac{2 \log p}{n}},\\
 \; \; 
\text{and} \; \; \;& \; \beta_{\min}(S) > b, \; \; \; \text{where} \; \; \; b =  \frac{\lambda}{2} \nnorm{\Sigma_{SS}^{-1} \bm{1}}_{\infty} + \frac{\lambda_M}{\sqrt{\phi_{\min}(S)}},
\end{split}
\end{align}
then $\{j:\plasso_j > 0 \} = S$  and $\nnorm{\plasso_S -   \beta_S^*}_{\infty} \leq b$
with probability at least $1 - 4 p^{-M^2}$.  
\end{theo}
There is some resemblance of the bound $b$ in \eqref{eq:nnlassocondition} and that of
Theorem \ref{theo:robustsparserecovery} for NNLS, with $\tau^2(S)$ playing a
role comparable to $1 - \iota(S)$ and $\nnorm{\Sigma_{SS}^{-1} \bm{1}}_{\infty}$ being a
lower bound on the quantity $K(S)$ defined in \eqref{eq:constantsp1}. On the
other hand, Theorem \ref{theo:plasso} yields a considerably stronger control
of the off-support coefficients ($\plasso_{S^c} = 0$) as does Theorem    
\ref{theo:robustsparserecovery}, which only provides an $\ell_1$-bound on
$\wh{\beta}_{S^c}$. Irrepresentable conditions as in above theorem are regarded
as rather restrictive in the literature \cite{ZhangHuang2008, MeinshausenYu2009, Zhang2009b}. Even in case the condition $\iota(S) < 1$ is
fulfilled, the choice of $\lambda$ in \eqref{eq:nnlassocondition} with
$\iota(S)$ possibly close to one may impose a rather stringent lower bound on
$\beta_{\min}(S)$ in order to achieve support recovery. At the same time, the choice $\lambda = 2 \sigma \sqrt{2 \log(p)/n}$ in combination with the restricted
eigenvalue condition (Condition \ref{cond:2}), which is regarded as far less 
restrictive than the irrepresentable condition, only yields a bound on
$\nnorm{\plasso - \beta^*}_q$ for $q \in [1,2]$, and it is no longer
guaranteed that $\plasso_{S^c} = 0$. As a result, two-stage
procedures like subsequent thresholding of $\plasso$ may be needed for support
recovery. However, this approach in general entails a sub-optimal condition on
$\beta_{\min}(S) = \Omega(\sqrt{s \log(p)/n})$ because of the term
$\nnorm{\Sigma_{SS}^{-1} \bm{1}}_{\infty}$ scaling as $\Theta(\sqrt{s})$ in the
worst case. In Appendix \ref{app:examplesup} this issue is illustrated by providing 
an explicit example of a design representing that worst case. We point out
that optimal rates for estimation in sup-norm of the order
$O(\sqrt{\log(p)/n})$ have been established e.g.~in \cite{Lounici2008,
  Bunea2008, CandesPlan2009} for the lasso under 'mutual
incoherence' conditions requiring fairly restrictive upper bounds
on the maximum inner product between two distinct columns of $X$.


\paragraph{NNLS vs. the non-negative lasso: pros and cons.}

To sum up, we list advantages and disadvantages of NNLS and the non-negative lasso,
thereby providing some guidance on when to use which approach in
practice. While NNLS can formally be seen as a special case of the
non-negative lasso with $\lambda = 0$, we suppose for the subsequent
discussion that $\lambda \geq \sigma\sqrt{2 \log(p)/n} $ as it is standard in the literature on the lasso.     

\begin{itemize}
\item As already stressed previously, reasonable performance of NNLS  
      in a high-dimensional regime is restricted to a specific class of
      designs, which excludes standard models such as random Gaussian design
      matrices (cf.~the discussion in Section \ref{sec:designs}). This 
      contrasts with the non-negative lasso, which has at least moderate performance  
      guarantees via a slow rate bound in spirit of Theorem \ref{theo:prediction}       
      \emph{without} further conditions on the design. 
\item As discussed in the previous paragraph, the non-negative lasso does
      not always attain the optimal rate for estimating $\beta^*$ in the sup-norm,
      in which case some room for improvement is left for NNLS. In Section  
      \ref{sec:empiricalresults}, we present two designs for which NNLS empirically  
      yields a better performance both with regard to estimation and support 
      recovery via thresholding, where the $\ell_{\infty}$-error in estimation enters 
      crucially. On the other hand, as asserted by Theorem \ref{theo:plasso} the non-negative lasso  
      succeeds in support recovery even without thresholding in certain
      regimes.  
\item From the point of view of a practitioner who is little familiar with
      theoretical results on how to set the regularization parameter of the
      (non-negative) lasso, NNLS has the advantage that it can be applied
      directly, without the need to specify or tune a regularization
      parameter.   
\end{itemize} 

\section{Discussion of the analysis of NNLS for selected designs}\label{sec:designs}

Our main results concerning the performance of NNLS as stated in Theorems
\ref{theo:prediction} to \ref{theo:ellinfsmin} are subject to the following
conditions: the self-regularizing property (Theorem \ref{theo:prediction}), a
combination of that property with a restricted eigenvalue condition (Theorem
\ref{theo:oracle}), a lower bound on the separating hyperplane constant
(Theorem \ref{theo:robustsparserecovery}), and sparsity of the NNLS solution 
(Theorem \ref{theo:ellinfsmin}). In the present section, we discuss to what
extent these conditions are fulfilled for selected designs, which we here
roughly divide into three classes. The first is the class of
\emph{non-self-regularizing designs} for which non-negativity constraints on
the regression coefficients do not seem to yield any significant
advantage. This is in contrast to the third class of \emph{equi-correlation-like
designs}, which are shown to be tailored to NNLS. The second class comprises
designs with a block or band structure arising in typical applications.

\subsection{Non-self regularizing designs}

In this paragraph, we provide several common examples of designs not having
the self-regularizing property of Condition \ref{cond:1}. Consequently, our
main results, which rely on that condition, do not apply. Those
designs can be identified by evaluating the quantity $\tau_0^2$
\eqref{eq:tau0sq} underlying Condition \ref{cond:1}. From   
\begin{equation}\label{eq:tau0sq_upper}  
\tau_0^2 = \min_{\lambda \in T^{p-1}} \lambda^{\T} \Sigma \lambda \leq \frac{1}{p^2}
\bm{1}^{\T} \Sigma \bm{1},
\end{equation}
we see that the sum of the entries of $\Sigma$ must scale as $\Omega(p^2)$ for
Condition \ref{cond:1} to be satisfied. In particular, this requires $\Sigma$
to have $\Omega(p^2)$ entries lower bounded by a positive constant, and 
a maximum eigenvalue scaling as $\Omega(p)$. Among others, this is not
fulfilled for the following examples.
\paragraph{Example 1: orthonormal design} \hfill\\
As already mentioned while motivating the self-regularizing property in
Section \ref{sec:nnlsdoesnotoverfit}, for $\Sigma = I$, $\tau_0^2$ attains the
upper bound in \eqref{eq:tau0sq_upper} which yields $\tau_0^2 = 1/p$.  

\paragraph{Example 2: power decay} \hfill\\
Let the entries of $\Sigma$ be given by $\sigma_{jk} = \rho^{|j-k|}$, $j,k=1,\ldots,p$ with $\rho
\in [0,1)$. From 
\begin{equation*}
\max_{1 \leq j \leq p} \sum_{k = 1}^p \sigma_{jk} \leq 2 \sum_{l = 0}^{p-1}
\rho^l \leq 2 (1 - \rho)^{-1}
\end{equation*}
and \eqref{eq:tau0sq_upper} it follows that $\tau_0^2 \leq 2 p^{-1} (1 - \rho)^{-1}$.
\paragraph{Example 3: random Gaussian matrices} \hfill\\
Consider a random matrix $X$ whose entries are drawn i.i.d.~from the standard
Gaussian distribution. We here refer to the work \cite{DonohoTanner2010} in
which NNLS is studied in the noiseless case. In that work, random Gaussian
matrices are considered as part of a broader class termed the centro-symmetric
ensemble. In a nutshell, that class encompasses random matrices with
independent columns whose entries have mean zero. The authors of
\cite{DonohoTanner2010} point out the importance
of Wendel's Theorem \cite{SchneiderWeil2008, Wendel1962}, which provides the
exact probability for the columns of $X$ being contained in a half-space,
i.e.~of having $\tau_0^2 > 0$. That result implies via Hoeffding's inequality that for $p/n > 2$,
$\tau_0^2 = 0$ with probability at least $1 - \exp(-n (p/n - 2)^2/2)$ so that 
the non-negativity constraints in NNLS become meaningless (cf.~the discussion
following \eqref{eq:halfspaceconstraint}).\hfill \vspace{8pt}\\
In all these three examples, a similar reasoning applies with regard to the
scaling of the separating hyperplane constant $\tau^2(S)$ \eqref{eq:taudual},
because its role is that of $\tau_0^2$ with respect to the matrix $Z =
\Pi_S^{\perp} X_{S^c}$. As a consequence, it is not hard to see that the scalings of the above examples
continue to hold (uniformly in $S$) with $p$ replaced by $p-s$.

\subsection{Designs with non-negative Gram matrices having a band or block structure}

\begin{center}
\begin{figure}[h!]%
\includegraphics[height = 0.15\textheight]{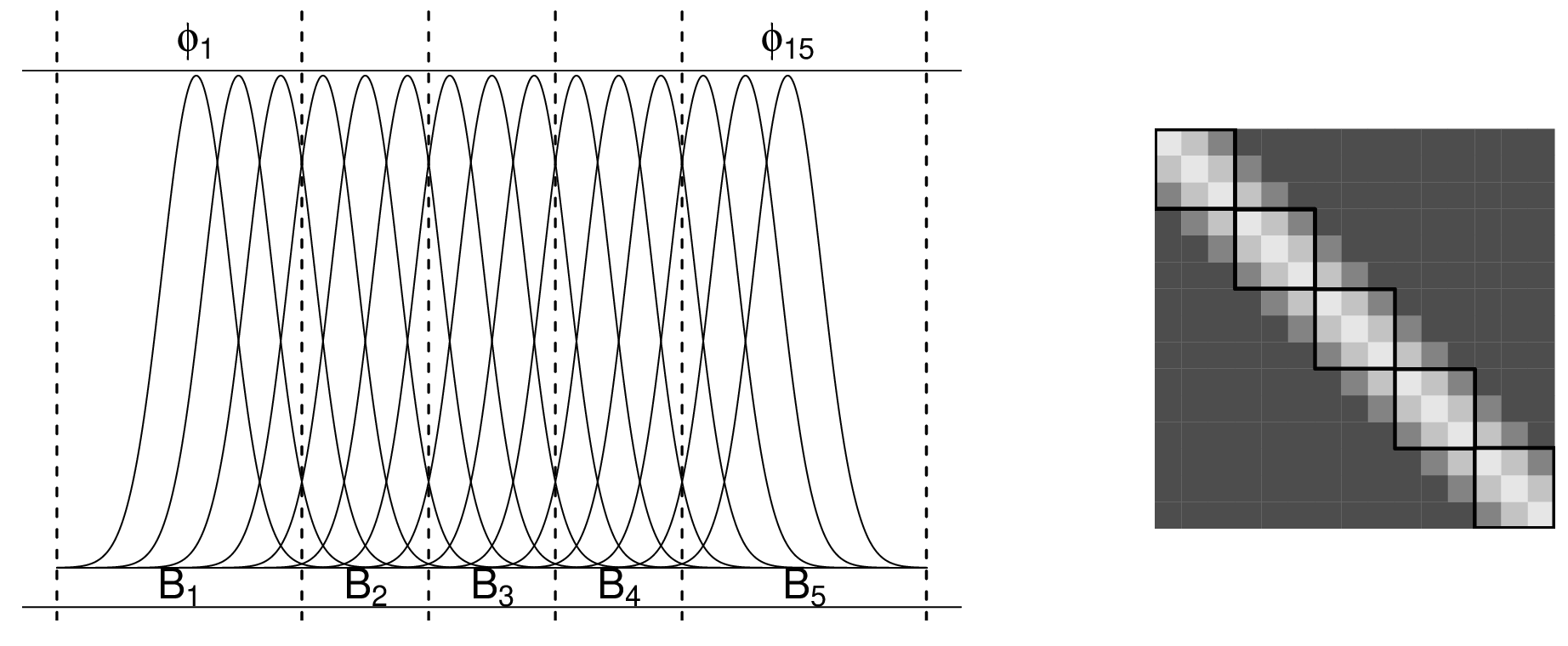}
\caption{Block partitioning of 15 Gaussians into $K = 5$ blocks. The right
  part shows the corresponding pattern of the Gram matrix.}
\label{fig:blockconstruction}
\end{figure}
\end{center}

We now present a simple sufficient condition for the self-regularizing
property to be satisfied, based on which we will identify a class of 
designs for which non-negativity of the regressions coefficients
may be a powerful constraint.\\
Suppose that the Gram matrix has the property that all its entries are lower
bounded by a positive constant $\sigma_0$. We then have the following lower
bound corresponding to the upper bound \eqref{eq:tau0sq_upper} above.   
\begin{equation}\label{eq:sufficient_lowerboundentry}
\tau_0^2 = \min_{\lambda \in T^{p-1}} \lambda^{\T} \Sigma \lambda \geq
\min_{\lambda \in T^{p-1}} \lambda^{\T} \left\{ \sigma_0 \bm{1} \bm{1}^{\T}
\right \} \lambda  = \sigma_0,
\end{equation}
i.e. Condition \ref{cond:1} is satisfied with $\tau^2 = \sigma_0$. More
generally, in case that $\Sigma$ has exclusively non-negative entries and the
set of variables $\{1, \ldots,p \}$ can be partitioned into blocks
$\{B_1,\ldots,B_K \}$ such that the minimum entries of the corresponding
principal submatrices of $\Sigma$ are lower bounded by a positive constant, 
then Condition \ref{cond:1} is satisfied with $\tau^2 = \sigma_0/K$:
\begin{align}\label{eq:sufficient_lowerboundentry_block}
\begin{split}
\tau_0^2 = \min_{\lambda \in T^{p-1}} \lambda^{\T} \Sigma \lambda &\geq
\min_{\lambda \in T^{p-1}} \sum_{l = 1}^K \lambda_{B_l}^{\T} \Sigma_{B_l \,
  B_l} \lambda_{B_l} \geq \sigma_0  \min_{\lambda \in T^{p-1}} \sum_{l = 1}^K (\lambda_{B_l}^{\T}
\bm{1})^2  = \sigma_0/K,  
\end{split}
\end{align}
where in the last equality we have used that the minimum of the map $x
\mapsto \sum_{l = 1}^K x_l^2$ over the simplex $T^{K-1}$ is attained for $x =
\bm{1}/K$.\\
As sketched in Figure \ref{fig:blockconstruction}, the lower bound \eqref{eq:sufficient_lowerboundentry_block}
particularly applies to design matrices whose entries contain the
function evaluations at points $\{u_i \}_{i=1}^n \subset [a, b]$ of
non-negative functions such as splines, Gaussian kernels and related 'localized' functions
traditionally used for data smoothing. If the points $\{u_i \}_{i =
  1}^n$  are placed evenly in $[a,b]$ then the corresponding
Gram matrix effectively has a band structure. For instance, suppose that $u_i = i/n, \;
i=1,\ldots,n$, and consider indicator functions of sub-intervals $\phi_j(u) = I\{u \in
[(\mu_j - h) \vee a, (\mu_j + h) \wedge b] \}$, where $\mu_j \in [0,1], \;
j=1,\ldots,p,$ and $h = 1/K$ for some positive integer $K$. Setting $X = (\phi_j(u_i))_{1 \leq i \leq n, \; 1 \leq j \leq p}$ and partitioning the $\{\mu_j \}$ by dividing $[0,1]$ into intervals
$[0,h], \; (h, 2h], \; \ldots, \; (1-h,1]$ and accordingly $B_l = \{j: \mu_j \in
((l-1) \cdot
h, l \cdot h] \}$, $l=1,\ldots,K$, we have that $\min_{1 \leq l \leq K}
\frac{1}{n} X_{B_l}^{\T} X_{B_l} \gec h$ such that Condition \ref{cond:1}
holds with $\tau^2 = h / K = 1/K^2$.
\emph{Applications.} As mentioned in the introduction, NNLS has been
shown to be remarkably effective in solving deconvolution problems
\cite{Speed2000, LinLeeSaul2004, SlawskiHein2010}. The
observations there are signal intensities measured over time, location etc.
that can be modelled as a series of spikes (Dirac impulses) convolved with a
\emph{point-spread function (PSF)} arising from a limited resolution of the
measurement device. The PSF is a non-negative localized function as 
outlined in the previous paragraph. Deconvolution of spike trains
is studied in more detail in Section \ref{subsec:sparsedeconvoultion}
below. Similarly, bivariate PSFs can be used to model blurring in 
greyscale images, and NNLS has been considered as a simple method for deblurring
and denoising \cite{Bardsley2006}.

\subsection{Equi-correlation-like designs}\label{eq:equicorrelationlike}

We first discuss equi-correlated design before studying random designs whose 
population Gram matrix has equi-correlation structure. While the population
setting is limited to having $n \geq p$, the case $n < p$ is possible for random
designs.  

\paragraph{Equi-correlated design.}
For $\rho \in (0,1)$, consider equi-correlated design with Gram matrix $\Sigma = (1 - \rho) I + \rho \bm{1} \bm{1}^{\T}$. We then have
\begin{equation}\label{eq:tau0_equicor}
\tau_0^2 = \min_{\lambda \in T^{p-1}} \lambda^{\T} \Sigma \lambda =\rho +
\min_{\lambda \in T^{p-1}} (1 - \rho) \norm{\lambda}_2^2 = \rho + \frac{1 -
  \rho}{p},
\end{equation}
so that the design has the self-regularizing property of Condition
\ref{cond:1}. Let $\emptyset \neq S \subset \{1,\ldots,p \}$ be
arbitrary. According to representation \eqref{eq:taucondcovariance}, the
corresponding separating hyperplane constant $\tau^2(S)$ can be evaluated
similarly to \eqref{eq:tau0_equicor}. We have
\begin{align}\label{eq:tauS_equicor}
\begin{split}
\tau^2(S) &=  \min_{\lambda \in T^{p-s-1}} \lambda^{\T} \left(\Sigma_{S^cS^c}
  - \Sigma_{S^cS} \Sigma_{SS}^{-1} \Sigma_{SS^c}\right)\lambda \\
&= \rho - \rho^2 \bm{1}^{\T} \Sigma_{SS}^{-1} \bm{1} + (1 - \rho)\min_{\lambda
  \in T^{p-s-1}} \nnorm{\lambda}_2^2  \\
&=\rho - \frac{s \rho^2}{1 + (s - 1) \rho} + \frac{1 - \rho}{p - s} =
\frac{\rho (1 - \rho)}{1 + (s - 1) \rho} + 
\frac{1 - \rho}{p - s} = \Omega(s^{-1}),
\end{split}
\end{align}
where from the second to the third line we have used that $\bm{1}$ is an
eigenvector of $\Sigma_{SS}$ corresponding to its largest eigenvalue $1 + (s -
1) \rho$. We observe that $\tau^2(S) = \tau^2(s)$,
i.e.~\eqref{eq:tauS_equicor} holds uniformly in $S$. We are not aware of any
design for which $\min_{S: \, |S|=s < p/2} \, \tau^2(S) \geq s^{-1}$, which
lets us hypothesize that the scaling of $\tau^2(S)$ in \eqref{eq:tauS_equicor}
uniformly over all sets of 
a fixed cardinality $s$ is optimal. On the other hand, when not requiring
uniformity in $S$, $\tau^2(S)$ can be as large as a constant independent of
$s$, as it is the case for the following example. Consider a Gram matrix of
the form
\begin{equation*}
\Sigma = \begin{pmatrix}  
\Sigma_{SS} & \Sigma_{SS^c} \\
\Sigma_{S^cS} & \Sigma_{S^cS^c} 
\end{pmatrix} = \begin{pmatrix}
\Sigma_{SS} & \mathbf{0} \\
\mathbf{0} & (1 - \rho) I + \rho \bm{1} \bm{1}^{\T}
\end{pmatrix} \quad \; \text{for} \; \; \rho \in (0,1). 
\end{equation*}
Combining \eqref{eq:tau0_equicor} and \eqref{eq:tauS_equicor}, we obtain that
$\tau^2(S) = \rho + \frac{1 - \rho}{p - s}$ independently of the specific form of  
$\Sigma_{SS}$. At the same time, this scaling does not hold uniformly over all
choices of $S$ with $|S| = s$ given the equi-correlation structure of the
block $\Sigma_{S^c S^c}$.

\paragraph{Sparsity of the NNLS solution for equi-correlated design.}
Exploiting the specifically simple structure of the Gram matrix, we are able
to derive the  distribution of the cardinality of the active set $F = \{j:\,
\wh{\beta}_j > 0 \}$ of the NNLS solution $\wh{\beta}$ conditional on the
event $\{ \wh{\beta}_S \succ 0 \}$. For the sake of better illustration, the result is 
stated under the assumption of Gaussian noise. Inspection of the proof shows
that, with appropriate modifications, the result remains valid for
arbitrary noise distributions.

\begin{prop}\label{prop:sparsity_equicor}
Consider the linear model $y = X \beta^* + \eps$, where 
$\beta^* \gec 0$, $\frac{1}{n} X^{\T} X = \Sigma = (1 - \rho) I + \rho \bm{1}
\bm{1}^{\T}$ for $\rho \in [0,1)$, and $\eps$ has i.i.d.~zero-mean, Gaussian
entries with variance $\sigma^2$. 
Let further $S = \{j:\, \beta_j^* > 0\}$. For any $M \geq 0$, if
$\beta_{\min}(S) > \frac{3 (1 + M) \sigma}{1 - \rho} 
\sqrt{2 \log(p)/n}$, then the event $\{\wh{\beta}_S \succ 0 \}$ occurs with
probability at least $1 - 4 p^{-M^2}$. Furthermore, let $z$ be a $(p-s)$-dimensional zero-mean
Gaussian random vector with covariance $(1 - \rho) I + \frac{\rho (1 -
  \rho)}{1 + (s -1) \rho} \bm{1} \bm{1}^{\T}$ and let $z_{(1)} \geq \ldots
\geq z_{(p-s)}$ denote the arrangement of the components of $z$ in decreasing order. Conditional  
on the event $\{ \wh{\beta}_S \succ 0 \}$, the cardinality of the active set $F =
\{j:\,\wh{\beta}_j > 0 \}$ has the following distribution: 
\begin{align}\label{eq:distn_activeset}
\begin{split}
&|F| \overset{\mc{D}}{=} s + I \left\{ z_{(1)} > 0 \right\} \left(1 + \max \left \{1 \leq j \leq p - s
  - 1:\, \zeta_j >  \theta(s, \rho) \right\} \right),\;\text{where} \\ 
&\zeta_j = \frac{z_{(j+1)}}{\sum_{k = 1}^j (z_{(k)} - 
  z_{(j+1)})}, \; \; j=1,\ldots,p-s-1, \; \; \text{and} \; \; \theta(s, \rho) =  \frac{\rho}{1 + (s - 1) \rho}.
\end{split}
\end{align}
\end{prop}

\begin{flushleft}
\begin{figure}[hb!]
\begin{tabular}{lll}
\hspace{-0.47cm}\includegraphics[height =
0.18\textheight]{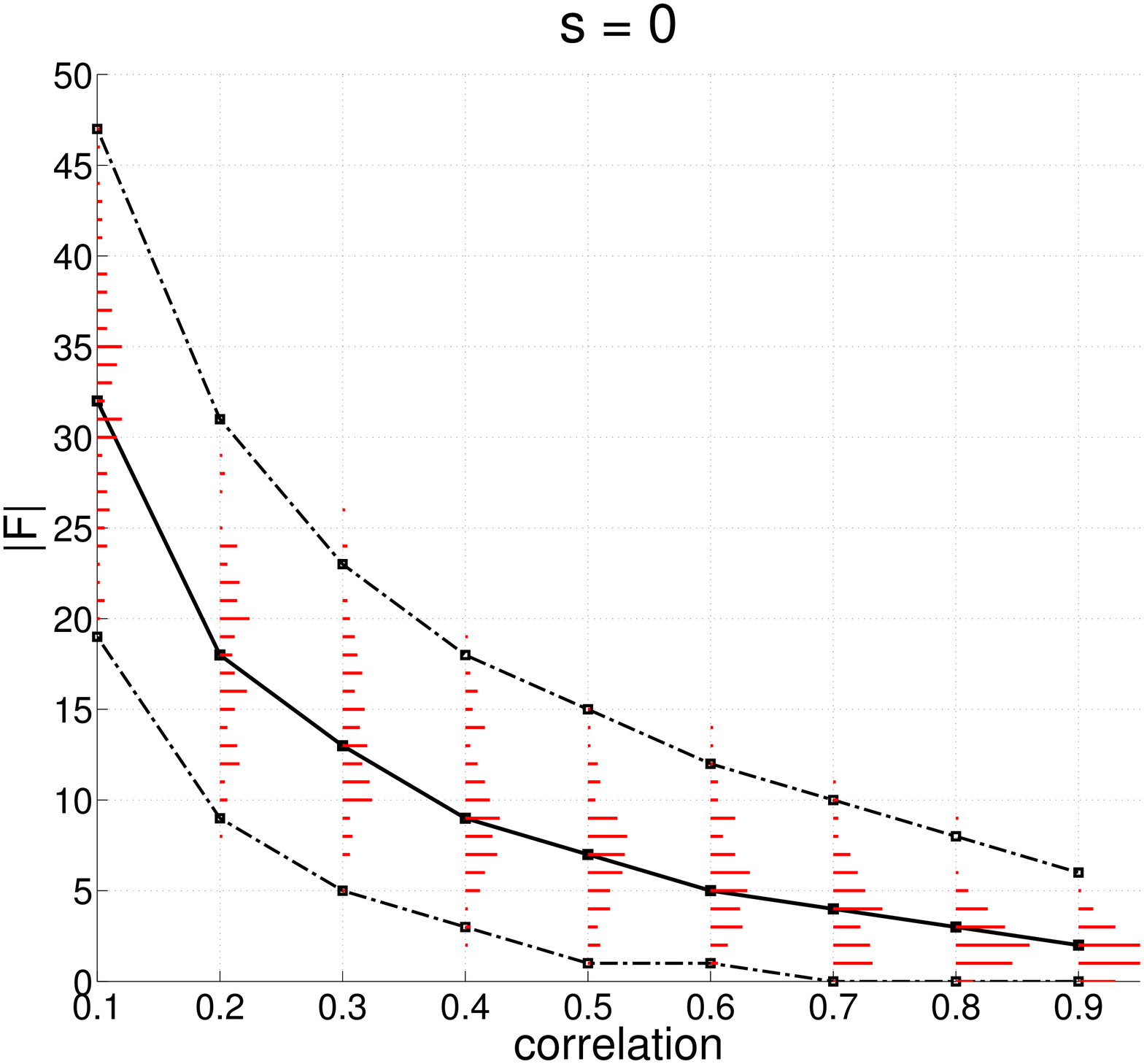} & \hspace{-0.51cm}
\includegraphics[height = 0.18\textheight]{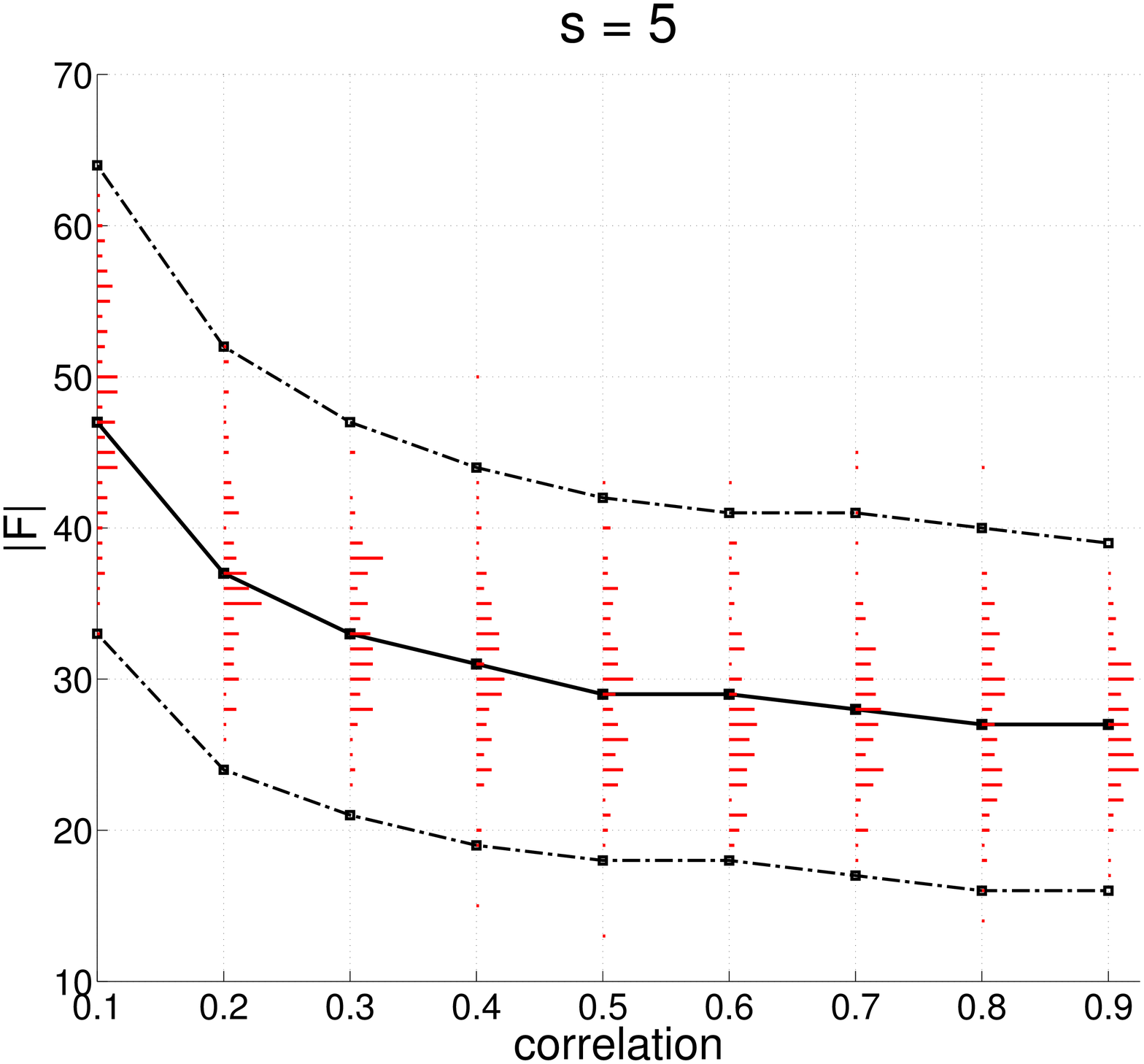} &  \hspace{-0.43cm}\includegraphics[height = 0.18\textheight]{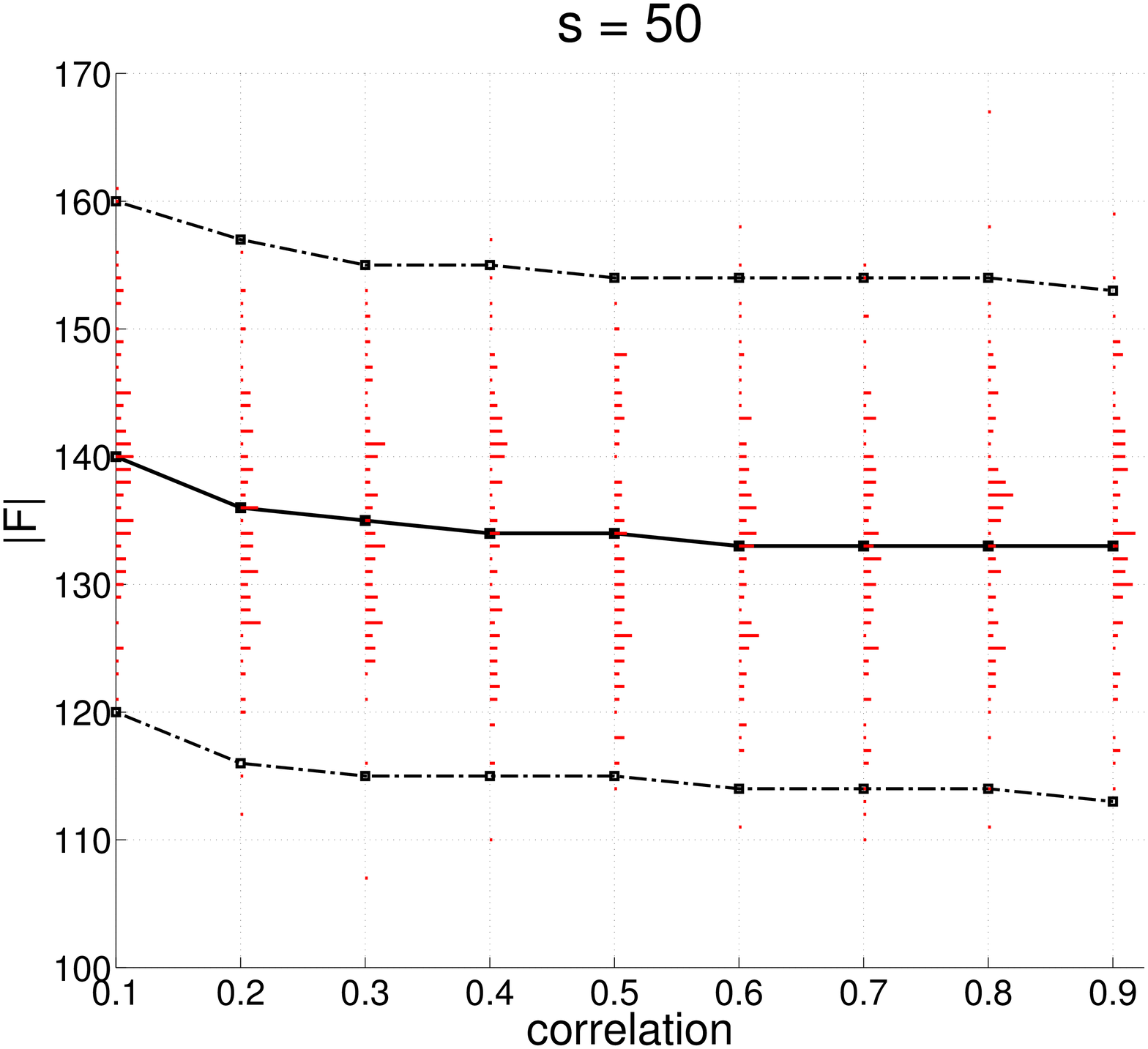}\\
\end{tabular}
\caption{Graphical illustration of Proposition \ref{prop:sparsity_equicor} for
  $p = 500$. The dotted lines represent the $\{0.01, 0.5, 0.99\}$-quantiles of the
distributions obtained from Proposition \ref{prop:sparsity_equicor} via Monte
Carlo simulation. The horizontal bars represent the corresponding relative
frequencies based on the solutions of $100$ random NNLS problems obtained for
each combination of $\rho$ and $s$.}\label{fig:cardF_equicor}
\end{figure}
\end{flushleft}

Proposition \ref{prop:sparsity_equicor} asserts that conditional on having the
support of $\beta^*$ included in the active set, the distribution of its
cardinality is $s$ plus an extra term, whose distribution depends on that of
the random variables $\{ \zeta\}_{j = 1}^{p-s-1}$ and a 'threshold' $\theta(s,
\rho)$. In order to better understand the role of these quantities, let us
first consider the case $\rho = 0$, i.e.~orthonormal design: since
$\theta(s, 0) = 0$, the distribution of $|F|$ is equal to $s$ plus the 
distribution of the number of non-negative components of a $(p-s)$-dimensional
Gaussian random vector, i.e.~a binomial distribution with $p-s$ trials and
a probability of success of $\frac{1}{2}$ (cf.~also Section \ref{sec:orthonormaldesign}).      
Once $\rho > 0$, the distribution of $|F|$ gets shifted towards $s$, noting
that $\{ \zeta\}_{j = 1}^{p-s-1}$ forms a non-increasing
sequence. Specifically, for $s = 0$, $\theta(0, \rho) = \frac{\rho}{1 -
  \rho}$, i.e.~the larger the correlation $\rho$, the stronger the
concentration of the distribution of $|F|$ near zero. The threshold
$\theta(s,\rho)$ is decreasing in $s$, i.e.~the number of extra variables 
increases with $s$. While the distribution of $\{ \zeta\}_{j = 1}^{p-s-1}$ is
not directly accessible, it can be approximated arbitrarily well by Monte Carlo
simulation for given $p$, $s$ and $\rho$ (note that the distribution does not
depend on the scale of the noise $\eps$). Figure \ref{fig:cardF_equicor} depicts the $0.01,
0.5, 0.99$-quantiles of the distribution of $|F|$ in virtue of
\eqref{eq:distn_activeset} for $p = 500$ and various choices of $\rho$ and
$s$. The results are based on $10,000$ Monte Carlo simulations for each value
of $s$. For comparison, for each pair $(s, \rho)$, we generate 100 datasets
$(X,y)$ with $n = p = 500$ according to the model of Proposition \ref{prop:sparsity_equicor}
with standard Gaussian noise (the components of $\beta_S^*$ are set to
the given lower bound on $\beta_{\min}(S)$ in to ensure that the event $\{\wh{\beta}_S \succ 0 \}$ has probability
close to one). We then solve the corresponding NNLS problems using the active
set algorithm of Lawson and Hanson \cite{LawsonHanson} and obtain the cardinalities of the active
sets. Figure \ref{fig:cardF_equicor} shows a strong agreement of the
predictions regarding the size of the active set based on the distribution of
Proposition \ref{prop:sparsity_equicor} and the empirical distributions.  

\paragraph{Non-negative random designs with equi-correlation structure.}
\hfill \\
We now consider random design matrices whose population Gram matrix is that
of equi-correlated design, but with the possibility that $n < p$. It is investigated to what extent these random design matrices inherit 
properties from the population setting studied in the previous paragraph. Specifically, we
consider the following ensemble of random matrices   
\begin{equation}\label{eq:Ensplus}
\text{Ens}_+: \;  X = (x_{ij})_{\substack{1 \leq i \leq n \\ 1 \leq j \leq
p}}, \;  \{ x_{ij} \} \; \text{i.i.d. from a sub-Gaussian distribution on} \; \R_+.
\end{equation}
All random designs from the class $\text{Ens}_+$ share the property that the
population Gram matrix $\Sigma^{\ast} = \E[\frac{1}{n} X^{\T} X]$ possesses
equi-correlation structure after re-scaling the entries of $X$ by a common
factor. Denoting the mean of the entries and their squares by $\mu$ and
$\mu_2$, respectively, we have
\begin{equation*}
\Sigma^{\ast} = \E \left[\ffon X^{\T} X \right] = (\mu_2 - \mu^2) I + \mu^2 \bm{1} \bm{1}^{\T},
\end{equation*}
such that re-scaling by $1/\sqrt{\mu_2}$ leads to equi-correlation structure
with parameter $\rho = \mu^2 / \mu_2$. Since applications of NNLS predominantly
involve non-negative design matrices, it is instructive to have a closer
look at the class \eqref{eq:Ensplus} as a basic model for such designs. Among
others, the class of sub-Gaussian random designs on $\R_+$ encompasses the
zero-truncated Gaussian distribution, all
distributions on a bounded subset of $\R_+$, e.g.~the family of beta
distributions (with the uniform distribution as special case) on $[0,1]$,
Bernoulli distributions on $\{0,1 \}$ or more generally multinomial distributions on
positive integers $\{0,1,\ldots,K\}$, as well as any finite mixture of these distributions.\\  
As shown in the sequel, the class \eqref{eq:Ensplus} provides instances of
designs for which Theorems \ref{theo:oracle} to Theorems \ref{theo:ellinfsmin}
yield meaningful results in the $n < p$ setting. Our reasoning hinges on both
theoretical analysis providing bounds on the deviation from population
counterparts as well as on numerical results.\\
\\
\emph{Self-regularizing property + restricted eigenvalue condition of
  Theorem \ref{theo:oracle}.}  \hfill \\  
Recall that Theorem \ref{theo:oracle} requires a combination of the
self-regularizing property (Condition \ref{cond:1}) and the restricted
eigenvalue condition (Condition \ref{cond:2}) to be satisfied. This turns out
to be the case for designs from $\text{Ens}_+$ in light of he following proposition. The
statement relies on recent work of Rudelson and Zhou \cite{RudelsonZhou2012}
on the restricted eigenvalue condition for random matrices with
independent sub-Gaussian rows.
\begin{prop}\label{prop:re_and_selfreg}
Let $X$ be a random matrix from $\text{Ens}_+$ \eqref{eq:Ensplus} scaled such that 
$\Sigma^* = \E[\frac{1}{n} X^{\T} X] = (1 - \rho) I
+ \rho \bm{1} \bm{1}^{\T}$ for some $\rho \in (0,1)$. Set $\delta \in (0,1)$. There exists
constants $C, c > 0$ depending only on $\delta$, $\rho$ and the sub-Gaussian parameter of
the centered entries of $X$ so that if $n \geq C \, s \log (p \vee n)$, then, with probability at least 
$1 - \exp(-c \delta^2 n ) - 6/(p \vee n)$,
$\Sigma = X^{\T} X / n$  has the self-regularizing property with $\tau^2 =
\rho/2$ and satisfies the $(3/\tau^2, s)$ restricted eigenvalue condition of Theorem 2 with $\phi(3/\tau^2, s) = (1 - \rho)(1 - \delta)^2$. 
\end{prop}
\emph{Scaling of $\tau^2(S)$.} \hfill\\
The next proposition controls the deviation of the separating hyperplane
constant $\tau^2(S)$ from its population counterpart as derived in
\eqref{eq:tauS_equicor}.
\begin{prop}\label{prop:scalingtauS}
Let $X$ be a random matrix from $\text{Ens}_+$ \eqref{eq:Ensplus} scaled such that 
$\Sigma^* = \E[\frac{1}{n} X^{\T} X] = (1 - \rho) I
+ \rho \bm{1} \bm{1}^{\T}$ for some $\rho \in (0,1)$. Fix $S
\subset \{1,\ldots,p\}$, $|S| = s$. Then there exists constants $c, c', C,
C' > 0$ depending only on $\rho$ and the sub-Gaussian parameter of the
centered entries of $X$ such that for all $n \geq C s^2 \log (p \vee n)$,  
\begin{equation*}
  \tau^2(S) \geq c s^{-1} - C' \sqrt{\frac{\log p}{n}}
\end{equation*}
with probability no less than $1 - 6/(p \vee n) -  3\exp(-c' (s \vee \log n))$.
\end{prop}

\begin{figure}[ht!]
\vspace{-4pt} 
\begin{center}
\includegraphics[height = 0.17\textheight]{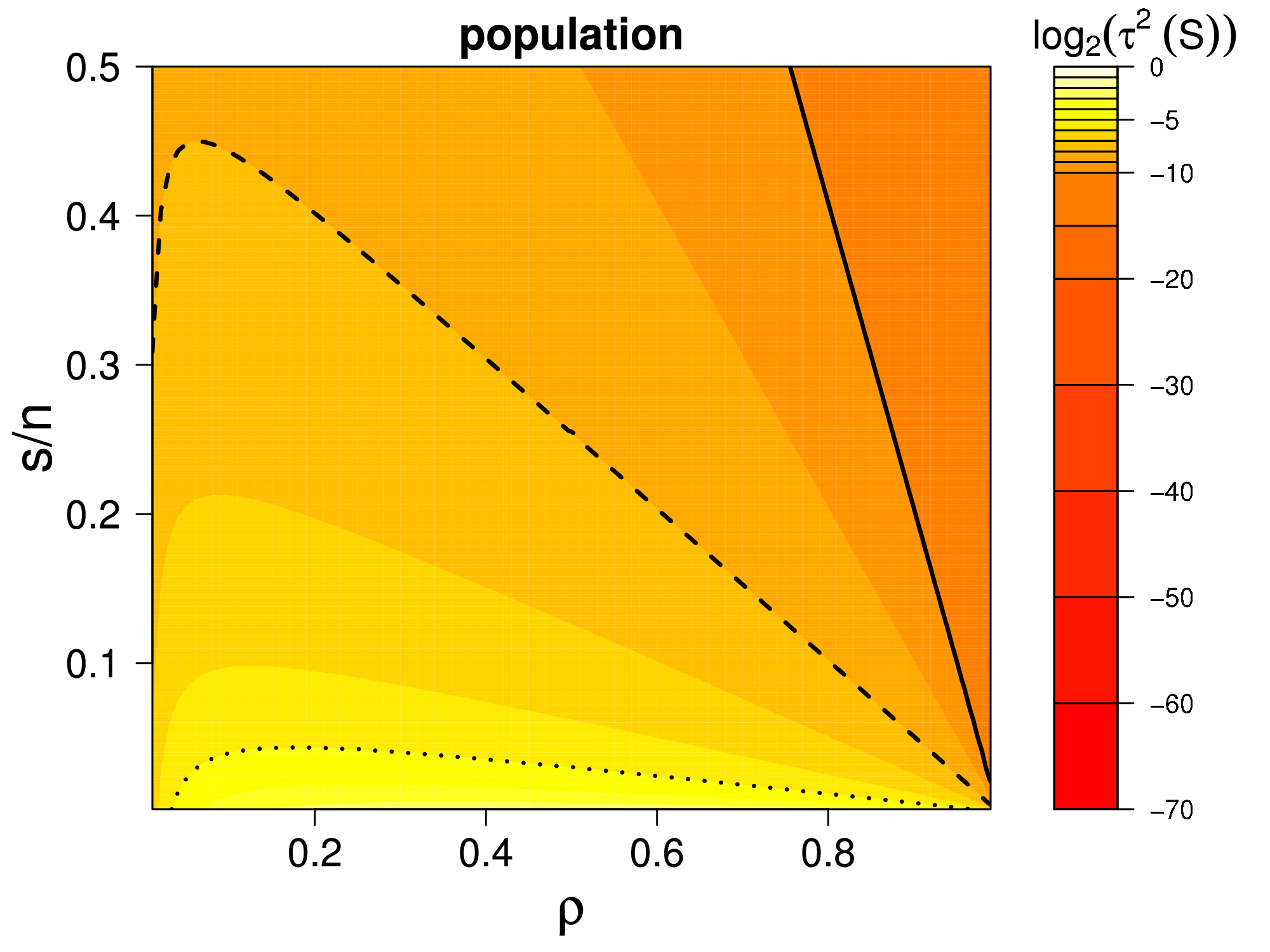}
\end{center}
\hfill\\
\hfill\\
\hfill\\
\begin{tabular}{ll}
 \hspace{.35cm}\includegraphics[height =
0.17\textheight]{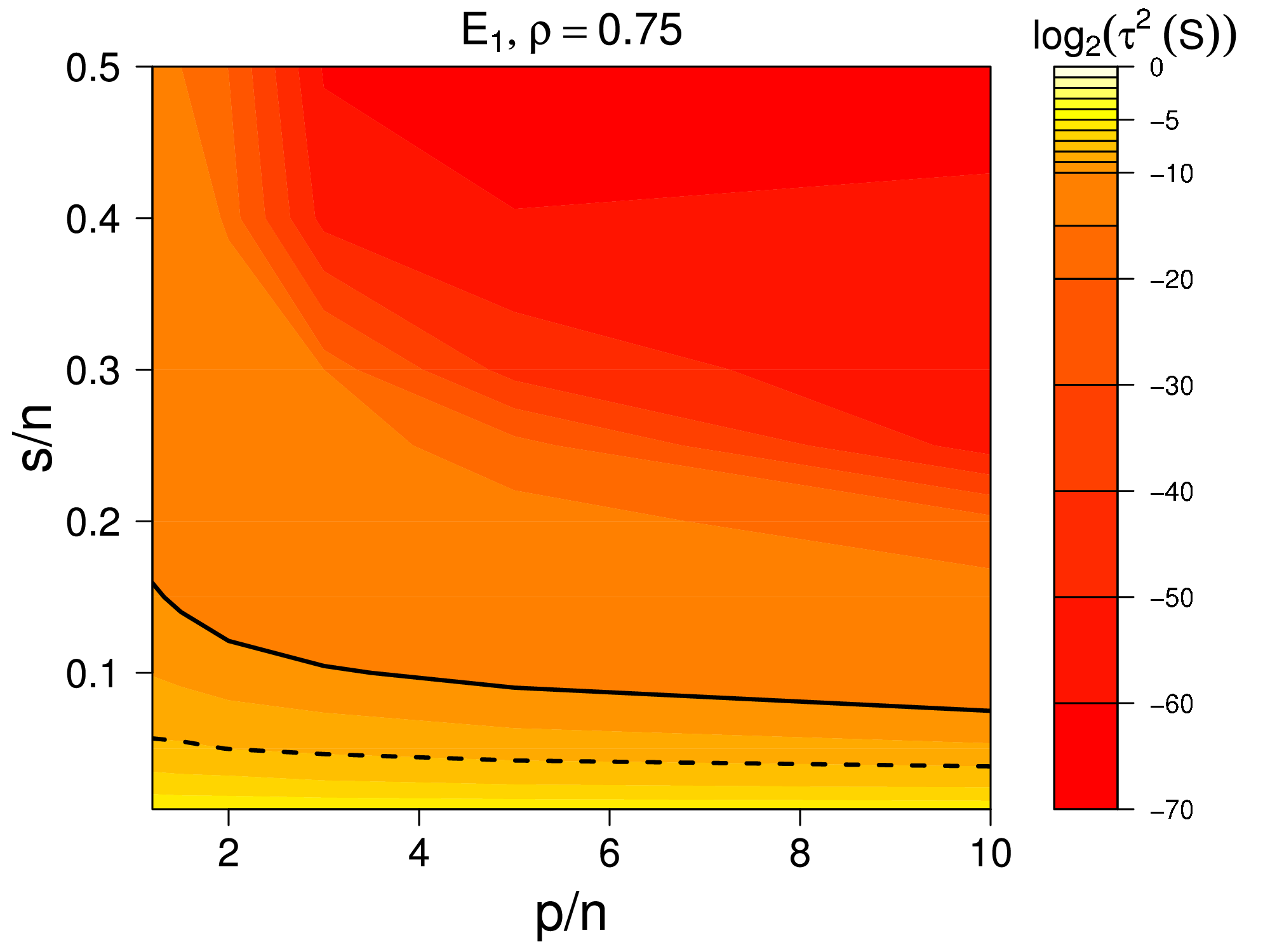}
& \hspace{2cm}\includegraphics[height =
0.17\textheight]{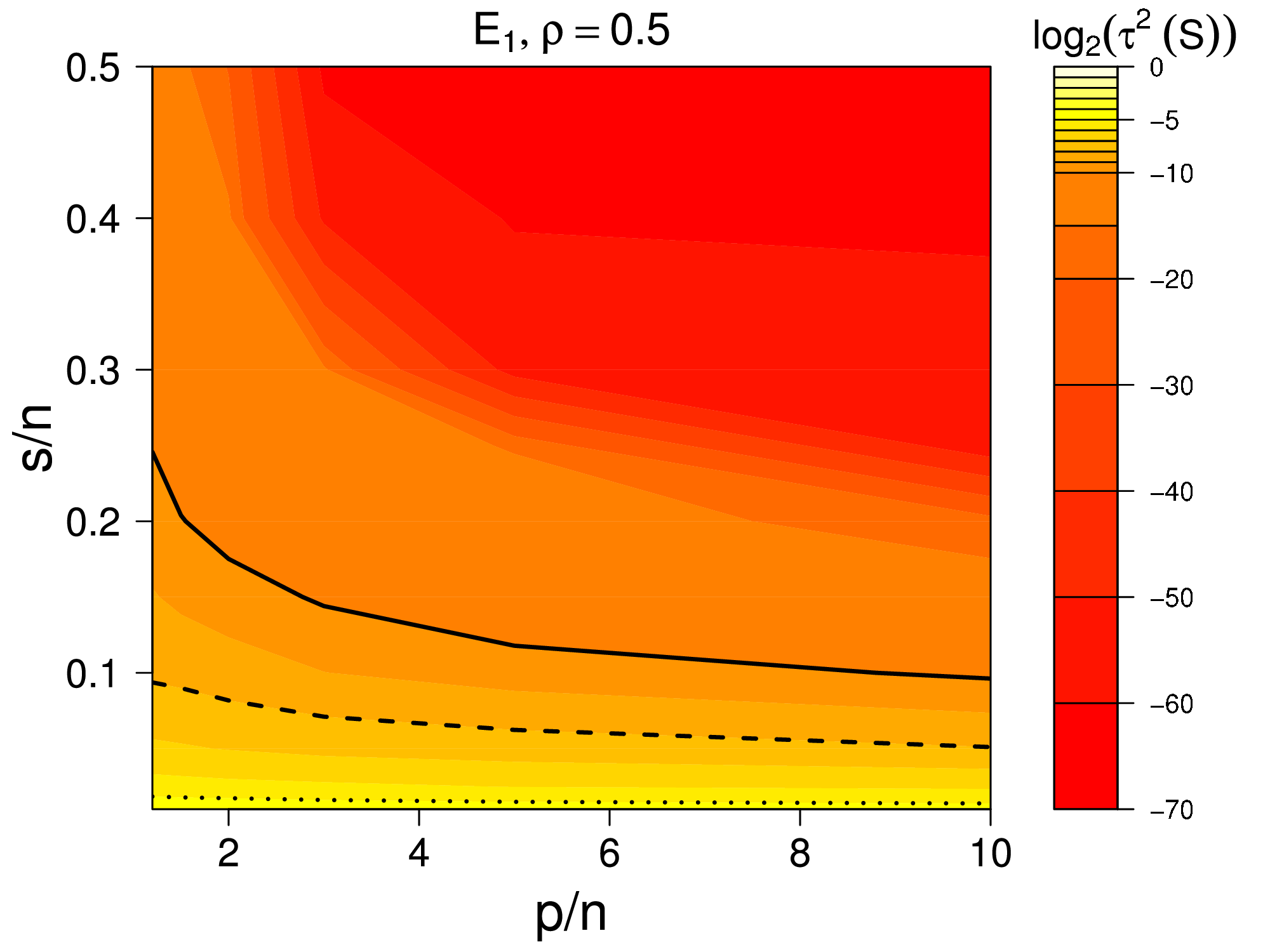} \\
  \\
\\
 \hspace{.35cm}\includegraphics[height =
0.17\textheight]{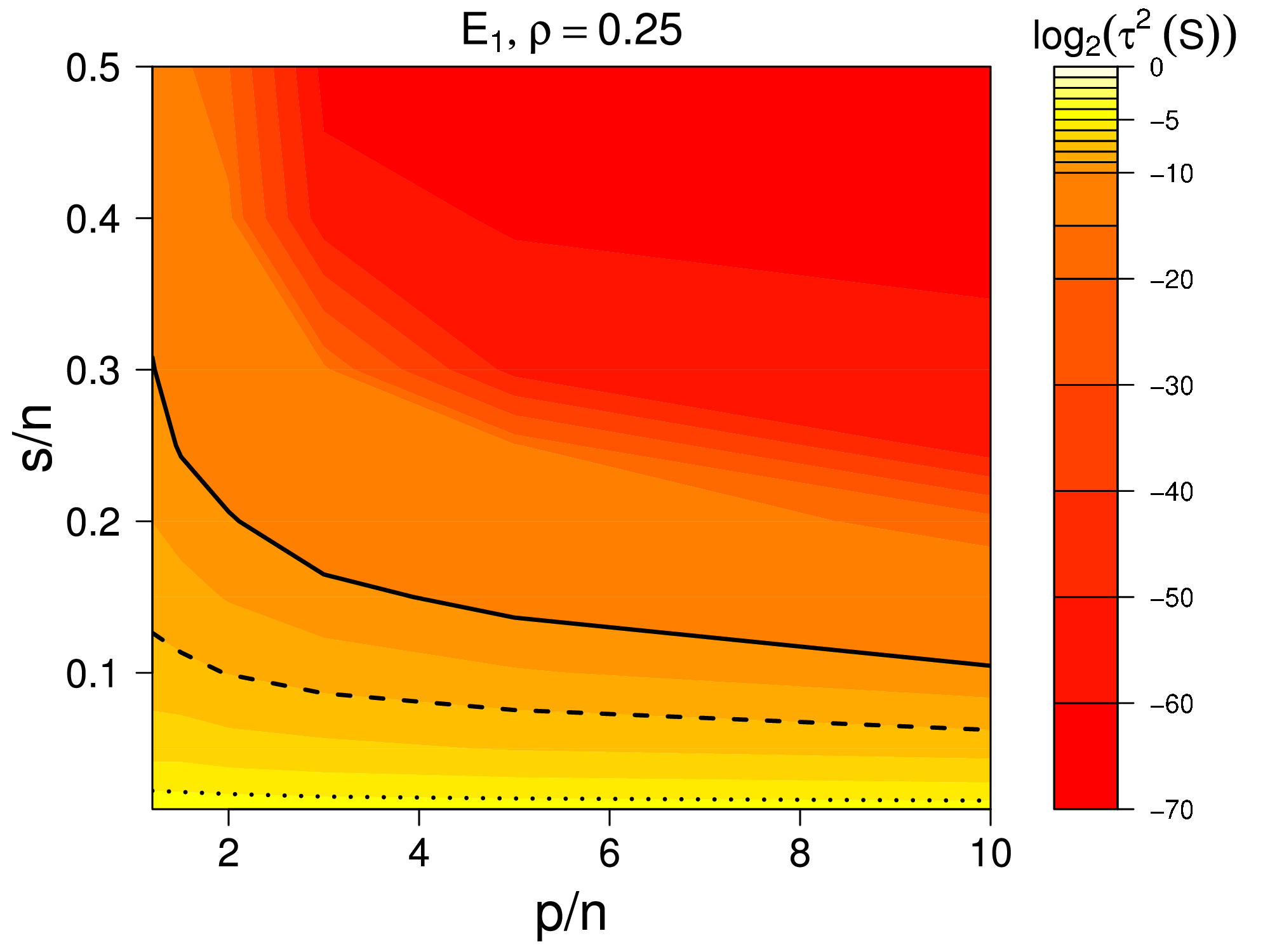} & \hspace{2cm}\includegraphics[height = 0.17\textheight]{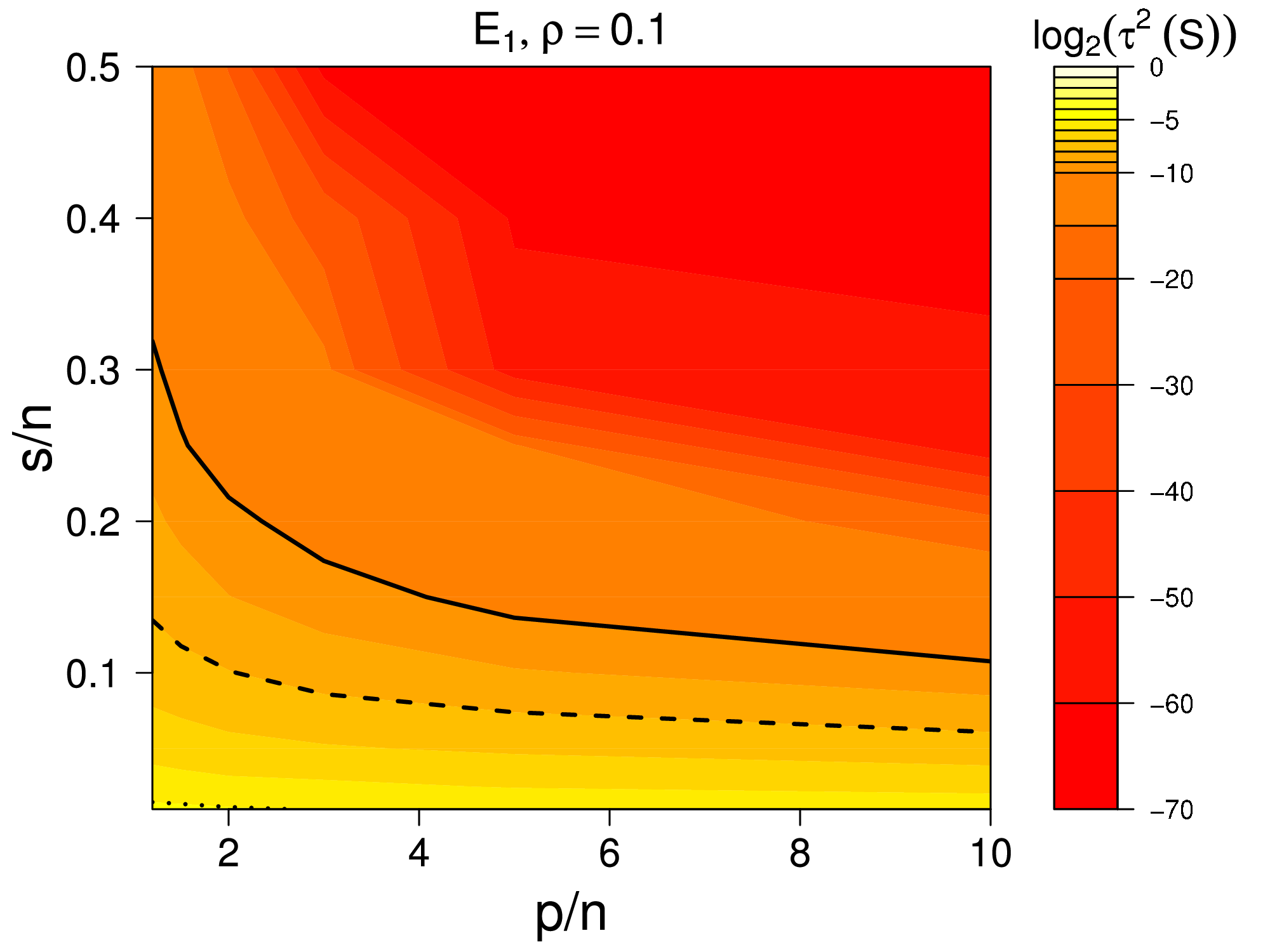}
\end{tabular}
\hfill\\
\hfill\\
\hfill\\
\vspace{-4pt}
\caption{Empirical scalings (0.05-quantiles over 100 replications) of the quantity $\log_2(\tau^2(S))$ for random
  design $E_1$ from the class $\text{Ens}_+$ in
  dependency of $s/n$ and $p / n$, displayed in form of a contour plot. 
The lines indicates the level set for $-10$
  (solid, $2^{-10} \approx 0.001$), $-8$ (dashed, $2^{-8} \approx
  0.004$) and $-5$ (dotted, $2^{-5} \approx 0.03$). The top plot displays
  $\log_2(\tau^2(S))$ for the population Gram matrix in dependency of $s/n$
  and $\rho \in (0,1)$.}\label{fig:tauS}
\end{figure}
It turns out that the requirement on the sample size as indicated by Proposition
\ref{prop:scalingtauS} is too strict in light of the results of complementary numerical experiments. For these experiments, 
$n = 500$ is kept fixed and $p \in (1.2,1.5,2,3,5,10) \cdot n$ and $s \in (0.01, 0.025, 0.05, 0.1, 0.15, 0.2,
0.25, 0.3, 0.4, 0.5) \cdot n$ vary. For each combination of $(p,s)$ and
several representatives of $\text{Ens}_+$ (re-scaled such that the population
Gram matrix has equi-correlation structure), 100 random design matrices are
generated. We set $S = \{1,\ldots,s \}$, compute $Z = (I - \Pi_S)
X_{S^c}$ using a QR decomposition of $X_S$ and then solve the quadratic
program $\min_{\lambda \in T^{p-s-1}} \lambda^{\T} \frac{1}{n} Z^{\T} Z \lambda$ with value
$\tau^2(S)$ by means of an interior point method
\cite{BoydVandenberghe2004}. As representatives of $\text{Ens}_+$, we have
considered matrices whose entries have been drawn from the following
distributions. In order to obtain population Gram matrices of varying
correlation $\rho$, we use mixture distributions with one of two mixture components being a point mass
at zero (denoted by $\delta_0$). Note that the larger the proportion $1-a$ of
that component, the smaller $\rho$.
{\small \begin{itemize}
\item[$E_1$:] $\{ x_{ij} \} \overset{\text{i.i.d.}}{\sim}  a \, \text{uniform}([0,\sqrt{3/a}]) + (1 - a) \delta_0$, $a \in \{1, \frac{2}{3},
  \frac{1}{3}, \frac{2}{15} \} \,$ ($\rho \in \{\frac{3}{4}, \frac{1}{2},
  \frac{1}{3}, \frac{1}{10} \}$)
\item[$E_2$:] $\{ x_{ij} \} \overset{\text{i.i.d.}}{\sim}
  \frac{1}{\sqrt{\pi}}\, \text{Bernoulli}(\pi), \; \pi \in
  \{\frac{1}{10}, \frac{1}{4}, \frac{1}{2}, \frac{3}{4}, \frac{9}{10}  \} \,$ 
  ($\rho \in \{\frac{1}{10}, \frac{1}{4}, \frac{1}{2}, \frac{3}{4},
  \frac{9}{10} \}$)
\item[$E_3$:] $\{ x_{ij} \} \overset{\text{i.i.d.}}{\sim} |Z|, \, Z \sim
  a \, \text{Gaussian}(0, 1) + (1 - a) \delta_0$, $a \in
  \{1,\frac{\pi}{4},\frac{\pi}{8},\frac{\pi}{20}\} \,$ ($\rho \in \{\frac{2}{\pi}, \frac{1}{2},
  \frac{1}{4}, \frac{1}{10} \}$)
\item[$E_4$:] $\{ x_{ij} \} \overset{\text{i.i.d.}}{\sim} \, a
  \text{Poisson}(3/\sqrt{12 \,a}) + (1 - a) \delta_0$, $a \in
  \{1,\frac{2}{3},\frac{1}{3},\frac{2}{15}\} \,$ ($\rho \in \{\frac{3}{4}, \frac{1}{2},
  \frac{1}{4}, \frac{1}{10} \}$)
\end{itemize}}   
For space reasons, we here only report the results for $E_1$. Regarding $E_2$ to $E_4$, the reader is referred
to the supplement; in brief, the results confirm
what is shown here. Figure \ref{fig:tauS} displays
the 0.05-quantiles of $\tau^2(S)$ over sets of 100 replications. It is
revealed that for $\tau^2(S)$ to be positive, $n$ does not need to be as large
relative to $s$ as suggested by Proposition \ref{prop:scalingtauS}. In fact, even for $s/n$
as large as $0.3$, $\tau^2(S)$ is sufficiently bounded away from zero as
long as $p$ is not dramatically larger than $n$ ($p/n = 10$).\\
\hfill \\
\begin{figure}
$$X:\;\bordermatrix{ & {\tiny \text{1}} & {\tiny \text{2}} & {\tiny \text{3}} & {\tiny \text{4}} & {\tiny \text{5}} &
  {\tiny \text{6}} & {\tiny \text{7}} & {\tiny \text{8}} & {\tiny \text{9}} & {\tiny \text{10}} \cr
{\tiny \text{1}} & 1 &  1 & 0 & 0 & 0 & 1 & 0 & 0 & 1 & 0 \cr
{\tiny \text{2}} & 0 &  1 & 0 & 0 & 0 & 0  & 1  & 0  & 0  & 1  \cr
{\tiny \text{3}} & 1 &  0 & 1 & 0 & 0 & 0  & 1  & 1 & 0 & 0 \cr 
{\tiny \text{4}} & 0 &  0 & 0 & 1 & 1 & 1 & 0 & 1& 0 & 0 \cr
{\tiny \text{5}} & 0 & 0 & 1 & 1 & 1 & 0 & 0 & 0 & 1 & 1 \cr
} \begin{array}{l}
         \text{find} \; \beta \gec 0 \; \text{s.t.} \; X \beta = X \beta^*: \\
         \Longrightarrow (1):\, \beta_1 + \beta_9 = 1 \, [\text{from} \, (2),(4)]\\
         \Longrightarrow (2):\,\beta_2 = \beta_7 = \beta_{10} = 0\\
         \Longrightarrow  (3):\,\beta_1 + \beta_3 = 1 \, [\text{from} \,
         (2),(4)]\\
         \Longrightarrow (4):\,\beta_4 = \beta_5 = \beta_6 = \beta_8 = 0\\
         \Longrightarrow (5):\,\beta_3 + \beta_9 = 2 \, [\text{from} \,
         (2),(4)] \\
         \hline
         \hline 
         \Longrightarrow \beta_1 = 0,\, \beta_3 = \beta_9 = 1 \,
         [\text{from} \, (1),(3), (5)]
\end{array}
$$
$$\beta^*:\;\;\,[\,0 \;\;\;\, 0  \;\;\;\, 1 \;\;\;\, 0 \;\;\;\, 0 \;\;\;\, 0 \;\;\;\, 0
\;\;\;\, 0 \;\;\;\, 1 \;\;\;\, 0\,]\qquad\quad \quad
\quad\qquad\qquad\qquad\qquad\qquad\qquad\quad\quad\;\,$$
\caption{Illustration of the setting. The rows of the sample measurement
  matrix $X$ represent five group assignments yielding (noiseless) measurements $y = X
  \beta^*$. A NNLS estimator is thus a solution of the systems of linear
  equations $y = X \beta$ subject to
  $\beta \gec 0$. A short calculation reveals that there is only a single
solution $\beta^*$ to this problem, even though it is severely underdetermined.}\label{fig:CS}
\end{figure}
\emph{Implications for compressed sensing-type problems.} \hfill \\ In compressed
sensing (CS) \cite{Donoho2006b, CandesTao2007}, the goal is to recover a sparse
vector $\beta^*$ from a limited number of non-adaptive measurements. In the
typical setup, the measurements are linear combinations of $\beta^*$,
contaminated with additive noise and thus fall under the linear model 
\eqref{eq:linearmodel}. CS is related to \emph{group testing} 
\cite{Dorfman1943, Du2006b}. Here, $\beta^* \in \{0,1\}^p$ indicates the
presence of a certain attribute of low prevalence in $p$ objects, e.g.~
presence of a rare disease in individuals or of a defect in manufactured
goods. An effective strategy for locating the affected entities is by
forming groups, testing for prevalence at the group level, discarding
groups with a negative test result, and repeating the procedure with
a new set of groups. The measurements obtained in this way are both
adaptive and non-linear and hence do not fit into the conventional framework
of CS. However, this is the case if it is possible to obtain aggregate
measurements (i.e.~sums) over arbitrary groups. As exemplified in 
Figure \ref{fig:CS}, the information of interest can be retrieved from few 
aggregated measurements and the associated group assignments. Proper
measurement design needs to achieve a proper amount of overlapping of the
groups. In fact, it is not possible to recover $\beta^*$ from a reduced 
number of measurements involving only disjoint groups. At the same time,
overlapping has to be limited to ensure that collections of $2s$ columns
of the measurement matrix are linearly independent; otherwise, recovery
of $\beta^*$ is not possible in general. Without further prior knowledge,
about the location of the non-zero entries of $\beta^*$, random group
assignments in which each entity $j$, $j=1,\ldots,p$, is assigned to
group $i$, $i=1,\ldots,n$, independently with probability $\pi$, appears
to be reasonable in light of the results of the previous paragraph.
Propositions \ref{prop:re_and_selfreg} and \ref{prop:scalingtauS} assert 
that, with high probability, the resulting measurement matrix suits well a
sparse recovery approach based on NNLS.\\
Network tomography as discussed in \cite{Meinshausen2013} arises
as a generalization of the above setting with measurements of the 
form $y = B A \beta^* + \eps$, where $B \in \R_+^{n \times p}$ represents the
measurement design and $A \in \R_+^{p \times p}$ is the adjacency matrix
of the $p$ nodes whose status of interest is contained in $\beta^*  \in
\R_+^p$. The goal is to spot sources of anomaly within the network. As
distinguished from a random design setting, the measurement matrix
$X = B A$ cannot be chosen freely, since $A$ is fixed and the choice of 
$B$ is subject to various constraints. Therefore, it is not clear a priori
whether $X$ satisfies our conditions, even though $\Sigma$ is potentially self-regularizing as $X$ is non-negative. 
\hfill \\
\hfill\\
\\
\emph{Sparsity of the solution.} \hfill \\
In Proposition \ref{prop:sparsity_equicor}, we have characterized the sparsity
in the population setting. It is of interest to investigate this aspect for
random design in the $p > n\,$ setup, particularly in light of Theorem
\ref{theo:ellinfsmin}, which implicitly relies on having a sparse NNLS
solution. We here provide a sketch of the empirical behaviour within the
experimental framework of the previous paragraph. We generate random design
matrices ($n \in \{250, 500, 750, 1000 \}$, $p/n \in \{2,5,10\}$) from $E_1$ for the four values of
the parameter $\rho$ as given above. For several values of $s/n$ ranging from
$0$ to $0.3$, we generate observations $y = X \beta^* + \eps$, where $\eps$ is a Gaussian noise
vector, and the components of $\beta_S^*$ are set to the lower bound in
Proposition \ref{prop:sparsity_equicor}. For each
combination of $(n, p/n, s/n, \rho)$, 100 replications are considered and the 
fraction of active variables $|F|/n$ is determined.
\begin{flushleft}
\begin{figure}[h!]
\begin{tabular}{lll}
\hspace{-0.6cm}\includegraphics[height =
0.18\textheight]{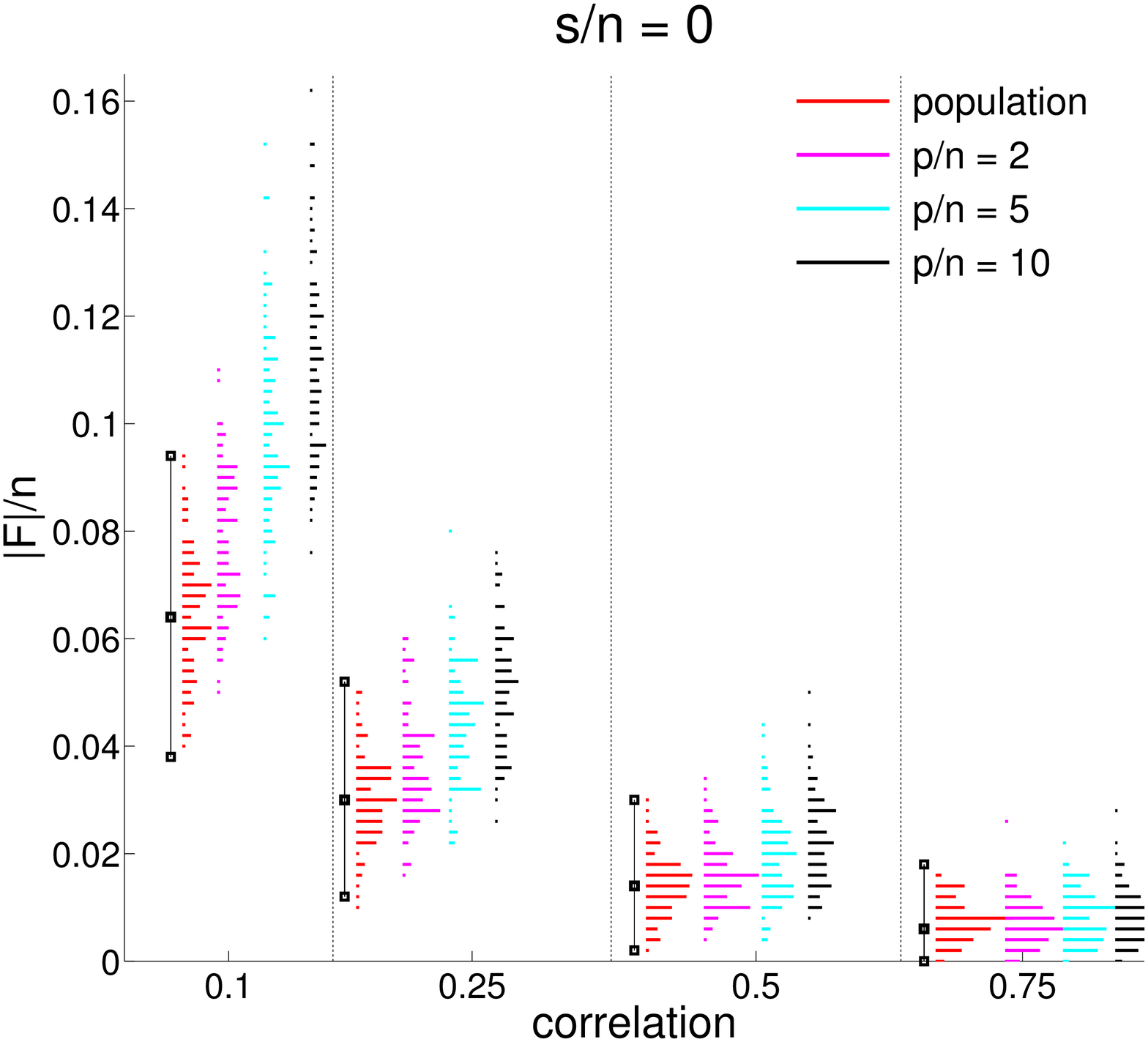} & \hspace{-0.7cm}
\includegraphics[height =
0.18\textheight]{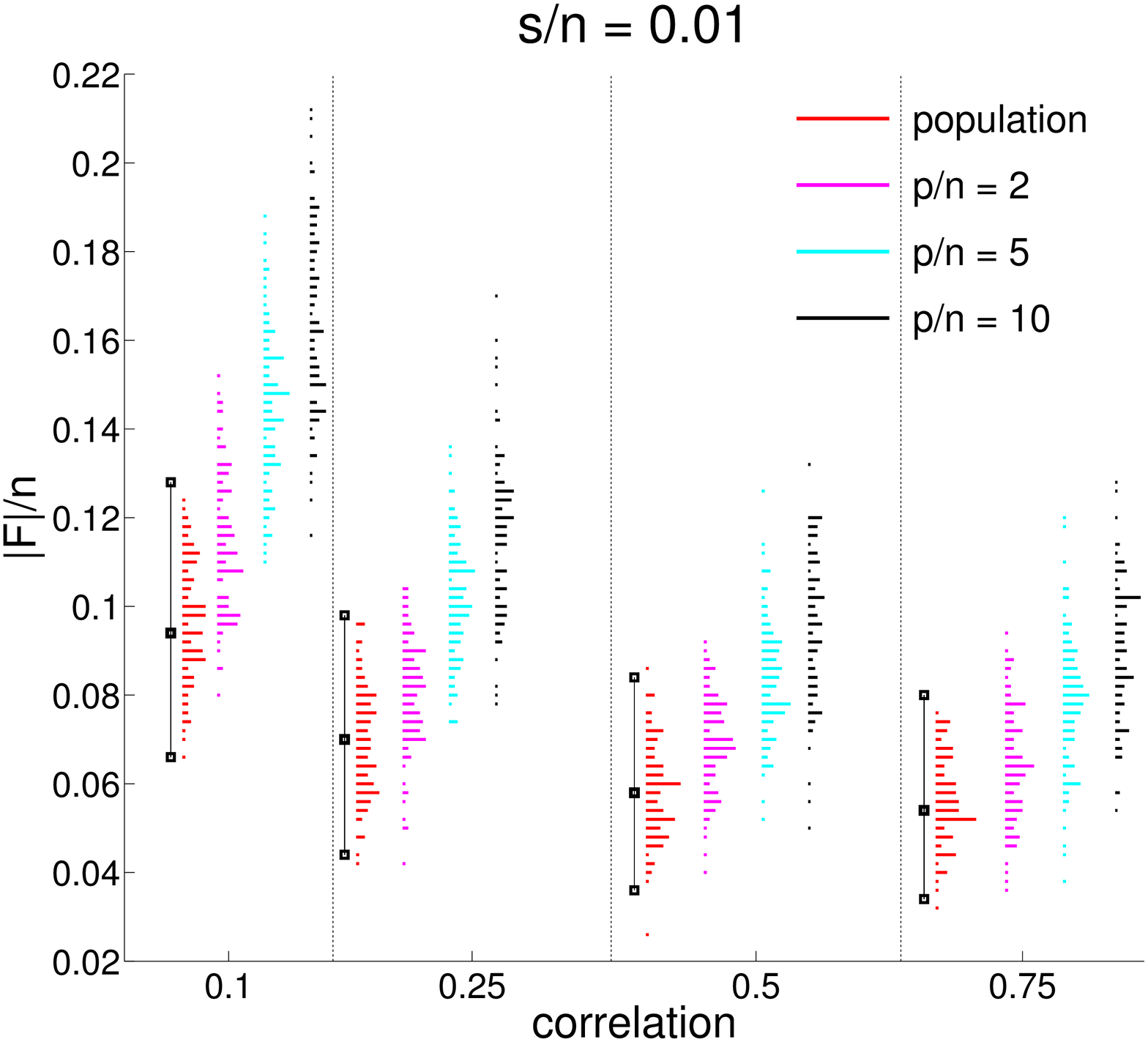} &
\hspace{-0.41cm}\includegraphics[height =
0.18\textheight]{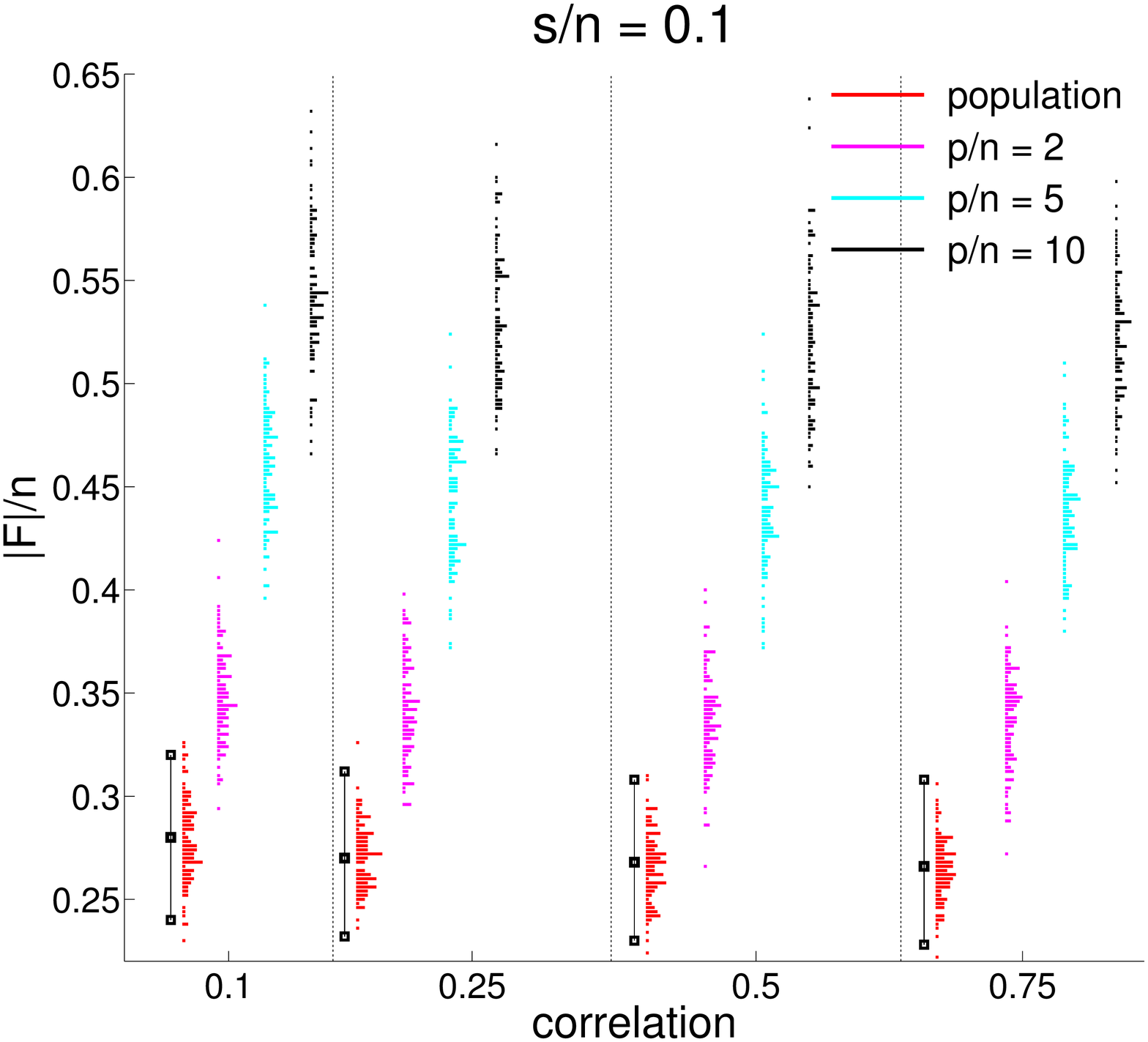}\\
\end{tabular}
\begin{tabular}{lll}
\hspace{-0.2cm}
\includegraphics[height =
0.18\textheight]{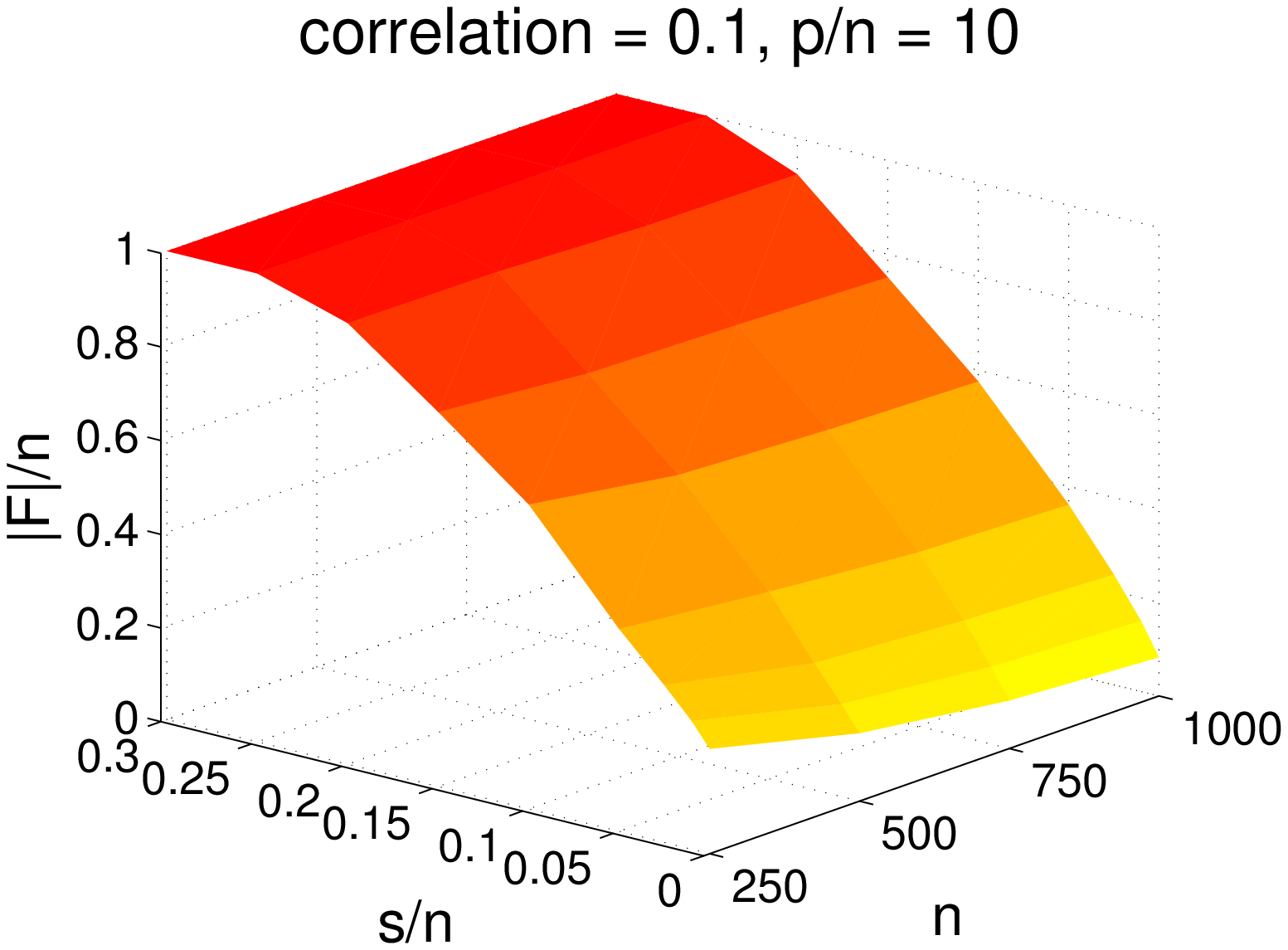} & 
\hspace{-1.45cm}
\includegraphics[height =
0.18\textheight]{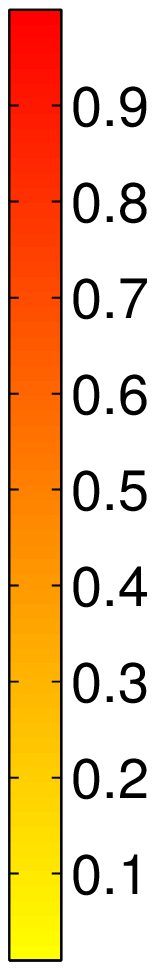}
& \hspace{-1.5cm}
\includegraphics[height =
0.18\textheight]{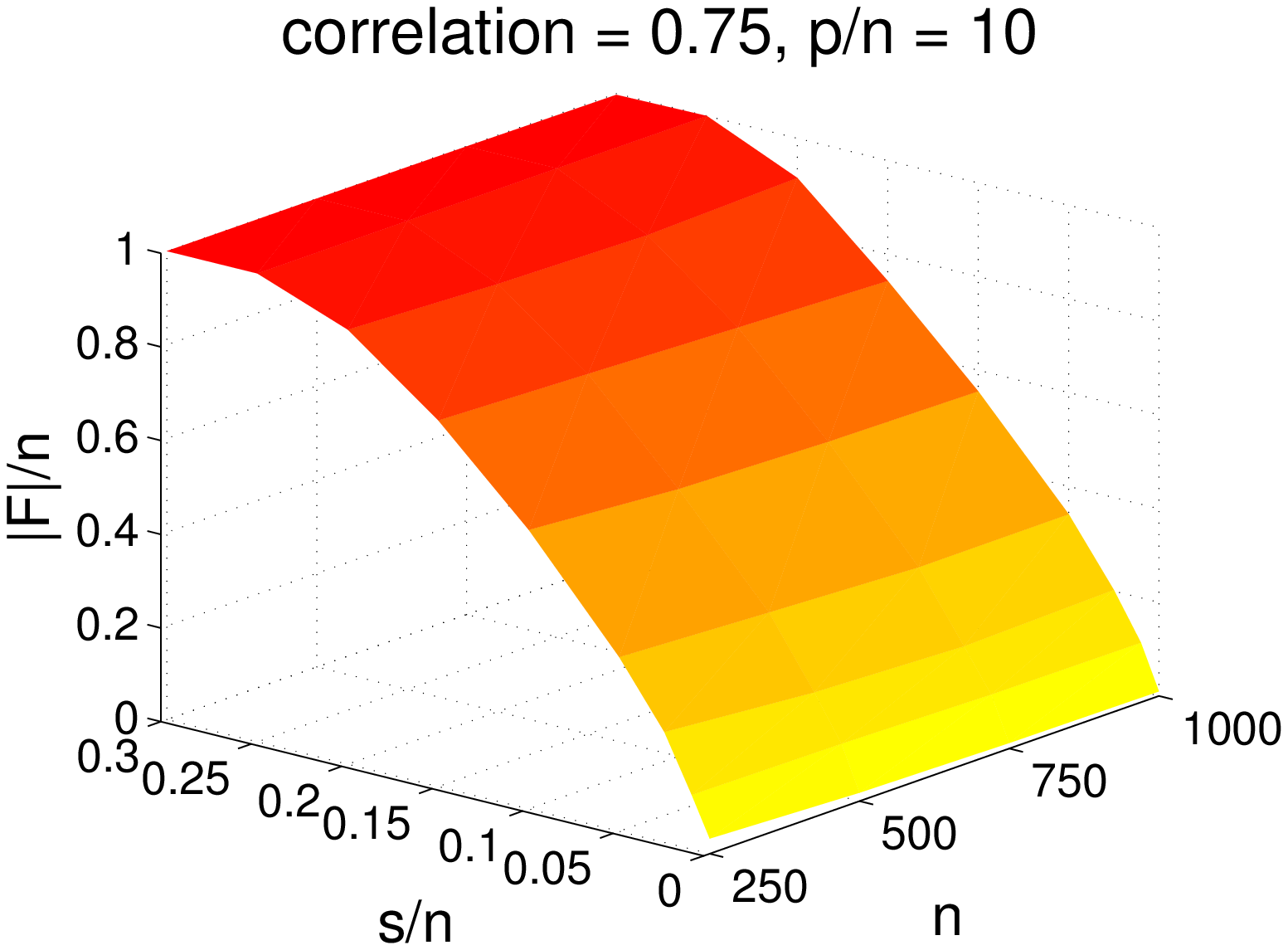}
\end{tabular}
\caption{Top: sparsity of the NNLS solution for random equi-correlation-like design
in the $n < p$ setup as compared to the population setting. The squares
represent the 0.01-, 0.5- and 0.99-quantiles of the (conditional) distribution
of the fraction of active variables $|F|/n$ in the population according to Proposition
\ref{prop:sparsity_equicor}. The vertical bars represent the empirical
distributions over 100 random datasets with $n = 500$, where the colours correspond to
different ratios $p/n$. Bottom: Surface plot of the 0.95-quantiles of $|F|/n$ over
100 random datasets for $n$ and $s/n$ varying.}\label{fig:cardF_randomdesign}
\end{figure}
\end{flushleft}
\vspace{-24pt}
Figure \ref{fig:cardF_randomdesign} summarizes the main findings of this experimental
study. For fixed $n = 500$, the top panel depicts the empirical
distributions of $|F|/n$  over the 100 replications in
comparison to the population setting (cf.~Figure \ref{fig:cardF_equicor}). We
observe that for all parameter configurations under consideration, the
cardinalities of the active sets stay visibly away from $1$ with 
$|F|/n$ being no larger than $2/3$. The cardinalities of the active sets are larger
than in the population case. The higher the sparsity level and the ratio
$p/n$, the more pronounced the shifts toward larger cardinalities: while for 
$s/n = 0$ and $\rho = 0.75$, the empirical distribution of $|F|/n$ is rather close
to that of the population, there is a consistent gap for $s/n = 0.1$. The
bottom panel displays how $|F|/n$ scales with $(n, s/n)$. For plotting and space reasons,
we restrict us to the $0.95$-quantiles over the 100 replications and 
$p/n = 10$, which, as indicated by the plots of the top panel, is the worst
case among all values of $p/n$ considered. The two surface plots for $\rho = 0.1$ and
$\rho = 0.75$ are of a similar form; a noticeable difference occurs only
for rather small $s/n$. It can be seen that for $s/n$ fixed, $|F|/n$ roughly
remains constant as $n$ varies. On the other hand, $|F|/n$ increases rather
sharply with $s/n$. For $s/n > 0.25$, we observe a breakdown, as $|F|/n = 1$.
We point out that as long as $|F|/n < 1$, it holds that the NNLS solution and the
active set are unique (with probability one), as follows from Lemma \ref{lem:uniqueness} 
in Appendix \ref{app:ellinfsmin}. 
\section{Empirical performance}\label{sec:empiricalresults}

We here present the results of simulation studies in order to compare the
performance of NNLS and the non-negative lasso in terms of prediction,
estimation and sparse recovery. 

\subsection{Deconvolution of spike trains}\label{subsec:sparsedeconvoultion}

We consider a positive spike-deconvolution model as in \cite{Speed2000}, as
it commonly appears in various fields of applications. The
underlying signal $f$ is a function on $[a,b]$ of the form
\begin{equation*}
f(u) = \sum_{k = 1}^s \beta_k^{\ast} \phi_k(u),  
\end{equation*}
with $\phi_k(\cdot) = \phi(\cdot \; - \mu_k)$, $k=1,\ldots,s$, where $\phi
\geq 0$ is given and the $\mu_k$'s define the locations of the
spikes contained in $[a,b]$. The amplitudes $\{ \beta_k^* \}_{j = 1}^s$ are
assumed to be positive. The goal is to determine the positions as well as the amplitudes of the
spikes from $n$ (potentially noisy) samples of the underlying signal $f$.
As demonstrated below, NNLS can be a first step towards deconvolution. The
idea is to construct a design matrix of the form $X = (\phi_j(u_i))$, where       
$\phi_j = \phi(\cdot \; - m_j)$ for candidate positions $\{m_j \}_{j = 1}^p$
placed densely in $[a,b]$ and $\{ u_i \}_{i = 1}^n \subset [a,b]$ are the
points at which the signal is sampled. Under an additive noise model with
zero-mean sub-Gaussian noise $\eps$, i.e.
\begin{equation}\label{eq:deconvolutionmodel}
y_i = \sum_{k = 1}^s \beta_k^{\ast} \phi_k(u_i) + \eps_i, \quad i=1,\ldots,n,
\end{equation}
and if $X$ has the self-regularizing property (cf.~Section
\ref{sec:nnlsdoesnotoverfit}), it follows immediately from Theorem
\ref{theo:prediction} that the $\ell_2$-prediction error of NNLS is bounded as 
\begin{equation}\label{eq:excessriskbound}  
\frac{1}{n} \nnorm{f - X \wh{\beta}}_2^2 \leq \mc{E}^* + C \sqrt{\frac{\log
    p}{n}}, \quad \; \; \text{where} \; \; \{f_i = f(u_i)
\}_{i = 1}^n,
\end{equation}
where $\mc{E}^{\ast} = \min_{\beta \gec 0} \frac{1}{n} \nnorm{f - X \beta}_2^2$.
Even though it is not realistic to assume that $\{\mu_k \}_{k =1}^s \subset \{m_j
\}_{j = 1}^p$, i.e.~that the linear model is
correctly specified, we may think of $\mc{E}^{\ast}$ being negligible as long
as the $\{m_j \}_{j = 1}^p$ are placed densely enough. This means that NNLS
may be suitable for de-noising. Furthermore, the bound \eqref{eq:excessriskbound}
implies that $\wh{\beta}$ must have large components only for those columns of $X$
corresponding to locations near the locations $\{\mu_k \}_{k=1}^{s}$ of the
spikes, which can then be estimated accurately by applying a simple form of
post-processing as discussed in \cite{SlawskiHein2010}.
On the other hand, the application of fast rate bounds such as that of Theorem 
\ref{theo:oracle} or corresponding results for the lasso is not adequate here, because the dense placement of 
the $\{ \phi_{j} \}_{j = 1}^p$ results into a tiny, if not zero, value
of the restricted eigenvalue of Condition \ref{cond:2}.
\begin{figure}[h!]
\begin{center}
\begin{tabular}{lr}
\includegraphics[height =
0.135\textheight]{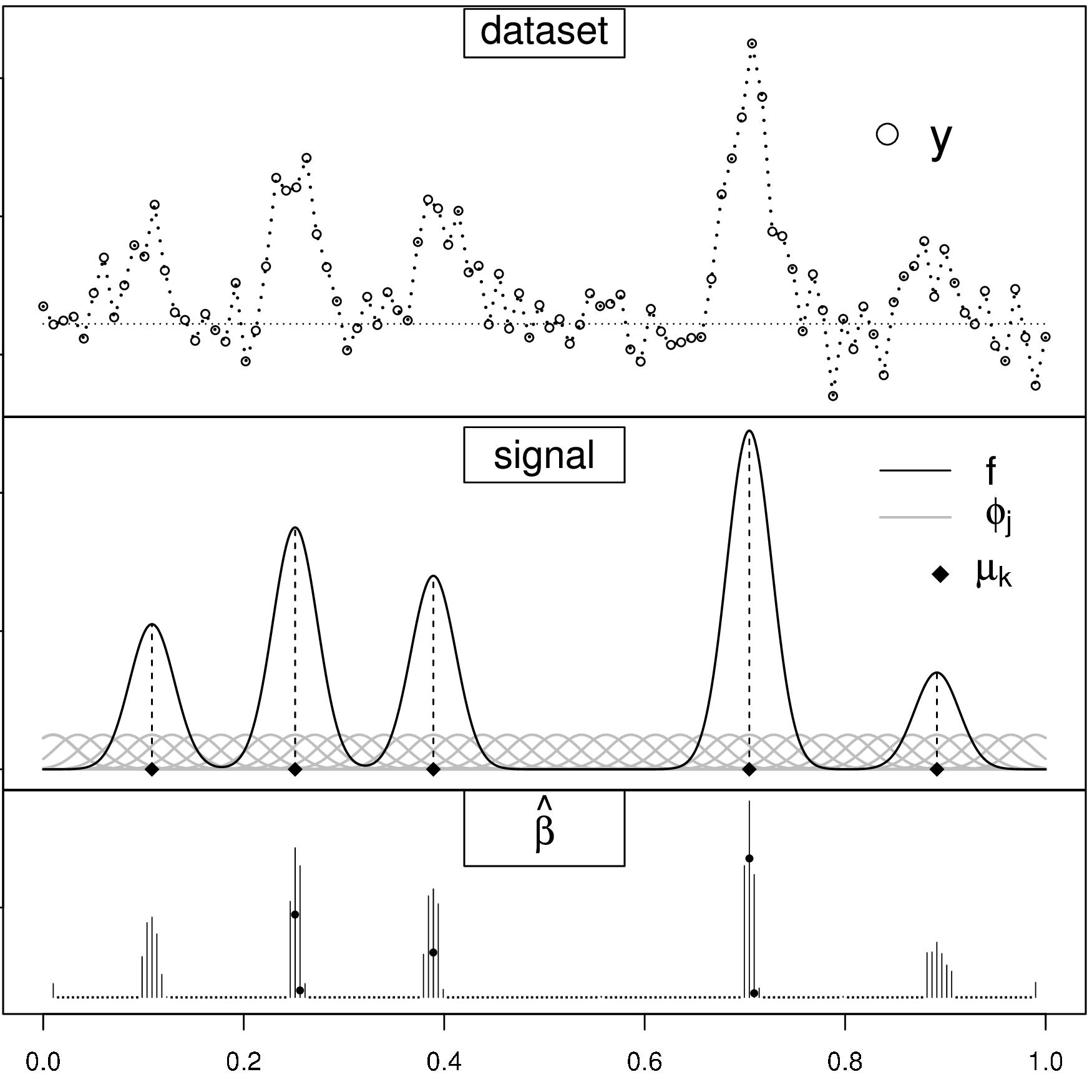}
&  \includegraphics[height = 0.25\textheight]{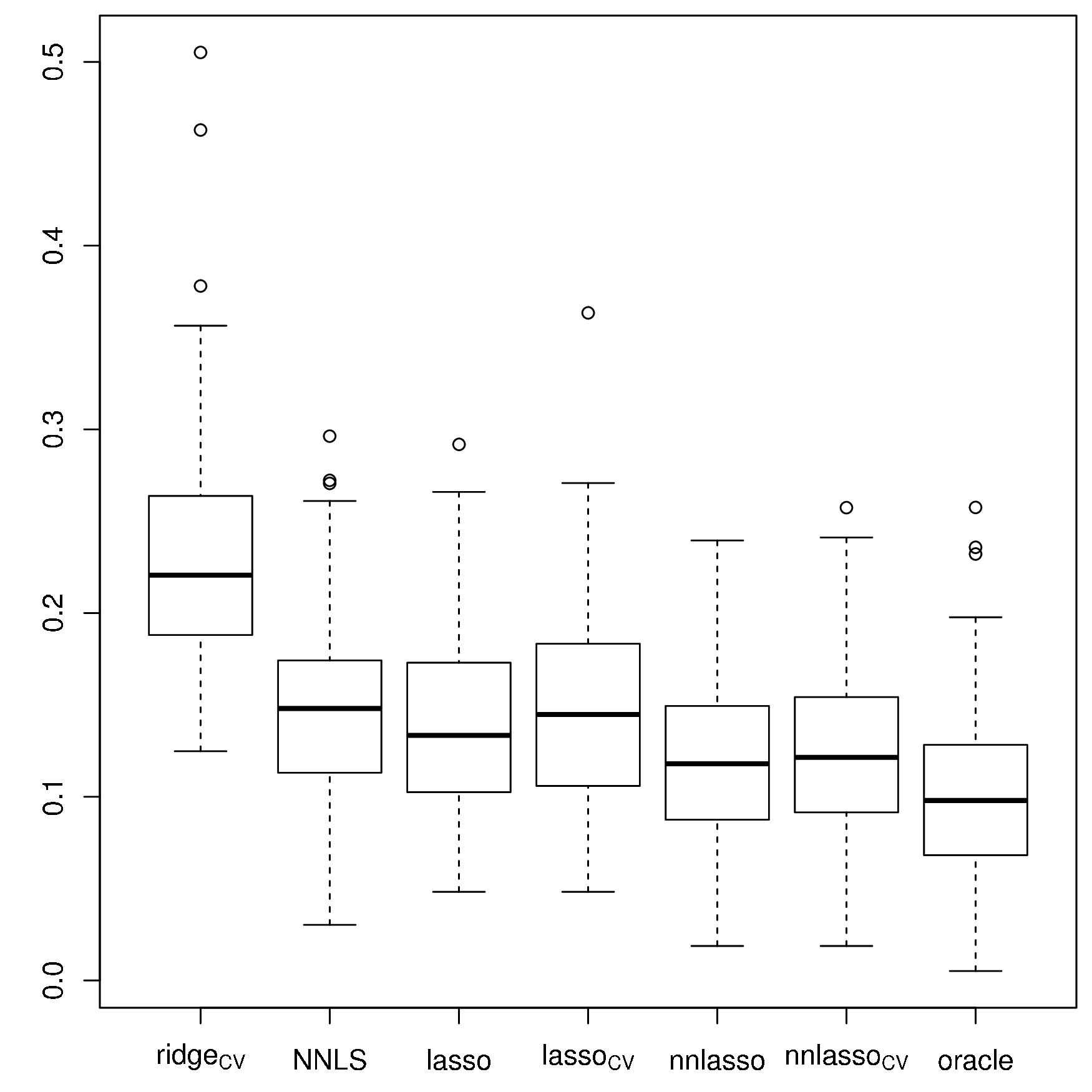} \\
\end{tabular}  
\caption{Left panel: Visualization of the experimental setting. 
The middle part depicts the underlying signal, a positive combination of five
Gaussians. The upper part depicts a sample dataset generated according to
model \eqref{eq:deconvolutionmodel}. The lower part provides a summary of the
coefficient vectors $\wh{\beta}$ returned by NNLS, the heights of the bars
representing the 0.9-quantiles and the dots the non-zero median coefficients at the
respective positions over 100 replications. Right panel: Boxplots of the mean squared
prediction errors (MSEs).}
\label{fig:4}
\end{center}
\end{figure}
For our simulation study, we consider model \eqref{eq:deconvolutionmodel} as
starting point. The signal is composed of five spikes of amplitudes between $0.2$ and $0.7$
convolved with a Gaussian function. The design matrix $X = (\phi_j(u_i))$ contains evaluations of $p = 200$ Gaussians
$\{ \phi_{j} \}_{j = 1}^p$ at $n = 100$ points $\{ u_i \}_{i = 1}^n$,
where both the centers $\{ m_j \}_{j = 1}^p$ of the $\{ \phi_j \}_{j = 1}^p$
as well as the $\{ u_i \}_{i =1}^n$ are equi-spaced in the unit
interval. We have by construction that $\{ m_j \}_{j = 1}^p \supset \{\mu_k
\}_{k=1}^{s}$ so that $\mc{E}^* = 0$. The standard deviation of the Gaussians is chosen such that it is roughly twice the spacing of the
$\{u_i \}$. At this point, it is important to note that the larger the
standard deviations of the Gaussians, the larger the constant
$\tau_0^2$ \eqref{eq:tau0}, which here evaluates as $\tau_0^2 = 0.2876$. According to that setup, we generate 100 different vectors
$y$ resulting from different realizations of the noise
$\eps$ whose entries are drawn i.i.d.~from a Gaussian distribution with
standard deviation $\sigma = 0.09$. The left panel of Figure \ref{fig:4} visualizes
the setting. We compare the performance of NNLS, lasso/non-negative lasso with
(i) fixed regularization parameter $\lambda$ fixed to $\lambda_0 = 2 \sigma \sqrt{2
  \log(p)/n}$ (ii) $\lambda$ chosen from the grid $\lambda_0 \cdot 2^k,
k=-5,-4\ldots,4$ by tenfold cross-validation, ridge regression (tuned by
tenfold cross-validation) and an oracle least squares
estimator based on knowledge of the positions $\{ \mu_k \}_{k = 1}^s$ of the
spikes. The right panel of Figure \ref{fig:4} contains boxplots of the MSEs
$\frac{1}{n} \nnorm{X \beta^* - X \wh{\beta}}_2^2$ over all 100
replications. The performance of NNLS is only slightly inferior to that of the
non-negative lasso, which is not far from the oracle,  and roughly as good as that of the lasso. All methods
improve substantially over ridge regression. The lower part of the left panel provides a summary of the       
NNLS estimator $\wh{\beta}$, which is remarkably sparse and concentrated near the positions of the
underlying spikes.

\subsection{Sparse recovery}

We now present the results of simulations in which we investigate the
performance of NNLS with regard to estimation and sparse recovery in
comparison to that of the non-negative lasso.\\
\textbf{Setup.} We generate data $y = X \beta^* + \eps$, where $\eps$ has i.i.d.~standard Gaussian
entries. For the design $X$, two setups are considered.\hfill\\
\emph{Design I: Equi-correlation like design.}\\
The matrix $X$ is generated by drawing its entries independently from the
uniform distribution on $[0,1]$ and rescaling them such that the population
Gram matrix is of equi-correlation structure with $\rho = 3/4$. Random matrices of
that form have been considered for the numerical results in Section
\ref{eq:equicorrelationlike} (cf.~Figures \ref{fig:cardF_equicor} to \ref{fig:cardF_randomdesign}). For given
$(n,p,s)$, the target $\beta^*$ is generated by setting $\beta_j^* = 6 b \cdot
\phi_{\min}^{-1/2 } \sqrt{2 \log(p)/n}(1 + U_j)$, $j=1,\ldots,s$, where
$\phi_{\min} = (1-\rho)$ denotes the smallest eigenvalue of the population Gram
matrix, the $\{U\}_{j = 1}^s$ are drawn uniformly from $[0,1]$, and we let
the parameter $b > 0$ vary. We further set $\beta_j^* = 0$, $j = (s+1),\ldots,p$.\\
\emph{Design II: Localized non-negative functions.}\\
The setup leading to the second class of designs can be regarded as a
simplification of the deconvolution problem discussed in 
the previous subsection to fit into the standard sparse recovery
framework. Indeed, in the experiments of Section \ref{subsec:sparsedeconvoultion}, 
recovery of the support of $\beta^*$ fails in the presence of noise, because the
$\{\phi_j\}$'s are placed too densely relative to the number of sampling
points; see \cite{Bunea2010} for a similar discussion concerning the
recovery of mixture components in sparse mixture density estimation. In order
to circumvent this issue, we use the following scheme. As in Section
\ref{subsec:sparsedeconvoultion}, we consider sampling points $u_i = i/n$,
$i=1,\ldots,n$, in $[0,1]$ and localized functions $\phi_j = \phi(\cdot -
m_j)$,where here $\phi(\cdot - m_j) = \exp(-|\cdot - m_j|/h)$,
$j=1,\ldots,p$ with $h = 2/n$. The centers $m_j$, $j=1,\ldots,p$, are taken from the interval
$[m_{\min}, m_{\max}]$, where $m_{\min} = u_1 - h \log(1/n)$ and $m_{\max} =
u_n + h \log(1/n)$. Given the sparsity level $s$, $[m_{\min}, m_{\max}]$ is
partitioned into $s$ sub-intervals of equal length and the centers
$m_1,\ldots,m_s$ corresponding to $S$ are drawn from the uniform distributions on these
intervals. The remaining centers $m_{s+1},\ldots,m_p$ corresponding to $S^c$
are drawn uniformly from $[m_{\min}, m_{\max}] \setminus \cup_{j = 1}^s
[m_j - \Delta, m_j + \Delta]$, where $\Delta > 0$ is set to enforce a
sufficient amount of separation of the $\{ \phi_j \}_{j = 1}^s$ from the $\{
\phi_j \}_{j = s+1}^p$. We here choose $\Delta = h = 2/n$. The design matrix
is then of the form $X_{ij} = \phi_j(u_i)/c_j$, $i=1,\ldots,n$,
$j=1,\ldots,p$, where the $c_j$'s are scaling factors such that $\nnorm{X_j}_2^2 = n \, \forall j$. As
for Design I, we generate observations $y = X \beta^* + \eps$, where $\beta_j^* = b \cdot \beta_{\min} (1 +
U_j)$, $j = 1,\ldots,s$ and $\beta_j^* = 0$, $j=s+1,\ldots,p$. The $\{U_j\}_{j
= 1}^s$ are random variables from the uniform distribution on $[0,1]$ and
the choice $\beta_{\min} = 4  \sqrt{6 \log(10)/n}$
has turned out to yield sufficiently challenging problems.\\
For both Design I and II, two sets of experiments are performed. In the first
one, the parameter $b$ controlling the magnitude of the coefficients of the
support is fixed to $b = 0.5$ (Design I) respectively $b = 0.55$ (Design II) , while the aspect ratio $p/n$ of $X$ and the
fraction of sparsity $s/n$ vary. In the second set of experiments,
$s/n$ is fixed to $0.2$ (Design I) and $0.05$ (Design II), while $p/n$ and $b$
vary. Each configuration is replicated 100 times for $n = 500$.\\
\textbf{Comparison.} Across these runs, we compare  thresholded NNLS, the non-negative lasso (NN$\ell_1$), the thresholded non-negative
lasso (tNN$\ell_1$) and orthogonal matching pursuit (OMP, \cite{Tropp2004, Zhang2009})
with regard to their performance in sparse recovery. Additionally, we compare 
NNLS and NN$\ell_1$ with $\lambda = \lambda_0$ as defined below (both
\emph{without} subsequent thresholding) with regard to
estimation of $\beta^*$ in $\ell_{\infty}$-norm (Tables
\ref{tab:ellinf_uniform} and \ref{tab:ellinf_block}) and
$\ell_2$-norm (the results are contained in the supplement).
The performance of thresholded NNLS with regard to sparse recovery is assessed in two
ways. For the first one (referred to as 'tNNLS$^*$'), success is reported whenever $\min_{j \in
  S} \wh{\beta}_j > \max_{j \in S^c} \wh{\beta}_j$, i.e.~there exists a
threshold that permits support recovery. Second, the procedure of Theorem
\ref{theo:Sdatadriven} (with $\sigma$ replaced by the naive estimate
$\frac{1}{n} \nnorm{y - X \wh{\beta}}_2^2$) is used to determine the threshold in a data-driven
manner without knowledge of $S$. This approach is referred to as tNNLS. For tNN$\ell_1$, both the regularization parameter $\lambda$
and the threshold have to be specified. Instead of fixing $\lambda$ to a single value, we 
give tNN$\ell_1$ a slight advantage by simultaneously considering all solutions $\lambda \in [\lambda_0
\wedge \wh{\lambda}, \lambda_0 \vee \wh{\lambda}]$ prior to thresholding,
where $\lambda_0 = 2 \sigma \sqrt{2\log(p)/n}$ equals the choice of the regularization parameter advocated in
\cite{Bickel2009} to achieve the optimal rate for the estimation of
$\beta^{\ast}$ in the $\ell_2$-norm and $\wh{\lambda} = 2 \nnorm{X^{\T} \eps/n}_{\infty}$ can be interpreted as
empirical counterpart to $\lambda_0$. The non-negative lasso modification of LARS
\cite{Efron2004} is used to obtain the solutions
$\{\wh{\beta}(\lambda): \lambda \in [\lambda_0 \wedge \wh{\lambda}, \lambda_0 \vee \wh{\lambda}] \}$; we
then report success of tNN$\ell_1$ whenever $\min_{j \in S} \wh{\beta}_j(\lambda) > \max_{j \in S^c}
\wh{\beta}_j(\lambda)$ holds for at least \emph{one} of these solutions. We
point out that specification of $\lambda_0$ is based on knowledge of the noise
variance, which constitutes a second potential advantage for tNN$\ell_1$.\\
Under the conditions of Theorem \ref{theo:plasso}, NN$\ell_1$ recovers the support
directly without thresholding. In order to judge the usefulness 
of subsequent thresholding of NN$\ell_1$, we obtain as well the set of
non-negative lasso solutions  $\{\wh{\beta}(\lambda): \lambda \geq \lambda_0
\wedge \wh{\lambda}  \}$  and check whether the sparsity pattern of any of these
solutions recovers $S$.\\
Given its simplicity, OMP serves  as basic reference
method. Success is reported whenever the
support has been recovered after $s$ steps.
\begin{figure}[h!] 
\begin{center}
\begin{tabular}{lll}
\hspace{-0.4cm}\includegraphics[height = 0.15\textheight]{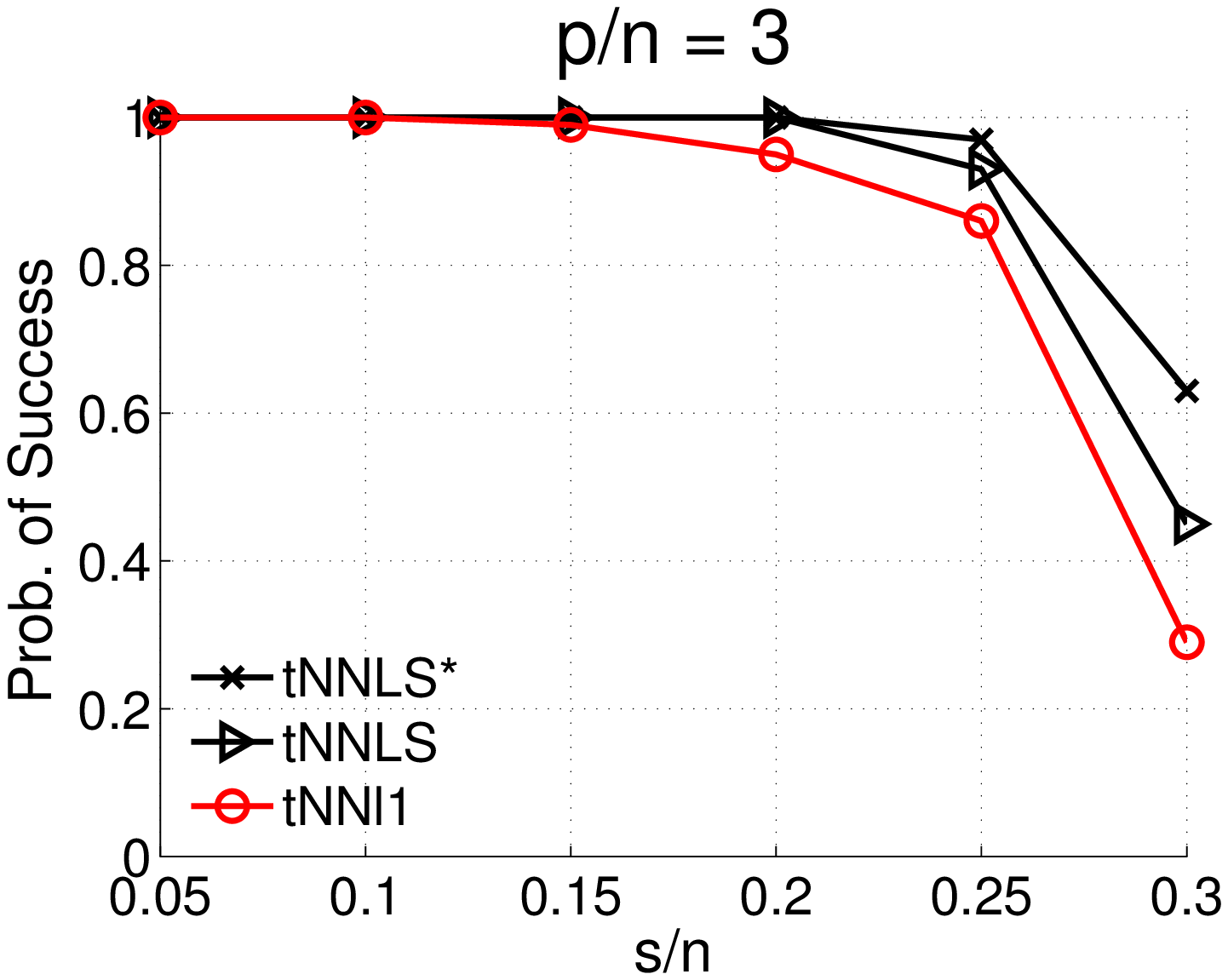}
& \hspace{-0.4cm}\includegraphics[height = 0.15\textheight]{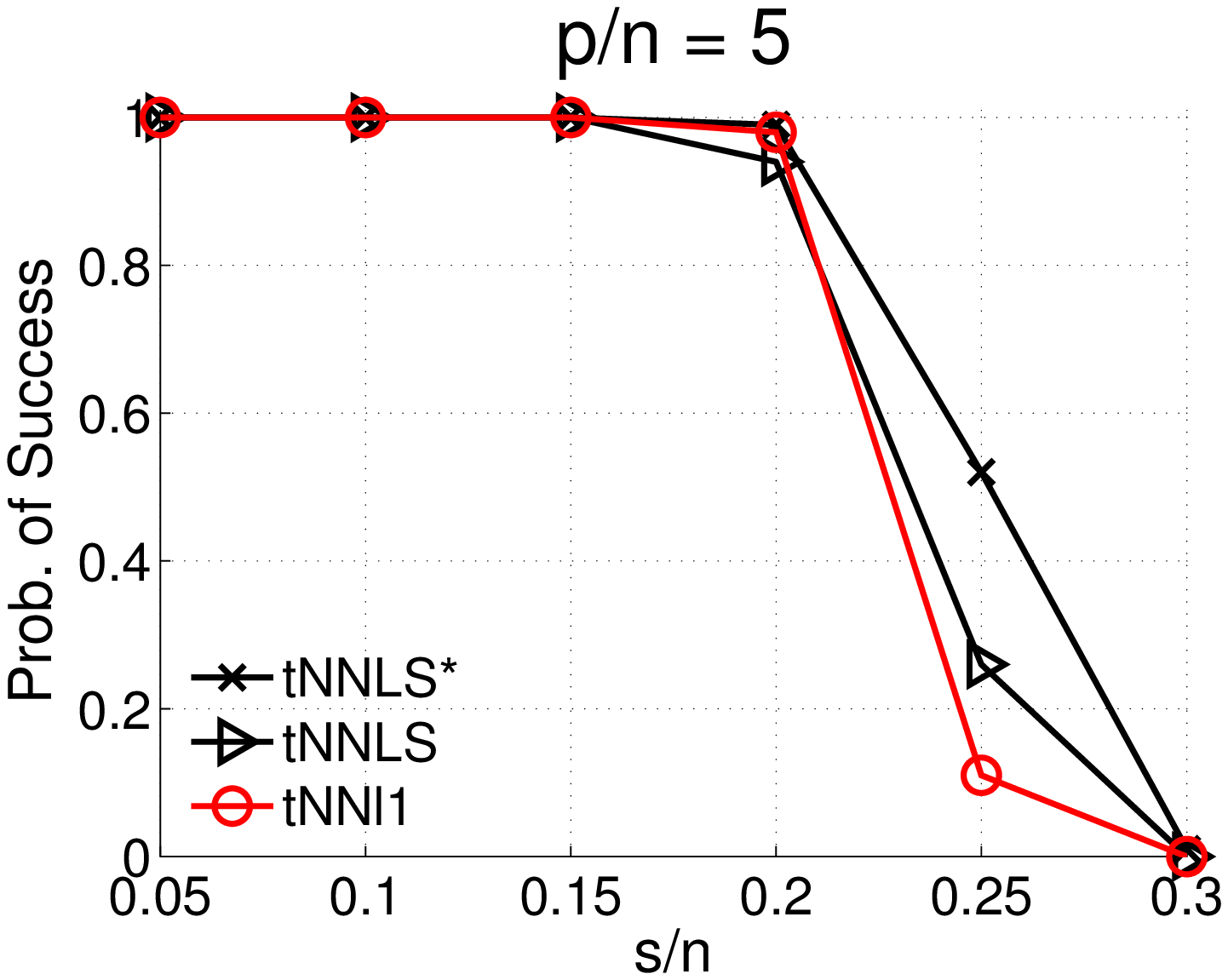}
& \hspace{-0.35cm}\includegraphics[height =
0.15\textheight]{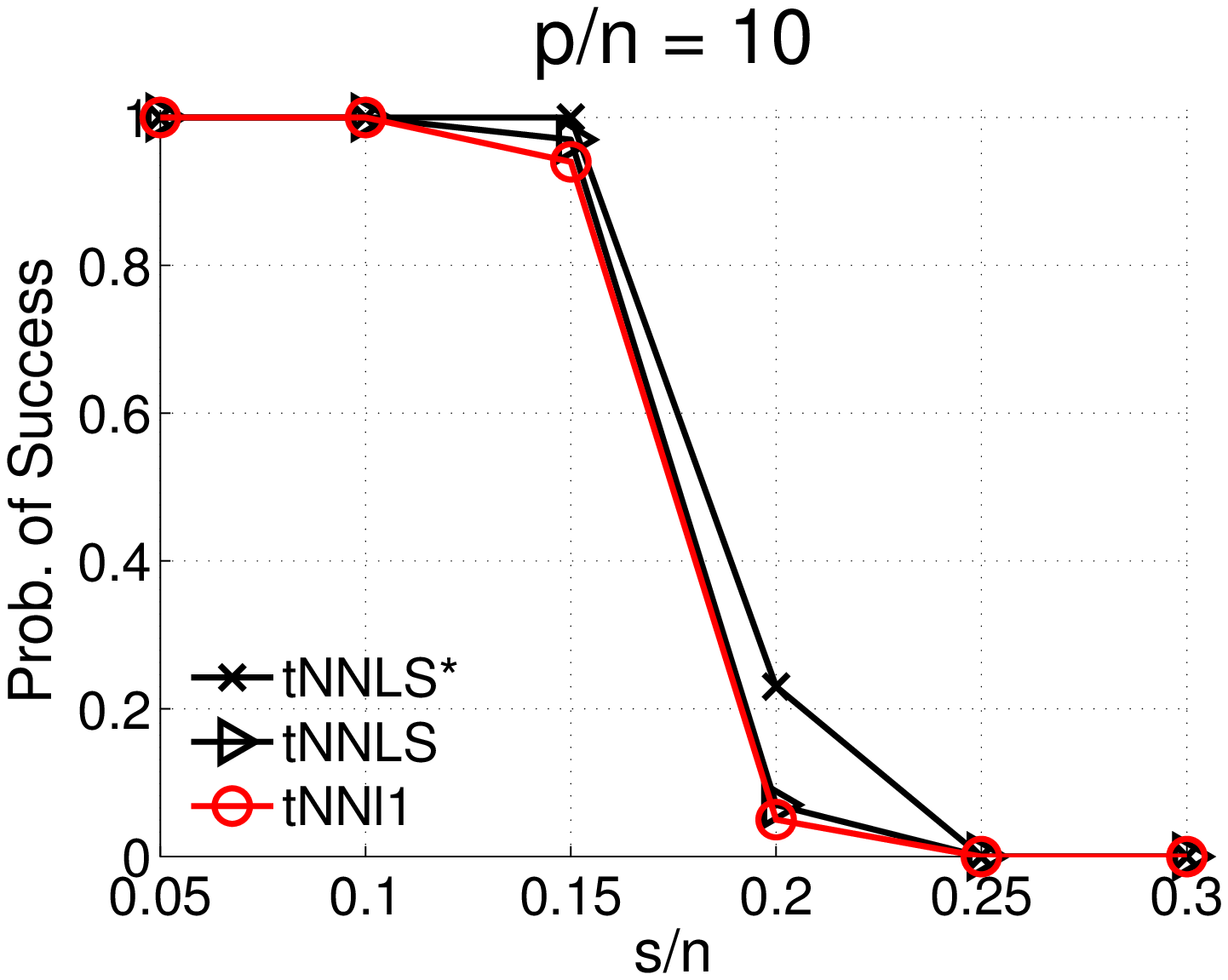}\\
\hspace{-0.4cm}\includegraphics[height = 0.15\textheight]{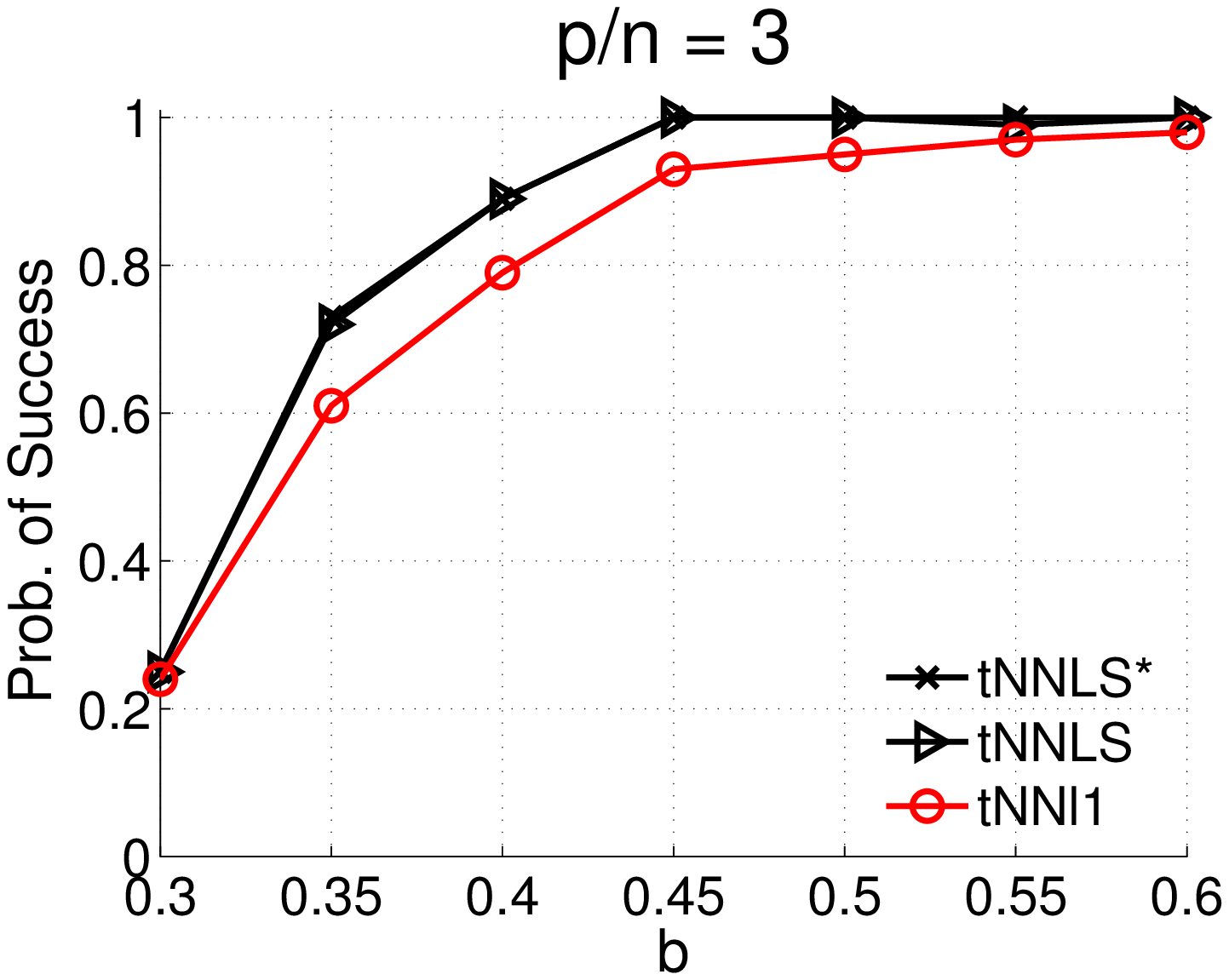}
& \hspace{-0.4cm}\includegraphics[height = 0.15\textheight]{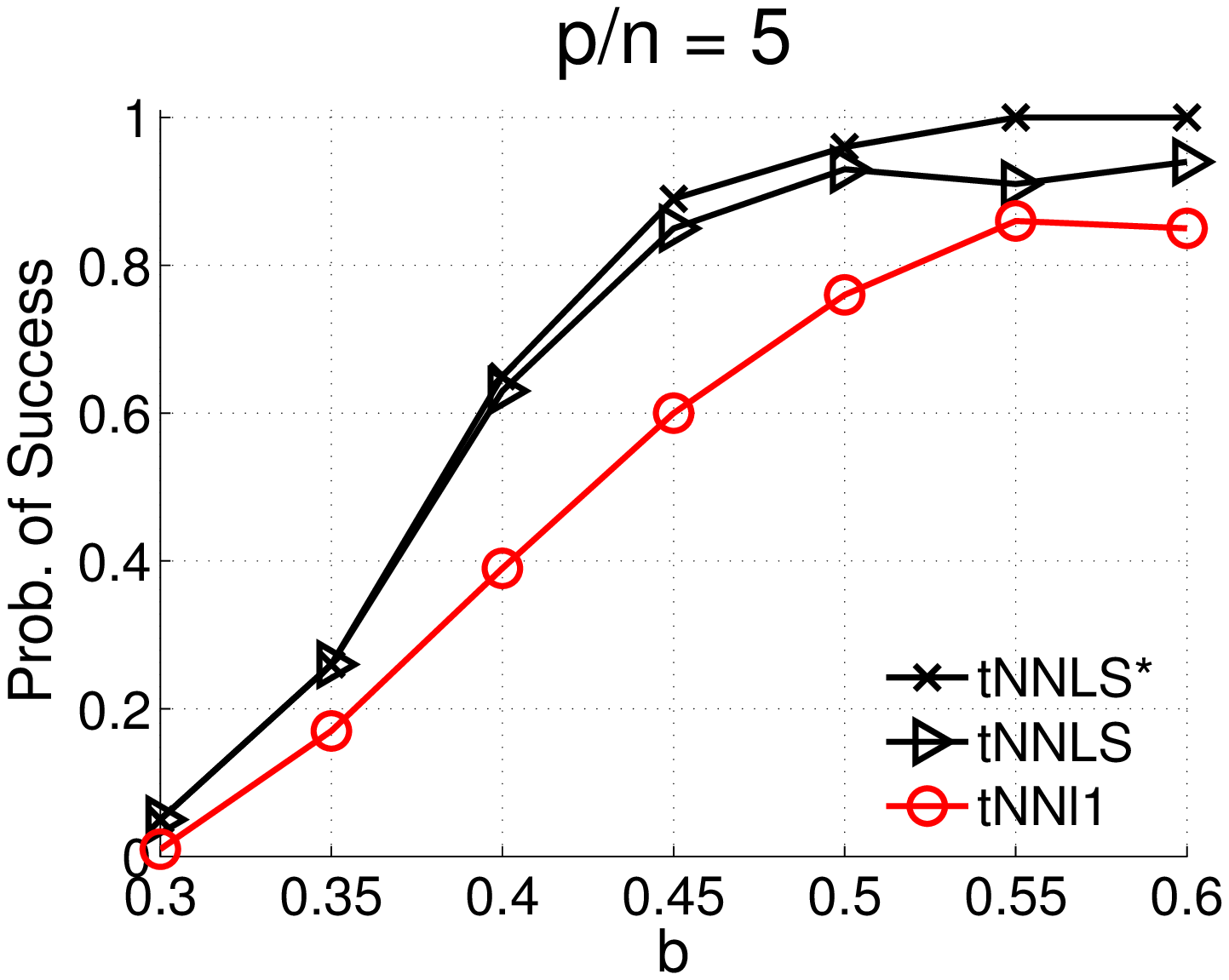}
& \hspace{-0.35cm}\includegraphics[height = 0.15\textheight]{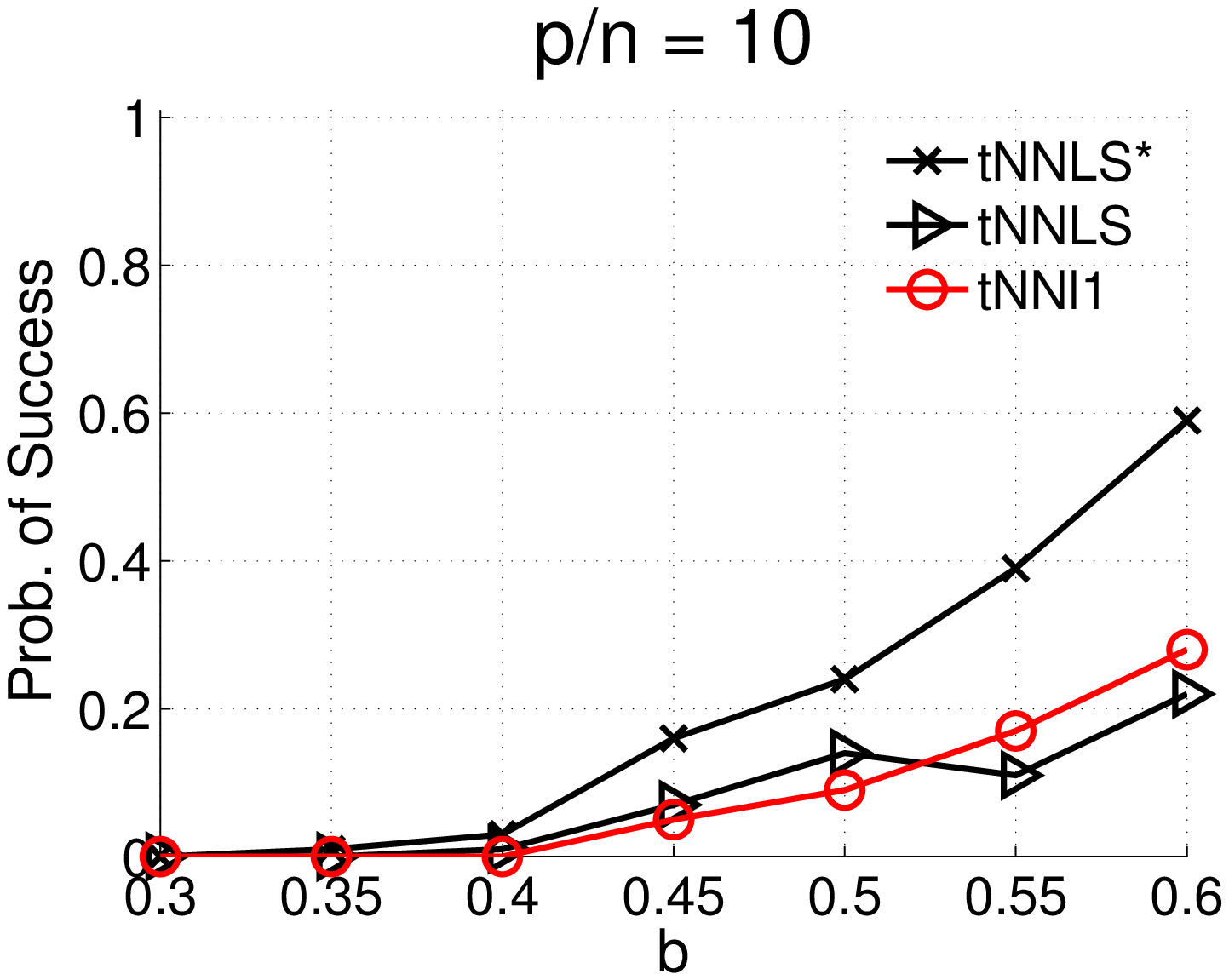}\\
\end{tabular}
\caption{Sparse recovery results for Design I. Top: Results of the set
  of experiments with fixed
  signal strength $b = 0.5$. Bottom: Results of the set of experiments
  with fixed fraction of sparsity $s/n = 0.2$. 'tNNLS$^*$' and 'tNN$\ell_1$'denote thresholded
NNLS and the thresholded non-negative lasso, where thresholding is done
with knowledge of $S$. 'tNNLS' denotes thresholded NNLS with data-driven choice
of the threshold. The results of the non-negative lasso without thresholding
and OMP are not displayed, because these two approaches fail in all instances.}\label{fig:sparserecovery_uniform}
\end{center}
\end{figure}
\begin{figure}[h!]
\begin{center}
\begin{tabular}{lll}
\hspace{-0.4cm}\includegraphics[height = 0.15\textheight]{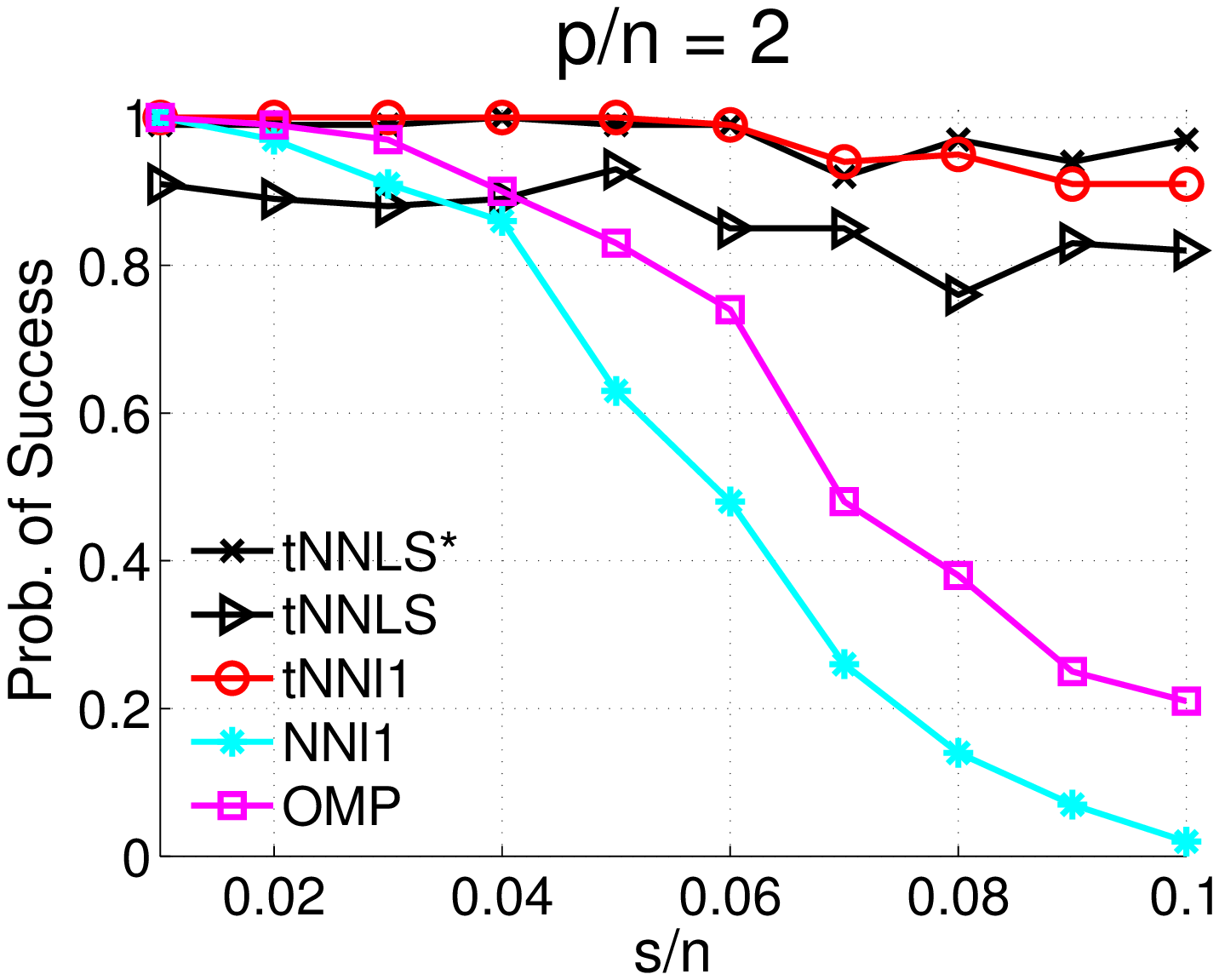}
& \hspace{-0.4cm}\includegraphics[height = 0.15\textheight]{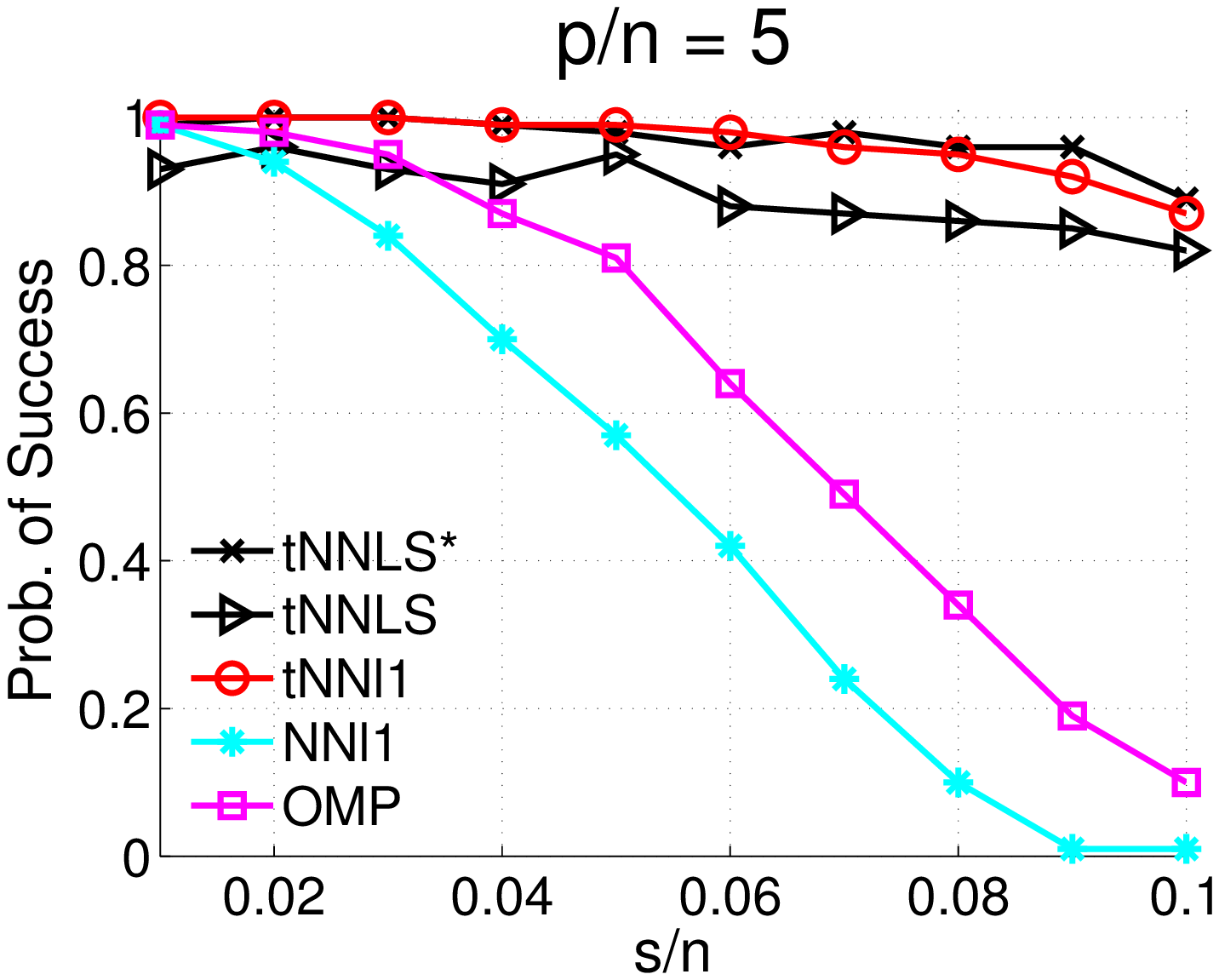}
& \hspace{-0.35cm}\includegraphics[height =
0.15\textheight]{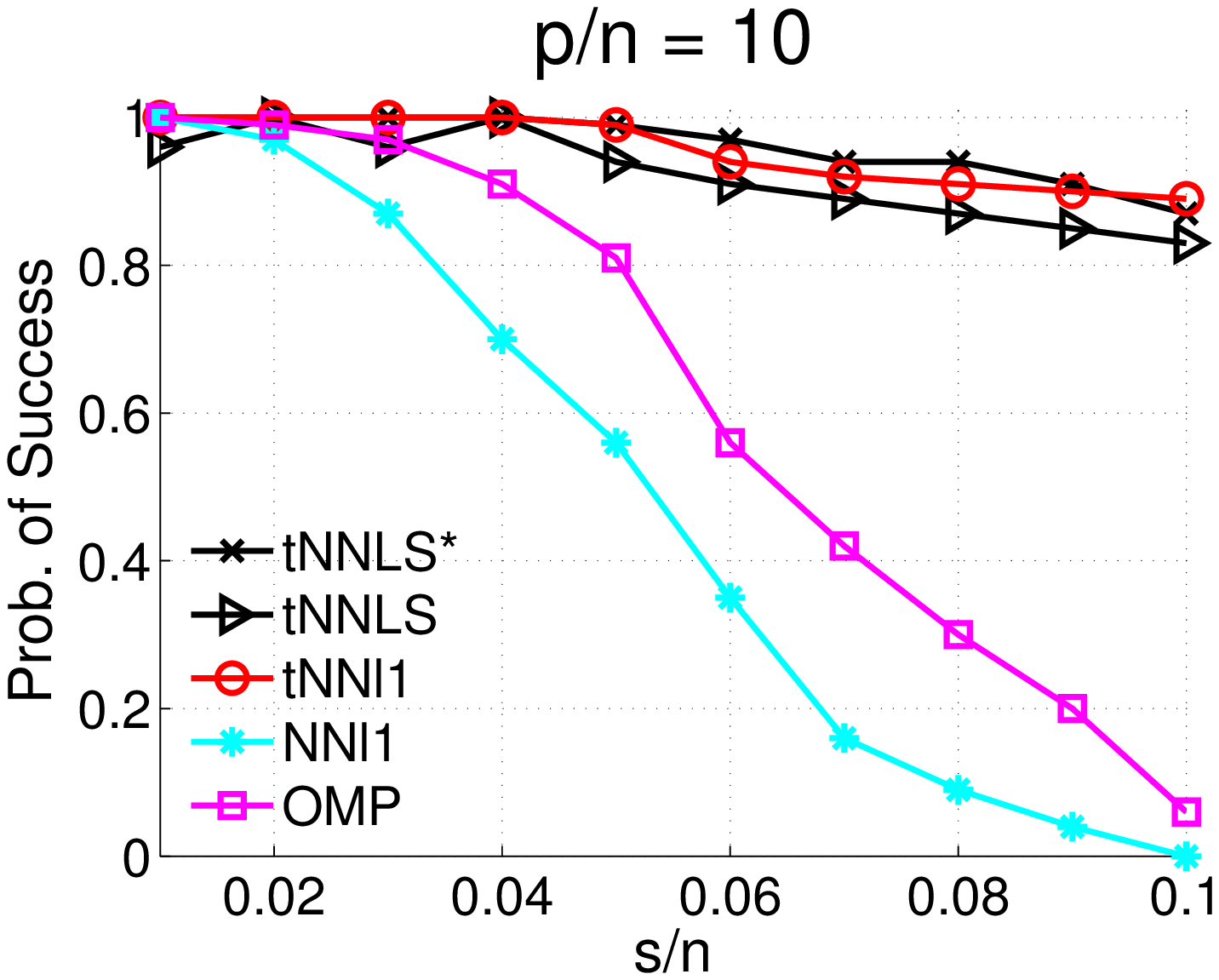}\\
\hspace{-0.4cm}\includegraphics[height = 0.15\textheight]{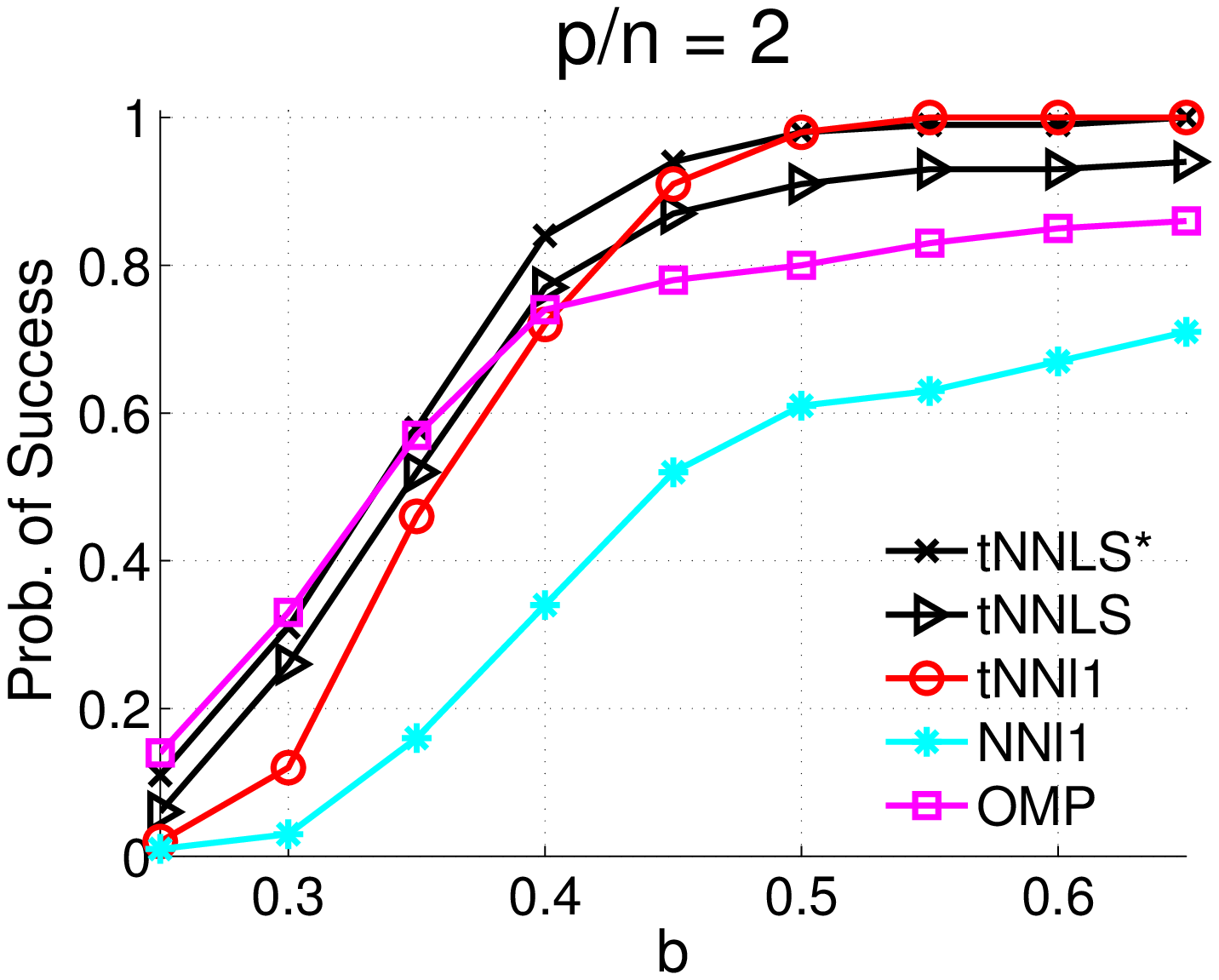}
& \hspace{-0.4cm}\includegraphics[height = 0.15\textheight]{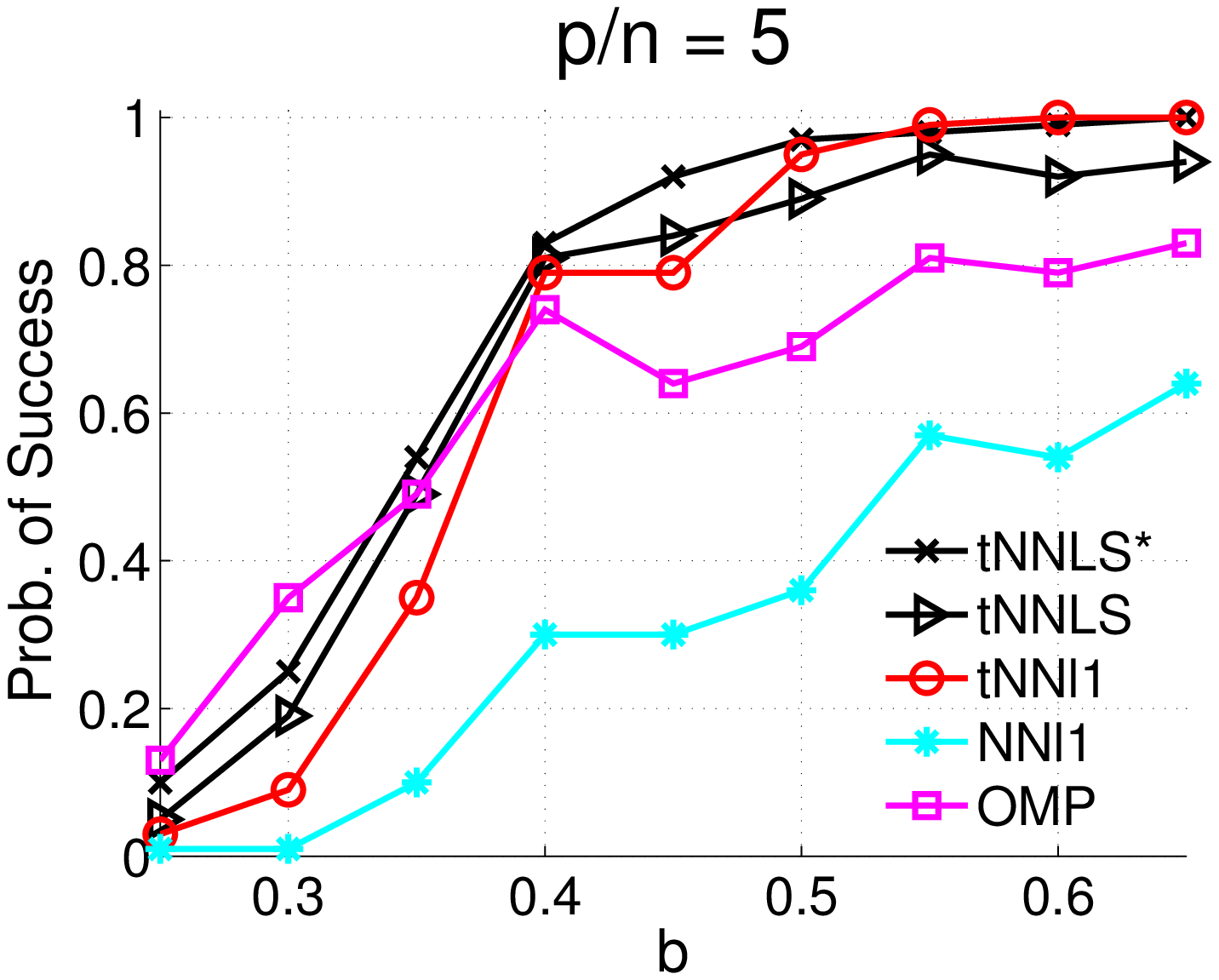}
& \hspace{-0.35cm}\includegraphics[height = 0.15\textheight]{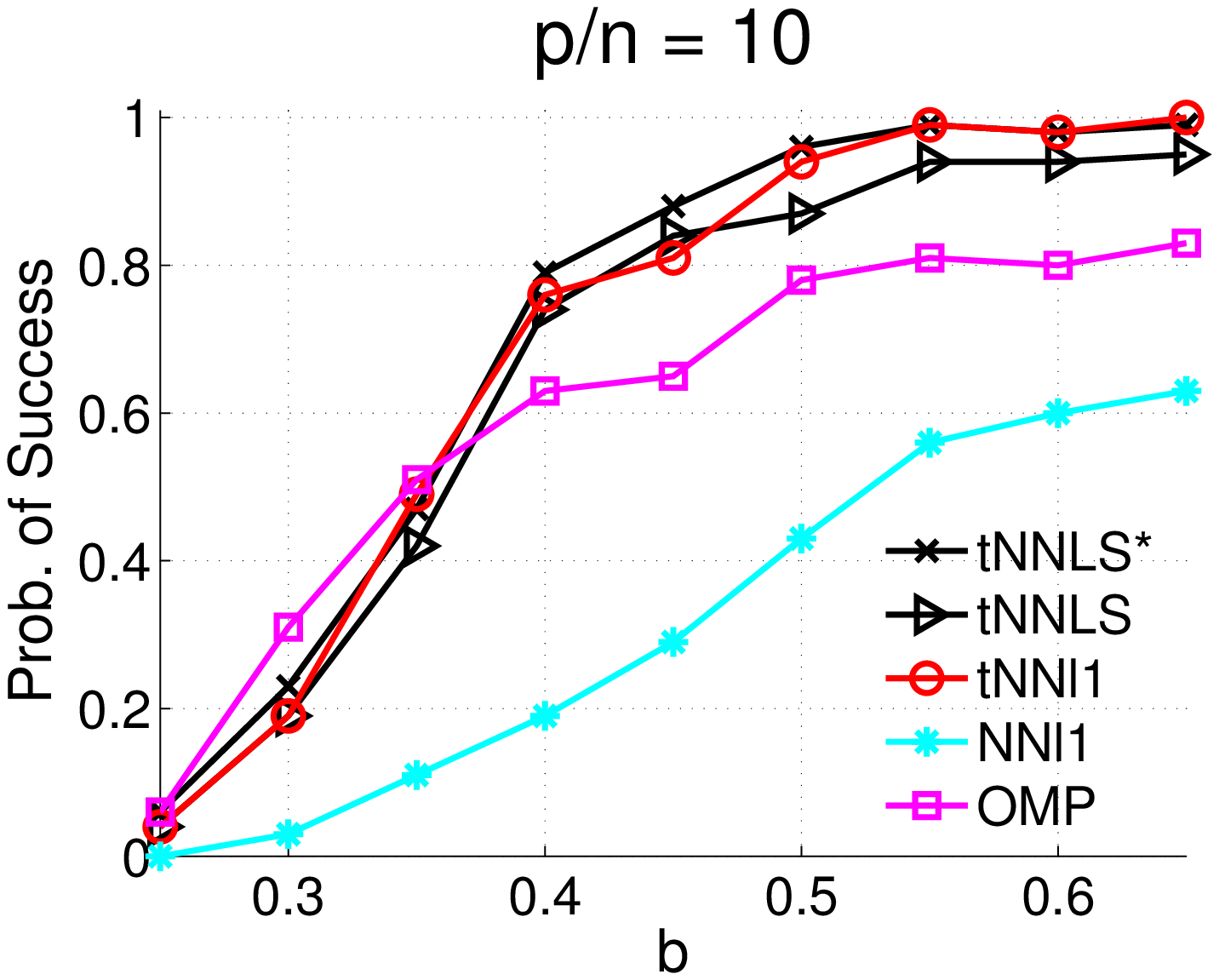}\\
\end{tabular}
\end{center}
\caption{Sparse recovery results for Design II. Top: Results of the set
  of experiments with fixed
  signal strength $b = 0.55$. Bottom: Results of the set of experiments
  with fixed fraction of sparsity $s/n = 0.05$. 'tNNLS$^*$' and 'tNN$\ell_1$'denote thresholded
NNLS and the thresholded non-negative lasso, where thresholding is done
with knowledge of $S$. 'tNNLS' denotes thresholded NNLS with data-driven choice
of the threshold. 'NN$\ell_1$' denotes the non-negative lasso without
thresholding and 'OMP' orthogonal matching pursuit.}\label{fig:sparserecovery_block}
\end{figure}
\\
\textbf{Discussion of the results.} In summary, Figures
\ref{fig:sparserecovery_uniform} and \ref{fig:sparserecovery_block} indicate that
for the two setups under consideration, NNLS and its thresholded version
exhibit excellent performance in sparse recovery. A superior performance relative to the thresholded non-negative lasso is
achieved particularly in more difficult parameter regimes characterized by
comparatively small signal strength $b$ or high fraction of sparsity. The
results of the experiments reveal that the non-negative lasso without
thresholding may perform well in estimation, but it is not competitive as far as sparse recovery is concerned. This
observation is in agreement with existing literature in which the
restrictiveness of the conditions for the lasso to select the correct set of
variables is pointed out and two stage procedures like thresholding are proposed     
as remedy \cite{Zhou2006, MeinshausenYu2009, Zhang2009b, Zhou2009, Tsybakov2010}. At this
point, we stress again that NNLS only requires one parameter
(the threshold) to be set, whereas competitive performance with regard to
sparse recovery based on the non-negative lasso entails specification of two
parameters. Let us now give some more specific comments separately for the two
designs. For Design I, thresholded NNLS visibly improves over tNN$\ell_1$,
predominantly even in case that the threshold is chosen adaptively without knowledge of $S$. For Design
II, noticeable differences between tNNLS$^*$ and tNN$\ell_1$ occur for small
values of $b$. Apart from that, the performance is essentially identical.  
Even though the results of tNNLS remain competitive, they fall behind
those of tNNLS$^*$ and tNN$\ell_1$. OMP partially        
keeps up with the other methods for $s/n$ and/or $b$ small, while NN$\ell_1$ succeeds as well in 
a substantial fraction of cases for small $s/n$. This is to be contrasted with the situation for Design I, in which both
OMP and NN$\ell_1$ do not even achieve success in a single trial. This outcome
is a consequence of the fact that the non-negative irrepresentable condition
(cf.~Section \ref{subsec:plasso}), which is necessary for
the success of OMP as well \cite{Zhang2009}, fails to hold in all these runs. 
The $\ell_{\infty}$-errors in estimating $\beta^*$ reported in Tables
\ref{tab:ellinf_uniform} and  \ref{tab:ellinf_block} are in accordance with the sparse recovery
results. The smaller $s/n$ and $p/n$, the closer NNLS and NN$\ell_1$
are in performance. An advantage of NNLS arises for more extreme combinations
of $(s/n, p/n)$. A similar conclusion can be drawn for the $\ell_2$-errors
(cf.~supplement). 
\begin{table}[h!]
  {\normalsize
\begin{tabular}{l|l|l||l|l||}
\cline{2-5}         
& \multicolumn{4}{c|}{$p/n$} \\
\cline{2-5}
& \multicolumn{2}{c||}{$2$} & \multicolumn{2}{c||}{$3$} \\
\cline{1-5}
\multicolumn{1}{|l|}{$s/n$} & \multicolumn{1}{|l|}{nnls}& nn${\ell_1}$ & \multicolumn{1}{|l|}{nnls}  &
nn${\ell_1}$ \\
\cline{1-5}
\multicolumn{1}{|l|}{\multirow{1}{*}{$0.05$}} &  .34\se{.005}& .34\se{.005} & .35\se{.005} & .36\se{.005} \\
\cline{1-5}
\multicolumn{1}{|l|}{\multirow{1}{*}{$0.1$}} & .37\se{.005}  & .37\se{.005} & .41\se{.005} & .40\se{.005}   \\
\cline{1-5}
\multicolumn{1}{|l|}{\multirow{1}{*}{$0.15$}} &  .41\se{.006} & .42\se{.009} & .44\se{.005} & .46\se{.012}   \\
\cline{1-5}
\multicolumn{1}{|l|}{\multirow{1}{*}{$0.2$}} & .43\se{.006} & .46\se{.012} & .50\se{.007} & .56\se{.023} \\
\cline{1-5}
\multicolumn{1}{|l|}{\multirow{1}{*}{$0.25$}} & .48\se{.006} & .54\se{.020} & .58\se{.009} & .72\se{.030}   \\
\cline{1-5}
\multicolumn{1}{|l|}{\multirow{1}{*}{$0.3$}} &  .55\se{.007} & .64\se{.027} & .70\se{.012} & 1.01\se{.04}\\
\cline{1-5}
\cline{2-5}         
& \multicolumn{4}{c|}{$p/n$} \\
\cline{2-5}
& \multicolumn{2}{c||}{$5$} & \multicolumn{2}{c||}{$10$} \\
\cline{1-5}
\multicolumn{1}{|l|}{$s/n$} & \multicolumn{1}{|l|}{nnls}& nn${\ell_1}$ & \multicolumn{1}{|l|}{nnls}  &
nn${\ell_1}$ \\
\cline{1-5}
\multicolumn{1}{|l|}{\multirow{1}{*}{$0.05$}} & .37\se{.005} & .38\se{.005} & .43\se{.006}& .43\se{.006}   \\
\cline{1-5}
\multicolumn{1}{|l|}{\multirow{1}{*}{$0.1$}} & .44\se{.005} & .45\se{.006}  & .51\se{.007} & .52\se{.007}   \\
\cline{1-5}
\multicolumn{1}{|l|}{\multirow{1}{*}{$0.15$}} & .52\se{.007} & .54\se{.007} & .66\se{.009} & .71\se{.012}    \\
\cline{1-5}
\multicolumn{1}{|l|}{\multirow{1}{*}{$0.2$}} & .61\se{.008} & .66\se{.009} & 1.01\se{.02} & 1.28\se{.03} \\
\cline{1-5}
\multicolumn{1}{|l|}{\multirow{1}{*}{$0.25$}} & .81\se{.014}  & 1.32\se{.04} & 1.91\se{.02} & 2.17\se{.02}   \\
\cline{1-5}
\multicolumn{1}{|l|}{\multirow{1}{*}{$0.3$}} &  1.36\se{.03} & 1.90\se{.03} & 2.32\se{.02} & 2.36\se{.03}\\
\cline{1-5}
\end{tabular}
}
\caption{Average sup-norm errors ($\pm$ standard errors) of $\nnorm{\wh{\beta} - \beta^*}_{\infty}$
(NNLS) and $\nnorm{\plasso - \beta^*}_{\infty}$ (NN$\ell_1$) for Design I with
$b = 0.5$.}\label{tab:ellinf_uniform}
\end{table}
\begin{table}
{\normalsize
\begin{tabular}{l|l|l||l|l||}
\cline{2-5}         
& \multicolumn{4}{c|}{$p/n$} \\
\cline{2-5}
& \multicolumn{2}{c||}{$2$} & \multicolumn{2}{c||}{$3$} \\ 
\cline{1-5}
\multicolumn{1}{|l|}{$s/n$} & \multicolumn{1}{|l|}{nnls}& nn${\ell_1}$ & \multicolumn{1}{|l|}{nnls}  &
nn${\ell_1}$ \\
\cline{1-5}
\multicolumn{1}{|l|}{\multirow{1}{*}{$0.02$}} &  .20\se{.005}& .32\se{.005} & .21\se{.005} & .32\se{.005} \\
\cline{1-5}
\multicolumn{1}{|l|}{\multirow{1}{*}{$0.04$}} & .23\se{.004}  & .34\se{.004} & .24\se{.007} & .35\se{.006} \\
\cline{1-5}
\multicolumn{1}{|l|}{\multirow{1}{*}{$0.06$}} &  .25\se{.006} & .36\se{.005} &
.27\se{.008} & .37\se{.006} \\
\cline{1-5}
\multicolumn{1}{|l|}{\multirow{1}{*}{$0.08$}} & .28\se{.010} & .37\se{.009} & .28\se{.009} & .37\se{.006}\\
\cline{1-5}
\multicolumn{1}{|l|}{\multirow{1}{*}{$0.1$}} & .29\se{.010} & .37\se{.007} & .32\se{.012} & .39\se{.010}\\
\cline{1-5}
\cline{2-5}         
& \multicolumn{4}{c|}{$p/n$} \\
\cline{2-5}
& \multicolumn{2}{c||}{$5$} & \multicolumn{2}{c||}{$10$} \\ 
\cline{1-5}
\multicolumn{1}{|l|}{$s/n$} & \multicolumn{1}{|l|}{nnls}& nn${\ell_1}$ & \multicolumn{1}{|l|}{nnls}  &
nn${\ell_1}$ \\
\cline{1-5}
\multicolumn{1}{|l|}{\multirow{1}{*}{$0.02$}} &  .21\se{.004} & .32\se{.004} & .22\se{.005}& .33\se{.006}  \\
\cline{1-5}
\multicolumn{1}{|l|}{\multirow{1}{*}{$0.04$}} &  .23\se{.005} & .34\se{.004}  & .24\se{.005} & .35\se{.005}\\
\cline{1-5}
\multicolumn{1}{|l|}{\multirow{1}{*}{$0.06$}} &   .27\se{.005} & .36\se{.005} & .27\se{.006} & .37\se{.006}  \\
\cline{1-5}
\multicolumn{1}{|l|}{\multirow{1}{*}{$0.08$}} & .29\se{.011} & .37\se{.009} & .30\se{.009} & .37\se{.006}\\
\cline{1-5}
\multicolumn{1}{|l|}{\multirow{1}{*}{$0.1$}} & .32\se{.011}  & .40\se{.010} & .32\se{.011} & .39\se{.008} \\
\cline{1-5}
\end{tabular}
\caption{Average sup-norm errors ($\pm$ standard errors) of $\nnorm{\wh{\beta} - \beta^*}_{\infty}$
(NNLS) and $\nnorm{\plasso - \beta^*}_{\infty}$ (NN$\ell_1$) for Design II with
$b = 0.55$.}\label{tab:ellinf_block}
}
\end{table}
\section*{Acknowledgements}
 We would like to thank Mike Davies, Rolf Schneider, Jared Tanner and Roman
 Vershynin for helpful discussions. We thank the Editor and the reviewers for
 thorougly reading the manuscript and providing valuable comments and
 suggestions. We are also indebted to a reviewer of an earlier draft,
 whose comments prompted us to organize Section \ref{subsec:oracleprediction} in its present
form.

\clearpage

\appendix

\section{Standard tail bounds for sub-Gaussian random variables}\label{app:subgaussian}

A zero-mean random variable $Z$ is called sub-Gaussian if there exists 
$\sigma > 0$ (referred to as sub-Gaussian parameter) so that the
moment-generating function obeys the bound $\E[\exp(t Z)] \leq \exp(\sigma^2
t^2/2), \; \forall t \in \R$. It follows that if $Z_1,\ldots, Z_n$ are i.i.d. copies of $Z$ and $v \in \R^n$, then $\su v_i Z_i$ is sub-Gaussian with parameter $\norm{v}_2^2
\sigma^2$. We have the well-known tail bound
\begin{equation}\label{eq:subgaussiantailbound}
\p(|Z| > z) \leq 2 \exp \left(-\frac{z^2}{2 \sigma^2} \right),\quad z \geq 0.
\end{equation}
Combining the previous two facts and using a union bound with $\mathbf{Z} =
(Z_1,\ldots,Z_n)^{\T}$ it follows that for any collection of vectors $v_j \in \R^n$, $j=1,\ldots,p$, 
\begin{equation}\label{eq:subgaussianmaximal}
\p \left(\max_{1 \leq j \leq p} |v_j^{\T} \mathbf{Z}| > \sigma \max_{1 \leq j
    \leq p} \norm{v_j}_2 \left(\sqrt{2 \log
    p} + z \right)\right) \leq 2 \exp \left(-\frac{1}{2} z^2 \right), \; \; z
\geq 0.
\end{equation}
To obtain the main results of the paper, \eqref{eq:subgaussianmaximal} is
applied with $Z = \eps$ and $v_j = X_j / n$ under the assumption $\nnorm{X_j}_2^2 = n$, $j=1,\ldots,p$, $j = 1,\ldots,p$, and the choice
$z = M \sqrt{2 \log p}$ for $M \geq 0$, which yields
\begin{equation}\label{eq:subgaussianmaximal_a}
\p \left(\max_{1 \leq j \leq p} \left| \frac{X_j^{\T} \eps}{n} \right| > \sigma (1 + M)\sqrt{\frac{2 \log
    p}{n}} \right) \leq 2 p^{-M^2}.
\end{equation}

\section{Proof of Proposition \ref{prop:selfreg}}

Since $X$ satisfies Condition \ref{cond:1}, by \eqref{eq:tau0}, there exists a
unit vector $w$ so that 
\begin{equation}\label{eq:defh}
\frac{X^{\T} w}{\sqrt{n}} = h, \quad \text{where} \; \; h \gec \tau \bm{1},
\end{equation}
for some constant $\tau > 0$. Setting $\Pi = I - ww^{\T}$ as the projection on
the subspace orthogonal to $w$, the least squares objective can be decomposed as
follows.
\begin{align*}
\frac{1}{n} \nnorm{\eps - X \beta}_2^2 &= \frac{\eps^{\T} \eps}{n} - \frac{2
  \eps^{\T} X \beta}{n} +  \frac{\beta^{\T} X^{\T} X \beta}{n} \\
&= \left( \frac{\eps^{\T} \eps}{n} - \frac{2  \eps^{\T} \Pi X \beta}{n}  +
  \frac{\beta^{\T} X^{\T} \Pi X \beta}{n}  \right) +  \frac{\beta^{\T} X^{\T}
  w w^{\T} X \beta}{n}- \\
&\quad - \frac{2 \eps^{\T} w
  w^{\T} X\beta}{n} \\
&=\frac{1}{n} \nnorm{\eps - \Pi X \beta}_2^2 + (h^{\T} \beta)^2  -  \frac{2
  \eps^{\T} w}{\sqrt{n}} h^{\T} \beta \\
&= \frac{1}{n} \nnorm{\eps - \overline{X} \beta}_2^2 + (h^{\T} \beta)^2 + O_{\p}
\left(\frac{1}{\sqrt{n}} \right)
h^{\T} \beta
\end{align*}
where $\overline{X} = \Pi X$. In the last line, we have invoked the assumptions made
for $\eps$. Writing $H$ for the diagonal matrix with the entries of $h/\tau$ on its diagonal and setting $D = H^{-1}$ and
$\wt{X} = \overline{X} D = (X \Pi) D$, we have
\begin{align*}
&\min_{\beta \gec 0} \frac{1}{n} \nnorm{\eps - \overline{X} \beta}_2^2 + (h^{\T} \beta)^2 + O_{\p}
\left(\frac{1}{\sqrt{n}} \right) h^{\T} \beta \\
= &\min_{\beta \gec 0} \frac{1}{n} \nnorm{\eps
  - \wt{X} \beta}_2^2 + \tau^2 (\bm{1}^{\T} \beta)^2 + O_{\p}
\left(\frac{1}{\sqrt{n}} \right) \tau \bm{1}^{\T} \beta,
\end{align*}
where we have used \eqref{eq:defh}. Note that by \eqref{eq:defh} and $\tau
\leq 1$, $D$ has the property claimed in the statement. In view of the
presence of the term $\tau^2 (\bm{1}^{\T} \beta)^2$, any minimizer $\beta^{\circ}$
of the r.h.s.~must obey $\bm{1}^{\T} \beta^{\circ} = O_{\p}(1)$. As a result, 
\begin{align*}
&\min_{\beta \gec 0} \frac{1}{n} \nnorm{\eps
  - \wt{X} \beta}_2^2 + \tau^2 (\bm{1}^{\T} \beta)^2 + O_{\p}
\left(\frac{1}{\sqrt{n}} \right) \tau \bm{1}^{\T} \beta \\
= &\min_{\beta \gec 0} \frac{1}{n} \nnorm{\eps
  - \wt{X} \beta}_2^2 + \tau^2 (\bm{1}^{\T} \beta)^2 + O_{\p}
\left(\frac{1}{\sqrt{n}} \right),
\end{align*}
which finishes the proof of the first claim of the proposition. To establish
the second claim, observe that any $\wh{\beta} \in \argmin_{\beta \gec 0} \frac{1}{n}
\nnorm{y - X \beta}_2^2$ satisfies 
\begin{equation*}
\frac{1}{n} \nnorm{\eps - X \wh{\beta}}_2^2 \leq \frac{1}{n} \nnorm{\eps}_2^2.
\end{equation*}
Expanding the square and re-arranging, we obtain
\begin{equation*}
\frac{1}{n} \nnorm{X \wh{\beta}}_2^2 \leq \frac{2 \eps^{\T} X \wh{\beta}}{n} 
\leq 2 \frac{\nnorm{X^{\T} \eps}_{\infty}}{n} \bm{1}^{\T} \wh{\beta}.
\end{equation*}
As established above, $\bm{1}^{\T} \wh{\beta} = O_{\p}(1)$, so that 
$\frac{1}{n} \nnorm{X \wh{\beta}}_2^2 = o_{\p}(1)$ as long as
$\frac{1}{n} \nnorm{X^{\T} \eps}_{\infty} = o_{\p}(1)$. 


\section{Proof of Theorem \ref{theo:prediction}}\label{app:theorem1}

Since $\wh{\beta}$ is a minimizer of the NNLS problem \eqref{eq:nnls} and since
$\beta^*$ is a feasible solution, we have that
\begin{align}\label{eq:basicgeneral}
\begin{split}
&\frac{1}{n} \nnorm{y - X \wh{\beta}}_2^2 \leq \frac{1}{n} \nnorm{y - X
  \beta^*}_2^2 \\
&\Leftrightarrow\;\frac{1}{n}\nnorm{(f + \eps - X \beta^*) + X \beta^* - X \wh{\beta}}_2^2 \leq \frac{1}{n} \nnorm{f
+ \eps - X \beta^*}_2^2 \\
&\Rightarrow\; \frac{1}{n} \nnorm{X \beta^* - X \wh{\beta}}_2^2 + \frac{2}{n}
(f + \eps - X \beta^*)^{\T} X (\beta^* - \wh{\beta}) \leq 0 \\
&\Rightarrow\; \frac{1}{n} \nnorm{X \beta^* - X \wh{\beta}}_2^2 \leq \frac{2}{n}
(f - X \beta^*)^{\T} X (\wh{\beta} - \beta^*) + \frac{2}{n} \eps^{\T} X
(\wh{\beta} - \beta^*).
\end{split}
\end{align}
Write $\wh{\delta} = \wh{\beta} - \beta^{\ast}$, $P = \{j: \; \wh{\delta}_j
\geq 0 \}$ and $N = \{j: \; \wh{\delta}_j < 0 \}$. 
We now lower bound $\frac{1}{n} \nnorm{X \wh{\delta}}_2^2 = \wh{\delta}^{\T}
\Sigma \wh{\delta}$ using the self-regularizing property of Condition
\ref{cond:1} according to \eqref{eq:tau0sq}. 
\begin{align}\label{eq:lowerboundselfreg}
\begin{split}
\frac{1}{n} \nnorm{X \wh{\delta}}_2^2 &= \wh{\delta}_P^{\T} \Sigma_{PP}
\wh{\delta}_P + 2 \wh{\delta}_P^{\T} \Sigma_{PN} \wh{\delta}_N +
\wh{\delta}_N^{\T} \Sigma_{NN} \wh{\delta}_N \\
&\geq \tau^2 (\bm{1}^{\T} \wh{\delta}_P)^2 - 2 \nnorm{\wh{\delta}_{P}}_1
\nnorm{\wh{\delta}_{N}}_1. 
\end{split}
\end{align}
Second, we bound the r.h.s.~of \eqref{eq:basicgeneral}. We set $A = \max_{1 \leq j \leq p} \left| \frac{1}{n} X_j^{\T} \eps
\right|$ and use the bound
\begin{equation*}
\max_{1 \leq j \leq p} \left| \frac{1}{n} X_j^{\T} (f - X \beta^{\ast})
\right| \leq \max_{1 \leq j \leq p} \frac{1}{\sqrt{n}} \norm{X_j}_2
\sqrt{\frac{1}{n} \norm{f - X \beta^*}_2^2} = \sqrt{\mc{E}^*},
\end{equation*}
obtaining that 
\begin{equation}\label{eq:upperboundbasicgeneral}
\frac{1}{n} \nnorm{X \wh{\delta}}_2^2 \leq 2 (A + \sqrt{\mc{E}^*}) \nnorm{\wh{\delta}}_1
\end{equation}
Inserting the lower bound \eqref{eq:lowerboundselfreg} into \eqref{eq:upperboundbasicgeneral}, we obtain
\begin{equation}\label{eq:combinewbasic}
  \tau^2 \nnorm{\wh{\delta}_P}_1^2 - 2 \nnorm{\wh{\delta}_{P}}_1
\nnorm{\wh{\delta}_{N}}_1  \leq 2 (A + \sqrt{\mc{E}^*})(\nnorm{\wh{\delta}_P}_1 + \nnorm{\wh{\delta}_N}_1). 
\end{equation}
We may assume that $\wh{\delta}_P \neq 0$, otherwise the assertion of the theorem
would follow immediately, because $\nnorm{\wh{\delta}_N}_1$ is already bounded
for feasibility reasons, see below. Dividing both sides by $\nnorm{\wh{\delta}_P}_1$ and
re-arranging yields
\begin{equation}\label{eq:dPbound}
  \nnorm{\wh{\delta}_P}_1 \leq \frac{4 (A  + \sqrt{\mc{E}^{\ast}}) + 2 \nnorm{\wh{\delta}_{N}}_1}{\tau^2},
\end{equation}
where we have assumed that $\nnorm{\wh{\delta}_N}_1 \leq
\nnorm{\wh{\delta}_{P}}_1$ (if that were not the case, one would obtain
$\nnorm{\wh{\delta}_P}_1 \leq \nnorm{\wh{\delta}_N}_1$, which is stronger than
\eqref{eq:dPbound}, since $0 < \tau^2 \leq  1$). We now substitute
\eqref{eq:dPbound} back into \eqref{eq:basicgeneral} and add $\mc{E}^{\ast} =
\frac{1}{n} \nnorm{X \beta^{\ast} - f}_2^2$ to both sides of the inequality in order to obtain  
\begin{align*}
\wh{\mc{E}} = \frac{1}{n} \nnorm{X \wh{\beta} - f}_2^2  &\leq \mc{E}^* + 2 A (\nnorm{\wh{\delta}_P}_1 + \nnorm{\wh{\delta}_N}_1) \\
&\leq \mc{E}^* + 2 A \left(\frac{4 (A + \sqrt{\mc{E}^{*}}) + 2 \nnorm{\wh{\delta}_{N}}_1}{\tau^2} +
  \nnorm{\wh{\delta}_N}_1 \right) \\
&\leq \mc{E}^* + \frac{6 A \nnorm{\beta^{\ast}}_1 + 8 (A^2 + A \sqrt{\mc{E}^*})}{\tau^2},  
\end{align*}
noting that by feasibility of $\wh{\beta}$, one has $\wh{\delta} \gec -\beta^{\ast}$ and hence
$\nnorm{\wh{\delta}_N}_1 \leq \nnorm{\beta^{\ast}}_1$. 
Using \eqref{eq:subgaussianmaximal_a}, the event $ \left \{A \leq (1 + M)\sigma
\sqrt{\frac{2 \log p}{n}} \right \}$ holds with probability no less than $1 -
2 p^{-M^2}$. The result follows.

\section{Proofs of Lemma  \ref{lem:lassocone} and Theorem \ref{theo:oracle}}
We build on ideas already used in the proof of Theorem
\ref{theo:prediction}. In particular, all notations introduced in the previous
proof are adopted. First note that $S^c \subseteq P$ and $N \subseteq S$. Hence, we obtain the
following analog to \eqref{eq:combinewbasic}.  
\begin{equation*}
  \tau^2 \nnorm{\wh{\delta}_{S^c}}_1^2 - 2 \nnorm{\wh{\delta}_{S^c}}_1
\nnorm{\wh{\delta}_{S}}_1  \leq 2 A (\nnorm{\wh{\delta}_S}_1 + \nnorm{\wh{\delta}_{S^c}}_1).
\end{equation*} 
Dividing both sides by $\nnorm{\wh{\delta}_{S^c}}_1$, assuming that $0 <
\nnorm{\wh{\delta}_S}_1 \leq \nnorm{\wh{\delta}_{S^c}}_1$ (otherwise, the claim 
$\wh{\delta} \in \mc{R}(3/\tau^2, S)$ as in the first event of Lemma
\ref{lem:lassocone} would follow trivially), we obtain
\begin{equation*}
  \tau^2 \nnorm{\wh{\delta}_{S^c}}_1 \leq 4 A  + 2  \nnorm{\wh{\delta}_{S}}_1.
\end{equation*}
If $4A \leq  \nnorm{\wh{\delta}_{S}}_1$, then the first event of Lemma
\ref{lem:lassocone} occurs. Otherwise, we conclude that the second event of
Lemma \ref{lem:lassocone} occurs by applying \eqref{eq:subgaussianmaximal_a}
to bound $A$ as in the proof of Theorem \ref{theo:prediction}.\\
In the latter case, the assertion of Theorem \ref{theo:oracle} follows
immediately. Thus, the rest of the proof is conditional on the first event.
In terms of Condition \ref{cond:2}, $\wh{\delta} \in \mc{R}(3/\tau^2, S)$, so
that one may invoke the restricted eigenvalue condition \eqref{eq:recondition},
which, when applied to \eqref{eq:basicgeneral}, yields
\begin{equation*} 
\phi \left(\frac{3}{\tau^2}, s \right) \nnorm{\wh{\delta}_S}_2^2 \leq \frac{1}{n} \nnorm{X \wh{\delta}}_2^2
\leq 2 A (\nnorm{\wh{\delta}_S}_1 + \nnorm{\wh{\delta}_{S^c}}_1) \leq 2 \left(1 +
\frac{3}{\tau^2} \right) A \nnorm{\wh{\delta}_S}_1
\end{equation*}
which implies that
\begin{equation*} 
\nnorm{\wh{\delta}_S}_1 \leq \frac{2 s}{\phi \left(\frac{3}{\tau^2}, s \right)} \left(1 +
\frac{3}{\tau^2} \right) A \; \; \Longrightarrow \; \nnorm{\wh{\delta}}_1 \leq 
\frac{2 s}{\phi \left(\frac{3}{\tau^2}, s \right)} \left(1 +
\frac{3}{\tau^2} \right)^2 A  
\end{equation*}
The preceding bound in turn implies  
\begin{equation*}
\frac{1}{n} \nnorm{X \wh{\delta}}_2^2 \leq \frac{4 s}{\phi \left(\frac{3}{\tau^2}, s \right)} \left(1 +
\frac{3}{\tau^2} \right)^2 A^2.   
\end{equation*}
Controlling $A$ as above, the $\ell_1$-bound and the bound on the prediction
error follow.  The bound on $\nnorm{\wh{\delta}}_q^q$ for $q \in (1,2]$ can be established
similarly. The proof is along the lines of the proof of Theorem 7.2 in 
\cite{Bickel2009} and is hence omitted.  

\section{Proof of Lemma \ref{lem:decoupling}}\label{app:lem:decoupling}

The proof of the lemma relies on the following auxiliary result, which
is immediate from the KKT optimality conditions of the NNLS problem. Its proof is hence omitted.   
\begin{lemma}\label{lem:NNLSopt} $\wh{\beta}$ is a minimizer of the NNLS
  problem \eqref{eq:nnls} if and only if
  there exists $F \subseteq \{1,\ldots,p \}$ such that
\begin{equation*}
\frac{1}{n} X_j^{\T} (y - X \wh{\beta}) = 0, \; \text{and} \; \wh{\beta}_j > 0, \; j \in F, \quad \frac{1}{n}
X_j^{\T} (y - X \wh{\beta}) \leq 0, \;  \text{and} \; \wh{\beta}_j = 0, \; j \in F^c. 
\end{equation*}
\end{lemma}
Lemma \ref{lem:NNLSopt} implies that any NNLS solution is a minimizer of a least squares problem
subject to the equality constraint $\beta_{F^c} = 0$ given the active set $F$,
that is 
\begin{equation*}
\frac{1}{n} \nnorm{y  - X \wh{\beta}}_2^2 = \min_{\beta} \frac{1}{n} \nnorm{y  - X \beta}_2^2 \; \; \; \text{subject to} \;
\beta_{F^c} = 0.
\end{equation*} 
\emph{Proof of Lemma \ref{lem:decoupling}.}
The NNLS objective can be split into two parts as follows. 
\begin{equation}\label{eq:decoupling}
\frac{1}{n} \norm{y - X \beta}_2^2 =  \frac{1}{n}  \norm{\Pi_S y - 
X_S \beta_S - \Pi_S X_{S^c} \beta_{S^c}}_2^2 +
\frac{1}{n} \norm{\xi - Z \beta_{S^c}}_2^2,  
\end{equation}
Separate minimization of the second summand on the r.h.s. of \eqref{eq:decoupling} yields
$\wh{\beta}^{(P1)}$. Substituting $\wh{\beta}^{(P1)}$ for $\beta_{S^c}$ in the
first summand, and minimizing the latter amounts to solving $(P2)$. In view of
Lemma \ref{lem:NNLSopt}, if $\wh{\beta}^{(P2)} \succ 0$, it coincides with an unconstrained least squares estimator
corresponding to problem $(P2)$. This implies that the
optimal value of $(P2)$ must be
zero, because the observation vector $X_S \beta_S^* + \Pi_S (\eps - X_{S^c} \wh{\beta}^{(P1)})$ of the non-negative least squares problem $(P2)$ is contained in the column space of
$X_S$. Since the second summand in \eqref{eq:decoupling} corresponding to $(P1)$ cannot be made smaller
than by separate minimization, we have minimized the non-negative least
squares objective. 

\section{Proof of Theorem \ref{theo:robustsparserecovery}}\label{app:robustsparserecovery}

We here state and prove a result that is slightly more general than 
Theorem \ref{theo:robustsparserecovery}, as it covers also the case
of an approximately sparse target. Writing $\beta_{(1)}^* \geq \ldots \geq \beta_{(p)}^* \geq 0$ for the
sequence of ordered coefficients, let $S = \{j: \; \beta_j^{\ast} \geq
\beta_{(s)}^* \}$ be the set of the $s$ largest coefficients of $\beta^*$ (for
simplicity, assume that there are no ties). For the result that follows,
we think of $\nnorm{\beta_{S^c}^*}_1$ being considerably smaller than the
entries of $\beta_S^{\ast}$. 
\begin{theo}\label{theo:robustsparserecovery_a}
Consider the linear model $y = X \beta^{\ast} + \eps$, where $\beta^{\ast} \gec
0$ and $\eps$ has i.i.d. zero-mean sub-Gaussian entries with sub-Gaussian
parameter $\sigma$. For $M \geq 0$, set 
\begin{align*}
&b = \frac{2 \left(\nnorm{\beta_{S^c}^{\ast}}_1 + (1 + M) \sigma
  \sqrt{\frac{2 \log p}{n}} \right)}{\tau^2(S)},\\
\text{and}\;\;&\wt{b} = (b +
\nnorm{\beta_{S^c}^*}_1) \cdot K(S) +
\frac{(1 + M)
  \sigma}{\sqrt{\phi_{\min}(S)}} \sqrt{\frac{2 \log p}{n}}.   
\end{align*}
If $\beta_{\min}(S) > \wt{b}$, then the NNLS estimator $\wh{\beta}$ has the following properties with probability no less
than $1 - 4 p^{-M^2}$:
\begin{equation*} 
\nnorm{\wh{\beta}_{S^c}}_1 \leq b \quad \text{and} \quad \nnorm{\wh{\beta}_S -
  \beta_S^*}_{\infty} \leq \wt{b}. 
\end{equation*}
\end{theo}
The special case of exact sparsity in which $S$
equals the support of $\beta^{\ast}$ is obtained for $\nnorm{\beta_{S^c}^*}_1
= 0$ (cf.~Theorem \ref{theo:robustsparserecovery}).
\begin{bew} First note that in the more general case with $\beta_{S^c}^* \neq
  0$, an analog of Lemma \ref{lem:decoupling} holds with   
\begin{align*}
&(P1): \min_{\beta^{(P1)} \gec 0} \frac{1}{n} \nnorm{\xi + Z \beta_{S^c}^* - Z
  \beta^{(P1)}}_2^2, \\
&(P2): \; \min_{\beta^{(P2)} \gec 0} \frac{1}{n} \nnorm{\Pi_S \eps + \Pi_S
  X_{S^c} \beta_{S^c}^* + X_S \beta_S^{*} - X_S \beta^{(P2)}  - \Pi_S X_{S^c} \wh{\beta}^{(P1)}}_2^2
\end{align*}
Consider problem $(P1)$.\\
\textbf{Step 1:} \emph{Controlling} $\nnorm{\wh{\beta}^{(P1)}}_1$ \emph{via}
$\tau^2(S)$. Since $\wh{\beta}^{(P1)}$ is a minimizer and $0$ is feasible
for $(P1)$, we have
\begin{equation*}
\frac{1}{n} \nnorm{\xi + Z \beta_{S^c}^{\ast} - Z \wh{\beta}^{(P1)}}_2^2 \leq \frac{1}{n} \nnorm{\xi + Z \beta_{S^c}^{\ast}}_2^2,
\end{equation*}
which implies that
\begin{align}\label{eq:ell1boundSc_1}
\begin{split}
(\wh{\beta}^{(P1)})^{\T} \frac{1}{n} Z^{\T} Z  \wh{\beta}^{(P1)} &\leq 
    \nnorm{\wh{\beta}^{(P1)}}_1 \left( A + 2 \norm{\fon Z^{\T} Z
      \beta_{S^c}^*}_{\infty} \right), \; \; \; A = \max_{j}  \frac{2}{n}
    |Z_j^{\T} \xi|.\\ 
&\leq \nnorm{\wh{\beta}^{(P1)}}_1 \left( A + 2 \nnorm{\beta_{S^c}^*}_1
  \max_{j,k} Z_j^{\T} Z_k / n
       \right) \\
       & \leq \nnorm{\wh{\beta}^{(P1)}}_1 \left( A + 2 \max_{j,k} \nnorm{Z_j/\sqrt{n}}_2 \nnorm{Z_k/\sqrt{n}}_2
      \nnorm{\beta_{S^c}^*}_1 \right) \\
&\leq \nnorm{\wh{\beta}^{(P1)}}_1 (A + 2 \nnorm{\beta_{S^c}^{\ast}}_1).
\end{split}  
\end{align}
In the last inequality, we have used that for all $j=1,\ldots,p$, it holds that
\begin{equation}\label{eq:columnormZ}
\nnorm{Z_j}_2 = \nnorm{\Pi_S^{\perp} X_j}_2 \leq \nnorm{X_j}_2.
\end{equation}
As observed in \eqref{eq:taudual}, $\tau^2(S) = \min_{\lambda \in T^{p -s
  - 1}} \lambda^{\T} \frac{1}{n} Z^{\T} Z \lambda$, s.t. the l.h.s. of
\eqref{eq:ell1boundSc_1} can be lower bounded via
\begin{equation}\label{eq:ell1boundSc_2}
(\wh{\beta}^{(P1)})^{\T} \frac{1}{n} Z^{\T} Z \wh{\beta}^{(P1)} \geq  \tau^2(S) \nnorm{\wh{\beta}^{(P1)}}_1^2. 
\end{equation}
Combining \eqref{eq:ell1boundSc_1} and \eqref{eq:ell1boundSc_2}, we have
$\nnorm{\wh{\beta}^{(P1)}}_1 \leq (A + 2 \nnorm{\beta_{S^c}^{\ast}}_1)/\tau^2(S)$.\\
\textbf{Step 2:} \emph{Back-substitution into (P2).} Equipped with the bound
just derived, we insert $\wh{\beta}^{(P1)}$ into problem $(P2)$ of Lemma \ref{lem:decoupling},
and show that in conjunction with the assumptions made for the minimum
support coefficient $\beta_{\min}(S)$, the \emph{ordinary} least squares
estimator corresponding to ($P2$)
\begin{equation*}
\bar{\beta}^{(P2)} = \argmin_{\beta^{(P2)}} \frac{1}{n} \nnorm{\Pi_S y - 
X_S \beta^{(P2)}  - \Pi_S X_{S^c} \wh{\beta}^{(P1)}}_2^2
\end{equation*}
has only positive components. Lemma \ref{lem:decoupling} then yields  $\bar{\beta}^{(P2)} =
\wh{\beta}^{(P2)} = \wh{\beta}_S$. Using the closed form expression for the ordinary least squares estimator, one
obtains
\begin{align}\label{eq:betaSclosed}
\begin{split}
\bar{\beta}^{(P2)} &=  \frac{1}{n} \Sigma_{SS}^{-1}  X_S^{\T} \Pi_S (y - X_{S^c} (\wh{\beta}^{(P1)} - \beta_{S^c}^*)) \\
&= \frac{1}{n} \Sigma_{SS}^{-1}  X_S^{\T} (X_S
\beta_S^{\ast} + \Pi_S \eps - \Pi_S X_{S^c} (\wh{\beta}^{(P1)} - \beta_{S^c}^*)) \\
&= \beta_S^{\ast} + \frac{1}{n} \Sigma_{SS}^{-1}  X_S^{\T} \eps -
\Sigma_{SS}^{-1} \Sigma_{SS^c}  (\wh{\beta}^{(P1)} - \beta_{S^c}^{\ast}).
\end{split}
\end{align}
It remains to control the two terms $A_S = \frac{1}{n} \Sigma_{SS}^{-1} X_S^{\T} \eps
$ and $\Sigma_{SS}^{-1} \Sigma_{SS^c} (\wh{\beta}^{(P1)} -
\beta_{S^c}^*)$. For the second term, we have
\begin{align}\label{eq:backsubsbound}
\begin{split}
\nnorm{\Sigma_{SS}^{-1} \Sigma_{SS^c} (\wh{\beta}^{(P1)} - \beta_{S^c}^*)}_{\infty} &\leq
\max_{\; \norm{v}_{\infty} = 1} \nnorm{\Sigma_{SS}^{-1} v}_{\infty}
\nnorm{\Sigma_{SS^c} (\wh{\beta}^{(P1)} - \beta_{S^c}^*)}_{\infty} \\
&\overset{\eqref{eq:constantsp1}}{\leq} K(S) \; (\nnorm{\wh{\beta}^{(P1)}}_1 +
\nnorm{\beta_{S^c}^*}_1).
\end{split}
\end{align}
\textbf{Step 3:} \emph{Putting together the pieces.}
The two random terms $A$ and $A_S$ are maxima of a finite collection of
linear combinations of sub-Gaussian random variables so that \eqref{eq:subgaussianmaximal} in Appendix
\ref{app:subgaussian} can be applied by estimating Euclidean norms. For $A$, we use
\eqref{eq:columnormZ}. Second, we have 
\begin{equation}\label{eq:AS}
A_S = \max_{1 \leq j \leq s} \frac{|v_j^{\T} \eps|}{n}, \quad v_j = X_S
\Sigma_{SS}^{-1} e_j, \; \; j=1,\ldots,s,   
\end{equation}
where $e_j$ denotes the $j$-th canonical basis vector. One has
\begin{equation*}
\max_{1 \leq j \leq s} \nnorm{v_j}_2^2 = \max_{1 \leq j \leq s} e_j^{\T}
\Sigma_{SS}^{-1} X_S^{\T} X_S \Sigma_{SS}^{-1} e_j
\overset{\eqref{eq:constantsp1}}{\leq} \frac{n}{\phi_{\min}(S)}.
\end{equation*} 
It follows that for any $M \geq 0$ the event
\begin{equation*}
\left \{A \leq 2 (1 + M)\sigma \sqrt{\frac{2 \log p}{n}} \right \} \cap \left \{A_S \leq
\frac{(1 + M) \sigma}{\sqrt{\phi_{\min}(S)}} \sqrt{\frac{2 \log p}{n}} \right \} 
\end{equation*}
holds with probability no less than $1 - 4 p^{-M^2}$. Conditional on that
event, it follows that with $b$ as in Theorem \ref{theo:robustsparserecovery},
we have 
\begin{equation*}
\nnorm{\beta_S^{\ast} - \bar{\beta}^{(P2)}}_{\infty} \leq  (b +
\nnorm{\beta_{S^c}^*}_1) K(S)  + \frac{(1 + M)
  \sigma}{\sqrt{\phi_{\min}(S)}} \sqrt{\frac{2 \log p}{n}}, 
\end{equation*}
and hence, using the lower bound on $\beta_{\min}(S)$, that
$\bar{\beta}^{(P2)} = \wh{\beta}_S \succ 0$ and thus also that
$\wh{\beta}^{(P1)} = \wh{\beta}_{S^c}$.
\end{bew}

\section{Proof of Theorem \ref{theo:ellinfsmin}}\label{app:ellinfsmin}

\paragraph{Part 1: proof of uniqueness.} To prove the first part of the 
theorem asserting uniqueness of the NNLS solution, we need two additional
lemmas. The first one is a concentration result which is a special case
of Theorem 2.5 in \cite{Latala2007}.
\begin{lemma}\label{lem:concprojection} 
Let $\Pi \in \R^{n \times n}$ be a projection matrix on a
  $d$-dimensional subspace of $\R^n$ and let
$\eps = (\eps_1,\ldots,\eps_n)^{\T}$ be a random vector whose entries are i.i.d.~zero-mean sub-Gaussian
random variables with parameter $\sigma$. Then
\begin{equation*}
  \p \left(\nnorm{\Pi \eps}_2^2  \leq \E[\eps_1^2] \, \frac{d}{4} \right) \leq
  2 \exp \left(-\frac{c}{\sigma^4} d \right),
\end{equation*}
where $c > 0$ is a universal constant.  
\end{lemma}
The second lemma provides two sufficient conditions for the NNLS solution to
be unique.
\begin{lemma}\label{lem:uniqueness} Let the columns of $X$ be in 
general linear position. Then the NNLS problem
has a unique solution if one of the following holds:\\
(i) $\; p \leq n$, $\quad$(ii) $\; \min_{\beta \gec 0} \frac{1}{n} \nnorm{y -
  X \beta}_2^2 > 0$.\\ Moreover, under (ii) the active set $F =
\{j:\,\wh{\beta}_j > 0 \}$ satisfies $|F| \leq  \min \{n-1, p\}$. Conversely,
if $y$ has a distribution that is absolutely continuous with respect to the
Lebesgue measure, then $|F| \leq \min \{n-1, p\}$ implies with probability one
that the NNLS problem has a unique solution.  
\end{lemma}
\begin{bew} 
  Suppose that (i) holds. The fact that the columns of $X$
are in general linear position implies that $\Sigma$ is strictly
positive definite so that the NNLS objective is strictly convex and
hence has a unique minimizer. We now turn to the case $p > n$. We first note
that $X \wh{\beta}$ is unique, because it is the projection of $y$ onto the polyhedral cone $\mc{C} = \{z \in \R^n: z = X
\beta, \, \beta \in \R_{+}^p\}$, which is a convex set. Moreover, under (ii),
$X \wh{\beta}$ must be contained in the boundary $\partial \mc{C}$ of $\mc{C}$ (by general linear position, the interior of $\mc{C}$ is non-empty).  
Note that $\partial \mc{C}$ equals the union of the facets of $\mc{C}$, that is
\begin{align*}
\partial \mc{C} &= \bigcup_{J \in \mc{F}} \mc{C}_J,  \; \; \text{where} \; \mc{C}_{J} =
\{z \in \mc{C}: z = X_{J} \beta, \; \beta \in \R_{+}^{n-1} \} \\
&\; \; \; \text{and} \; \mc{F} = \{J \subset \{1,\ldots,p\}, \; |J|=n-1:\\
&\qquad \qquad  \quad \; \exists w \in
\R^n \; \text{s.t.} \; z^{\T} w = 0 \, \forall z \in \mc{C}_J \; \text{and} \;
z^{\T}w > 0 \;\forall z \in \mc{C} \setminus \mc{C}_J \}. 
\end{align*}
From $X \wh{\beta} \in \partial \mc{C}$, it hence follows that the active set $F =
\{j:\,\wh{\beta}_j > 0 \}$ has cardinality at most $n-1$. General linear
position implies that the linear system $X_F \beta = X \wh{\beta}$ has exactly
one solution $\beta = \wh{\beta}_F$.\\
Concerning the second part of the lemma, the assertion $|F| \leq \min\{n-1,p\}$ is
trivial for $p \leq (n-1)$. Conversely, if $p \geq n$, the fact that $\min_{\beta \gec 0} \frac{1}{n} \nnorm{y -
  X \beta}_2^2 > 0$ allows us to conclude that $X \wh{\beta} \in \partial \mc{C}$
so that the assertion follows from the reasoning above. For the third part, we
note that the fact that $y$ has a distribution which is absolutely continuous
with respect to the Lebesgue measure implies that $y$ is not contained in any
subspace of dimension smaller than $n$ with probability one, so that
$\min_{\beta \gec 0} \frac{1}{n} \nnorm{y - X \beta}_2^2 > 0$, and the claim
follows from part (ii).  
\end{bew}
Using Lemma \ref{lem:uniqueness} and the condition \eqref{eq:objectivepositivecond}
\begin{equation*}
\frac{32 (1 + M)^2 \sigma^2}{\E[\eps_1^2]}
\; \, \frac{\log p}{\tau^2(S) \, n} \leq \left(1 - \frac{s}{n} \right),
\end{equation*}
we will show that for $p \geq n$, condition (ii) of Lemma \ref{lem:uniqueness} holds with the
stated probability, from which we will conclude the proof of the first part
of the theorem. Note that for $p \leq n-1$, uniqueness follows from general
linear position while the claim $|F| \leq \min\{n-1, p\}$ is trivial.  
Let us recall the decomposition of Lemma \ref{lem:decoupling}. Note that 
\begin{equation*}
\min_{\beta \gec 0} \frac{1}{n} \nnorm{y - X \beta}_2^2 \geq \min_{\beta^{(P1)} \gec 0}
\frac{1}{n} \nnorm{\xi - Z \beta^{(P1)}}_2^2,
\end{equation*} 
hence it suffices to show that the right hand side is strictly positive. Suppose conversely that  
$\xi = Z \wh{\beta}^{(P1)}$, then $\frac{1}{n} \nnorm{\xi}_2^2 = \frac{1}{n}
\nnorm{Z \wh{\beta}^{(P1)}}_2^2$. Since $\wh{\beta}^{(P1)}$ is a minimizer of $(P1)$, 
$\frac{1}{n} \nnorm{Z \wh{\beta}^{(P1)}}_2^2 \leq \frac{2}{n} \xi^{\T} Z
\wh{\beta}^{(P1)}$, which, by the definition of $\tau^2(S)$, implies that 
\begin{equation*}  
\nnorm{\wh{\beta}^{(P1)}}_1 \leq \frac{1}{\tau^2(S)} \frac{2}{n}
\nnorm{Z^{\T} \xi}_{\infty}      
\end{equation*}
and in turn 
\begin{equation*}
  \frac{1}{n} \nnorm{\Pi_S^{\perp} \eps}_2^2 = \frac{1}{n} \nnorm{\xi}_2^2 =
  \frac{1}{n} \nnorm{Z \wh{\beta}^{(P1)}}_2^2 \leq \frac{1}{\tau^2(S)} \left( \frac{2}{n}
    \nnorm{Z^{\T} \xi}_{\infty} \right)^2
\end{equation*}      
Hence, conditional on the event
\begin{equation}\label{eq:condevent}
\left\{ \nnorm{\Pi_S^{\perp} \eps}_2^2 > \E[\eps_1^2] \, \frac{n-s}{4}
\right\} \cap \left \{ \left( \frac{2}{n}
    \nnorm{Z^{\T} \xi}_{\infty} \right)^2 \leq  8 (1 + M)^2 \sigma^2 \,
  \frac{\log p}{n}   \right \}
\end{equation}
it holds that
\begin{equation*}
\frac{\E[\eps_1^2]}{4} \left(1 - \frac{s}{n} \right) < \frac{1}{n}
\nnorm{\xi}_2^2 \leq 8 (1 + M)^2 \sigma^2 \frac{\log p}{\tau^2(S) \, n},
\end{equation*}
which contradicts \eqref{eq:objectivepositivecond}.  As a result, $\min_{\beta \gec 0} \frac{1}{n} \nnorm{y - X
  \beta}_2^2 > 0$ as was to be shown. Invoking Lemma \ref{lem:concprojection}
with $\Pi = \Pi_S^{\perp}$ so that $d = n-s$ by general linear position and
treating the second event in \eqref{eq:condevent} as in step 3 of Appendix \ref{app:robustsparserecovery}, the probability of
the event \eqref{eq:condevent} is no less than $1 - \exp(-c (n - s)/\sigma^4)
- 2 p^{-M^2}$.
\paragraph{Part 2: proof of the bound on $\nnorm{\wh{\beta} - \beta^*}_{\infty}$.} 
Given uniqueness of the NNLS solution and in turn of its active set
$F = \{j:\,\wh{\beta}_j > 0 \}$, the stated bound on $\nnorm{\wh{\beta} -
  \beta^*}_{\infty}$ follows readily once it holds that $S \subseteq F$. In
fact, the optimality conditions of the NNLS problem (cf.~Lemma
\ref{lem:NNLSopt}) then yield that $\wh{\beta}_F$ can be recovered from the linear
system
\begin{equation*}
\Sigma_{FF} \wh{\beta}_F = \frac{X_F^{\T} (X_S \beta_S^* + \eps)}{n} = \frac{X_F^{\T} (X_F \beta_F^* + \eps)}{n},
\end{equation*}
where the second equality results from $S \subseteq F$. As an immediate
consequence, we have that 
\begin{equation*}
\nnorm{\wh{\beta} - \beta^*}_{\infty} = \nnorm{\wh{\beta}_F -
  \beta_F^*}_{\infty} = \nnorm{\Sigma_{FF}^{-1} X_F^{\T} \eps / n}_{\infty}.   
\end{equation*}
In order to control the random term, we may follow the reasoning 
below \eqref{eq:AS} to conclude that for any $M \geq 0$, the event 
\begin{equation*}
\{ \nnorm{\Sigma_{FF}^{-1} X_F^{\T} \eps / n}_{\infty} \leq (1 + M)
\sigma \left\{ \phi_{\min}(F) \right \}^{-1/2} \sqrt{2 \log(p)/n}\}
\end{equation*}
has probability at least $1 - 2 p^{-M^2}$. It remains to show that under
the conditions of the theorem, we indeed have that $S \subseteq F$. This
is done by referring to the scheme in Appendix \ref{app:robustsparserecovery}
. Given the lower bound on $\beta_{\min}(S)$, it is established that the
event $\{ \wh{\beta}_S = \wh{\beta}^{(P2)} \succ 0 \}$ and in turn $\{S
\subseteq F \}$ occurs with probability at least $1 - 4 p^{-M^2}$.
This finishes the proof.
\section{Proof of Theorem \ref{theo:Sdatadriven}}

We first recall that the analysis is conditional on the event
\begin{equation}\label{eq:correctranking}
E = \{r_j \leq s \; \, \text{for all} \; j \in S \}, \quad \text{where} \;
\; r_j = k \, \Leftrightarrow \; \wh{\beta}_{j} = \wh{\beta}_{(k)}. 
\end{equation}
Our proof closely follows the corresponding proof in \cite{Genovese2012}. 
We show in two steps that both $S \setminus \wh{S} = \emptyset$ and $\wh{S}
\setminus  S = \emptyset$. For both steps, we shall need the following
observations. Let $V_{k}$ denote the linear space spanned by
the top $k$ variables according to the given ranking, $k=1,\ldots,p$, and 
let $V_0 = \{ 0 \}$. Let further $U_k = V_k^{\perp} \cap V_{k+1}$,
$k=0,\ldots,p-1$, which are subspaces of $\R^n$ of dimension at most $1$.
In case that the dimension of $U_k$ is one, let $u_k$ be the unit vector spanning  
$U_k$ and let $u_k = 0$ otherwise, $k=0,\ldots,p-1$. Note that $\Pi(k+1) -
\Pi(k)$ as appearing in the definition of the $\delta(k)$'s equals the projection   
on the $U_k$, $k=0,\ldots,p-1$. In particular, we have 
\begin{equation}\label{eq:projectednoise}
\nnorm{(\Pi(k + 1) - \Pi(k)) \eps}_2 = |\scp{u_k}{\eps}|, \; \; \,
k=0,\ldots,p-1. 
\end{equation}


\subsubsection*{Step 1: no false negatives} 
In the sequel, let $\Delta$ denote the threshold of the procedure so that
\begin{equation*}
\wh{s} = \max \left \{ 0 \leq k \leq (p-1): \delta(k) \geq \Delta \right\} + 1.
\end{equation*}
Later in the proof, it  will be verified that $\Delta$ can be chosen
as asserted in the theorem. We first note that conditional on $E$, by definition of $\wh{s}$,  
it holds that the event $\{S \setminus \wh{S} = \emptyset \}$ is contained in the
event $\{ \delta(s - 1) \geq \Delta \}$. Hence it suffices to upper bound the 
probability of the event $\{ \delta(s - 1) < \Delta \}$. We have
\begin{align}\label{eq:lowerdelta1}
\begin{split}
\p(\delta(s - 1) < \Delta) &= \p \left(\nnorm{(\Pi(s) - \Pi(s - 1)) y}_2 <
  \Delta \right) \\
                           &\leq \p \left(\nnorm{(\Pi(s) - \Pi(s - 1)) X_S
                             \beta_S^*}_2 < \Delta + \nnorm{(\Pi(s) - \Pi(s -
                             1)) \eps}_2 \right) \\
                           &\stackrel{\eqref{eq:projectednoise}}{=}
                             \p \left(\nnorm{(\Pi(s) - \Pi(s - 1)) X_S
                             \beta_S^*}_2 < \Delta + |\scp{u_{s-1}}{\eps}| \right) \\
                           &\leq \p \left(\min_{j \in S} \nnorm{(\Pi_{S} -
                               \Pi_{S \setminus j}) X_j \beta_j^{*}}_2 <
                             \Delta + |\scp{u_{s-1}}{\eps}| \right),                           
\end{split}
\end{align}
where $\Pi_S$ and $\Pi_{S \setminus j}$ denote the projection on the linear
spaces spanned by the columns of $X$ corresponding to $S$ respectively $S
\setminus j$, $j=1,\ldots,s$. In order to obtain the second inequality, we
have used again that we work conditional on the event $E$. As will be established
at the end of the proof, we further have  
\begin{align}\label{eq:lowerdelta2}
\min_{j \in S} \nnorm{(\Pi_S - \Pi_{S \setminus j}) X_j \beta_j^{\ast}}_2 \geq
\sqrt{n} \left\{ \phi_{\min}(S) \right \}^{1/2}  \beta_{\min}(S). 
\end{align}
%
Combining \eqref{eq:lowerdelta1} and \eqref{eq:lowerdelta2}, it suffices to upper bound
\begin{equation}\label{eq:upperbound_fn}
\p \left(|\scp{u_{s-1}}{\eps}| > \sqrt{n} \left\{ \phi_{\min}(S)
  \right\}^{1/2} \beta_{\min}(S) - \Delta \right)
\end{equation}
as will be done below after fixing $\Delta$. 
\subsubsection*{Step 2: no false positives}
Conditional on $E$, the probability of having a false positive selection is upper bounded as
\begin{align}\label{eq:upperboundfalsepositive}
  \p(\cup_{k=s}^{p-1} \{ \delta(k) \geq \Delta \}) &= \p \left(\max_{s \leq
      k \leq p-1}
\nnorm{(\Pi(k+1) - \Pi(k)) y}_2 \geq \Delta \right) \notag\\
& = \p \left(\max_{s \leq
      k \leq p-1} \nnorm{(\Pi(k+1) - \Pi(k)) \eps}_2 \geq
  \Delta \right) \notag\\
&= \p \left(\max_{s \leq
      k \leq p-1} |\scp{u_k}{\eps}| \geq
  \Delta \right).
\end{align}
Choosing $\Delta = (1 + M) \sigma \sqrt{2 \log(p)}$ for an arbitrary $M
\geq 0$, using the assumption on $\beta_{\min}(S)$, and controlling 
$\max_{0 \leq  k \leq p-1} |\scp{u_k}{\eps}|$ according to
\eqref{eq:subgaussianmaximal} in the usual way, the probabilities 
\eqref{eq:upperbound_fn} and \eqref{eq:upperboundfalsepositive}  
do not exceed $2 p^{-M^2}$. The assertion of the theorem then follows. To conclude the proof, it remains to establish
\eqref{eq:lowerdelta2}. Let us fix an arbitrary $j \in S$. We have 
\begin{equation*}
\nnorm{(\Pi_S - \Pi_{S \setminus j}) X_j}_2 = \nnorm{X_j - \Pi_{S \setminus
    j} X_j}_2 = \sqrt{\nnorm{X_j}_2^2 -  \nnorm{\Pi_{S \setminus j} X_j}_2^2}
\end{equation*}
Write $\theta$ for the vector of regression coefficients for
the linear regression of $X_j$ on $\{X_k\}_{k \in S \setminus \{ j \}}$, i.e.~
\begin{equation*}
\theta = (X_{S \setminus j}^{\T} X_{S \setminus j})^{-1} X_{S \setminus j}^{\T} X_j, 
\end{equation*}
and note that, according to a block decomposition of the matrix $X_S^{\T} X_S$ 
\begin{align*}
\begin{pmatrix}
                 -\theta \\
                 1  
                 \end{pmatrix}^{\T}  (X_S^{\T} X_S) \begin{pmatrix}
                 -\theta \\
                 1  
                 \end{pmatrix} &=  \begin{pmatrix}
                 -\theta \\
                 1  
                 \end{pmatrix}^{\T} \begin{pmatrix}
                                   X_{S \setminus j}^{\T} X_{S \setminus j} &
                                   X_{S \setminus j}^{\T} X_j  \\
                                  X_j^{\T} X_{S \setminus j} & \norm{X_j}_2^2 
                                  \end{pmatrix} \begin{pmatrix}
                 -\theta \\
                 1  
                 \end{pmatrix} \\
&= \norm{X_j}_2^2 - X_j^{\T} X_{S \setminus j} (X_{S \setminus j}^{\T} X_{S
  \setminus j})^{-1} X_{S \setminus j}^{\T} X_j \\
&= \norm{X_j}_2^2 - \nnorm{\Pi_{S \setminus j} X_j}_2^2.
\end{align*}
We conclude the proof from $X_S^{\T} X_S = n \Sigma_{SS}$ and 
\begin{equation*}
\begin{pmatrix}
-\theta \\
                 1  
                 \end{pmatrix}^{\T}  (X_S^{\T} X_S) \begin{pmatrix}
                 -\theta \\
                 1  
                 \end{pmatrix} \geq n \phi_{\min}(S) \norm{\begin{pmatrix}
                 -\theta \\
                 1
               \end{pmatrix}
             }_2^2 \geq n \phi_{\min}(S). 
\end{equation*}

\section{Proof of Theorem \ref{theo:plasso}}\label{app:positivelasso}

Consider the non-negative lasso problem \eqref{eq:nnlasso}. It follows from the
KKT optimality conditions that any minimizer $\plasso$ of \eqref{eq:nnlasso} satisfies
\begin{align}\label{eq:nnlassopt}
\begin{split}
&\; \; \; \; \frac{2}{n} X_j^{\T} (y - X \plasso) = \lambda \; \; \text{and} \; \; \plasso_j > 0, \\ 
\text{or}& \quad \frac{2}{n} X_j^{\T} (y - X \plasso) \leq \lambda \; \; \text{and} \; \;
\plasso_j = 0, \quad j=1,\ldots,p.
\end{split}
\end{align}
Following the technique employed in \cite{Wain2009}, we establish that under
the conditions of the theorem, the unique minimizer of the non-negative lasso
problem is given by $\plasso_S = \wh{\alpha}_S \succ 0$ and
$\plasso_{S^c} = 0$ with the specified probability, where $\wh{\alpha}$ denotes the minimizer of
the following constrained non-negative lasso problem
\begin{equation}\label{eq:constrainednnlasso}
\min_{\beta_S \gec 0, \; \beta_{S^c} = 0} \frac{1}{n} \norm{y - X \beta}_2^2 +
\lambda \bm{1}^{\T} \beta.  
\end{equation}
To this end, in view of \eqref{eq:nnlassopt}, it suffices to show that the following system of inequalities is satisfied
\begin{equation}\label{eq:systemnnlasso}
\frac{2}{n} \left[ \begin{array}{c} X_S^{\T} \eps \\ X_{S^c}^{\T} \eps \end{array} \right]+ 2 \left[ \begin{array}{cc}
\Sigma_{SS}  & \Sigma_{SS^c}  \\
\Sigma_{S^cS}  & \Sigma_{S^cS^c}    
\end{array} \right]     \left[\begin{array}{c}
                               \beta_S^{\ast} - \wh{\alpha}_S\\
                               0 
                              \end{array} \right]
                                                     \left[ \begin{array}{c} 
                                                       =   \\
                                                       \prec 
                                                      \end{array} \right]
\left[
                                                       \begin{array}{c}
                                                      \lambda \bm{1} \\
                                                      \lambda \bm{1}
                                                     \end{array}
                                                    \right].
\end{equation}
In view of the required lower bound on $\beta_{\min}(S)$ in
\eqref{eq:nnlassocondition}, we have 
\begin{equation}\label{eq:alphahatS}
0 \prec \wh{\alpha}_S = \beta_S^{\ast} - \frac{\lambda}{2} \Sigma_{SS}^{-1}
\bm{1} + \Sigma_{SS}^{-1} \frac{1}{n} X_S^{\T} \eps
\end{equation}
with probability at least $1 - 2 p^{-M^2}$, handling the random term as
\eqref{eq:AS} in the proof of Theorem \ref{theo:robustsparserecovery}.  
Substituting \eqref{eq:alphahatS} back into \eqref{eq:systemnnlasso}, we find that the following
system of inequalities must hold true:
\begin{equation}\label{eq:offsupportzerocond}
\frac{\lambda}{2} \left \{ \Sigma_{S^c S} \Sigma_{SS}^{-1}
\bm{1} \right \}+ \frac{1}{n} X_{S^c}^{\T} (I - \Pi_S) \eps \prec \frac{\lambda}{2} \bm{1}. 
\end{equation}
In light of the maximal inequality \eqref{eq:subgaussianmaximal}, the event
\begin{equation}\label{eq:eventplasso1}
\left \{ \max_{j \in S^c}  \frac{1}{n} X_j^{\T} (I - \Pi_S) \eps \leq \lambda_M \right\},
\quad \lambda_M = \sigma (1 + M) \sqrt{2 \log(p)/n},
\end{equation}
occurs with probability at least $1 - 2 p^{-M^2}$, noting that $\nnorm{(I - \Pi_S) X_j}_2 \leq \nnorm{X_j}_2$ for all $j=1,\ldots,p$. Hence, conditional on
the events $\{ \wh{\alpha}_S \succ 0 \}$ and \eqref{eq:eventplasso1}, each
component of the left hand side of \eqref{eq:offsupportzerocond} is no larger than $\frac{\lambda}{2} \iota(S) +
\lambda_M$ (cf.~definition \eqref{eq:iotaS}), so that for $\lambda > 2 \lambda_M / (1 - \iota(S))$, the system of
inequalities \eqref{eq:offsupportzerocond} and hence also
\eqref{eq:systemnnlasso} are indeed fulfilled.

\section{Example of a design for which the non-negative lasso performs sub-optimally}\label{app:examplesup}

As indicated in the discussion following Theorem \ref{theo:plasso}, the
non-negative lasso does not always attain the optimal rate for estimating
$\beta^*$ with respect to the $\ell_{\infty}$-norm. For the data-generating
model of Theorem \ref{theo:plasso} we give an example of a
design so that the non-negative lasso estimator has the property 
\begin{equation}\label{eq:supnormlowerbound}
\nnorm{\plasso - \beta^*}_{\infty} = \Omega\left(\sqrt{s \log(p)/n} \right),
\end{equation}
with high probability, provided the regularization parameter $\lambda = \Omega(\sqrt{\log(p)/n})$ (as
conventionally suggested in the literature) and the minimum support
coefficient $\beta_{\min}(S) = \Omega(\sqrt{s \log(p)/n})$.
For the same design, if $\beta_{\min}(S) = \Omega(\sqrt{\log(p)/n})$, the NNLS estimator obeys the bound 
\begin{equation}\label{eq:supnormupperboundnnls}
\nnorm{\wh{\beta} - \beta^*}_{\infty} = O \left(\sqrt{\log(p)/n} \right),
\end{equation}
with high probability. In order to establish \eqref{eq:supnormlowerbound}, we
shall build on the scheme used for the proof of Theorem \ref{theo:plasso} in
the previous paragraph. Consider a design whose Gram matrix is of the form
\begin{equation*}
\Sigma =    \left[ \begin{array}{cc}
             \Sigma_{SS}   &  0            \\
                          &                \\
              0           & I_{p-s}
             \end{array}  \right],
\end{equation*}
where 
\begin{equation*}
\Sigma_{SS} =  \left[ \begin{array}{cc}
              1         &   -1/\sqrt{2(s-1)}\bm{1}_{s-1}^{\T}             \\
                        &                     \\
              -1/\sqrt{2(s-1)} \bm{1}_{s-1}            &  I_{s-1}
             \end{array}  \right].
\end{equation*}
The constant $\sqrt{2}$ in the denominator is chosen for convenience; any other constant
larger than $1$ would do as well. Using Schur complements, one computes that  
\begin{equation*}
\Sigma_{SS}^{-1} = \left[ \begin{array}{cc}
 2  &  \sqrt{2/(s-1)} \bm{1}_{s-1}^{\T}   \\
    &                                     \\   
 \sqrt{2/(s-1)} \bm{1}_{s-1}   &  I_{s-1} + \frac{1}{s-1} \bm{1}_{s-1}
 \bm{1}_{s-1}^{\T}
\end{array} \right].
\end{equation*}
As a result, we have that
\begin{equation}\label{eq:lowerboundSigmaSSone}
e_1^{\T} \Sigma_{SS}^{-1} \bm{1} = 2 + \sqrt{2 (s-1)} = \Omega(\sqrt{s}),
\end{equation}
where $e_1$ is the first canonical basis vector. Furthermore, the 
sequence of eigenvalues of $\Sigma_{SS}$ (ordered decreasingly) is given by
\begin{equation*}
\phi_1 = 1 + \frac{1}{\sqrt{2}}, \quad \phi_2 = \ldots = \phi_{s-1} =
1, \quad \phi_s = 1 - \frac{1}{\sqrt{2}}. 
\end{equation*}
From the proof of Theorem \ref{theo:plasso}, we know that given the active set $Q = \{j: \plasso_j
> 0\}$, the non-negative lasso estimator has the following closed form
expression. 
\begin{equation}\label{eq:nnlassotoactiveset}
\plasso_Q = \beta_Q^{\ast} - \frac{\lambda}{2} \Sigma_{QQ}^{-1}
\bm{1} + \Sigma_{QQ}^{-1} \frac{1}{n} X_Q^{\T} \eps, \quad \text{and} \; \; \, \plasso_{Q^c} = 0.
\end{equation}
If $S \nsubseteq Q$, the claim \eqref{eq:supnormlowerbound} follows trivially
from the assumption $\beta_{\min}(S) = \Omega(\sqrt{s
  \log(p)/n})$. Conversely, if $S \subseteq Q$, \eqref{eq:nnlassotoactiveset}
and the block structure of $\Sigma$ imply that 
\begin{align*}   
\nnorm{\plasso - \beta^*}_{\infty} &\geq \nnorm{\plasso_Q - \beta_Q^*}_{\infty} \\
                                   &\geq \nnorm{\plasso_S -
                                     \beta_S^*}_{\infty} \\
&\geq \frac{\lambda}{2} e_1^{\T} \Sigma_{SS}^{-1} \bm{1} - \left| e_1^{\T} \Sigma_{SS}^{-1}
\frac{1}{n} X_S^{\T} \eps \right| \\
&\geq \frac{\lambda}{2} e_1^{\T} \Sigma_{SS}^{-1} \bm{1} - \frac{1}{\sqrt{\phi_s}}  \max_{1 \leq j
  \leq p} \left| 
 \frac{X_j^{\T} \eps}{n} \right| = \Omega \left(\sqrt{s \log(p)/n} \right),
\end{align*}
using \eqref{eq:lowerboundSigmaSSone} and the maximal inequality \eqref{eq:subgaussianmaximal} to upper bound the
second term after the third inequality. The latter scales as $O(\sqrt{\log(p)/n})$ with
probability at least $1 - O((p \vee n)^{-1})$. Furthermore, we have used that
$\lambda = \Omega(\sqrt{\log(p) / n})$ and $\phi_s = \Omega(1)$.\\ 
Regarding the upper bound \eqref{eq:supnormupperboundnnls} for NNLS, we note
that the optimality conditions of the NNLS problem (Lemma \ref{lem:NNLSopt})
in conjunction with the block structure of $\Sigma$ imply that 
\begin{equation*}
\wh{\beta}_S \succ 0 \; \; \, \text{and} \; \; \nnorm{\wh{\beta}_S -
  \beta_S^*}_{\infty} \leq \nnorm{\Sigma_{SS}^{-1} X_S^{\T} \eps/n}_{\infty} \leq
\frac{1}{\sqrt{\phi_s}} \nnorm{X_S^{\T} \eps / n}_{\infty} = O(\sqrt{\log(p)/n}),
\end{equation*}
with probability at least $1 - O((p \vee n)^{-1})$, provided $\beta_{\min}(S)$
exceeds the bound on the right hand side. Similarly, 
\begin{equation*}
\nnorm{\wh{\beta}_{S^c}}_{\infty} \leq \nnorm{X_{S^c}^{\T} \eps / n}_{\infty} = O(\sqrt{\log(p)/n})
\end{equation*}
with probability at least $1 - O((p \vee n)^{-1})$, so that
\eqref{eq:supnormupperboundnnls} follows.  

\section{Proof of Proposition \ref{prop:sparsity_equicor}}

We start by noting that $\Sigma$ is strictly positive definite so that the NNLS
problem is strictly convex. Thus, the NNLS solution and its active set $F = \{j:\,\wh{\beta}_j > 0 \}$ are unique.
Let us first consider the case $s > 0$. Using a slight modification of the scheme used in the proofs of Theorem
\ref{theo:robustsparserecovery} and \ref{theo:ellinfsmin}, we will show that 
under the required condition on $\beta_{\min}(S)$, the event $\{\wh{\beta}_S
= \wh{\beta}^{(P2)} \succ 0 \}$ holds with the stated probability, which
proves the first statement of the proposition. Following the proof of Theorem
\ref{theo:robustsparserecovery}, we have that 
\begin{equation*}
\nnorm{\wh{\beta}^{(P1)}}_1 \leq \frac{2 (1 + M) \sigma \sqrt{2 \log(p)/n}}{\tau^2(S)} 
  \stackrel{\eqref{eq:tauS_equicor}}{\leq} \frac{2 (1 + (s - 1) \rho) (1 + M) \sigma \sqrt{2
      \log(p)/n}}{\rho (1 - \rho)},
\end{equation*}
with probability at least  $1 - 2 p^{-M^2}$, where we have used the closed form expression for $\tau^2(S)$ in
\eqref{eq:tauS_equicor}. In order to verify that $\wh{\beta}_S =
\wh{\beta}^{(P2)} \succ 0$, we follow the back-substitution step (step 2 in
Appendix \ref{app:robustsparserecovery}) apart from the following
modification. In place of \eqref{eq:backsubsbound}, we bound
\begin{align*}
\nnorm{\Sigma_{SS}^{-1} \Sigma_{SS^c} (\wh{\beta}^{(P1)} -
  \beta_{S^c}^*)}_{\infty} &= \rho \nnorm{\Sigma_{SS}^{-1} \bm{1}}_{\infty}
\nnorm{\wh{\beta}^{(P1)}}_1 \\
&\leq \frac{\rho}{1 + (s - 1) \rho} \nnorm{\wh{\beta}^{(P1)}}_1
                           \leq \frac{2 (1 + M) \sigma \sqrt{2 \log(p)/n}}{1
                             - \rho}
\end{align*}
For the first equality, we have used that $\beta_{S^c}^* = 0$ and the
fact that the matrix $\Sigma_{SS^c}$ has constant entries equal to $\rho$. For
the second inequality, we have used that $\bm{1}$ is an eigenvector of
$\Sigma_{SS}$ corresponding to its largest eigenvalue $1 + (s - 1)
\rho$. Turning to step 3 in Appendix \ref{app:robustsparserecovery}, we note
that with $\phi_{\min}(S) = (1 - \rho)$, 
\begin{align*}
\nnorm{\beta_S^* - \bar{\beta}^{(P2)}}_{\infty} &\leq \frac{2 (1 + M) \sigma \sqrt{2 \log(p)/n}}{1
                             - \rho} +  \frac{(1 + M) \sigma \sqrt{2
                               \log(p)/n}}{\sqrt{1 - \rho}} \\
&\leq \frac{3 (1 + M) \sigma \sqrt{2 \log(p)/n}}{1
                             - \rho}
\end{align*}
so that $\bar{\beta}^{(P2)} = \wh{\beta}^{(P2)} = \wh{\beta}_S \succ 0$ with probability at least $1 - 4 p^{-M^2}$ as claimed.\\       
We now turn to the second statement of the proposition concerning the
(conditional) distribution of the cardinality of the active set.  
Conditional on the event $\{\wh{\beta}_S  \succ 0 \}$, the KKT optimality conditions of the NNLS
problem as stated in Lemma \ref{lem:NNLSopt} imply that the following block system of
inequalities holds.
\begin{equation}\label{eq:KKTblock}
\left[ \begin{array}{cc} \Sigma_{SS} & \Sigma_{SS^c} \\ 
                                     & \\                         \Sigma_{S^cS} & \Sigma_{S^cS^c}
                      \end{array}
                     \right] \left[ \begin{array}{c}
                                               \wh{\beta}_S  \\
                                                       \\
                                               \wh{\beta}_{S^c} 
                                               \end{array}
                                               \right]    \left[\begin{array}{c}
                                                                = \\
                                                                 \\
                                                                \lec
                                                              \end{array}
                                                            \right]  
                                                          \left[ \begin{array}{c}
                                               \Sigma_{SS} \beta_S^*  +
                                               \frac{X_S^{\T} \eps}{n}  \\
\\
                                               \Sigma_{S^c S} \beta_S^* +
                                               \frac{X_{S^c}^{\T} \eps}{n} 
                                               \end{array}
                                               \right].
\end{equation}
Resolving the top block for $\wh{\beta}_S$, we obtain 
\begin{equation*}
\wh{\beta}_S = \beta_S^* + \Sigma_{SS}^{-1} \left(\frac{X_S^{\T} \eps}{n}  -
  \Sigma_{SS^c} \wh{\beta}_{S^c} \right).
\end{equation*}
Back-substituting that expression into the bottom block of inequalities
yields the following system of inequalities.
\begin{equation}\label{eq:KKTblockSc}
\left( \Sigma_{S^c S^c} - \Sigma_{S^c S} \Sigma_{SS}^{-1} \Sigma_{S S^c}
\right) \wh{\beta}_{S^c} \lec \frac{X_{S^c}^{\T} (I - X_S
  (X_S^{\T} X_S)^{-1} X_S^{\T}) \eps}{n} = \frac{Z^{\T} \eps}{n},  
\end{equation} 
where $Z = \Pi_S^{\perp} X_{S^c}$.
For equi-correlated design with $\Sigma = (1 - \rho) I + \rho \bm{1}
\bm{1}^{\T}$, we have that 
\begin{equation}\label{eq:ZtZequicor}
\Sigma_{S^c S^c} - \Sigma_{S^c S} \Sigma_{SS}^{-1} \Sigma_{S S^c} = (1 - \rho) I + \underbrace{\frac{\rho (1 - \rho)}{1 + (s - 1)
    \rho}}_{\gamma(s, \rho)}
\bm{1} \bm{1}^{\T} = (1 - \rho) I + \gamma(s, \rho) \bm{1} \bm{1}^{\T},
\end{equation}
cf.~the derivation in \eqref{eq:tauS_equicor}. Denote $\wh{\alpha} = \wh{\beta}_{S^c}$,
and $G = \{k:\; \wh{\alpha}_k > 0 \}$. Using Lemma \ref{lem:NNLSopt} and \eqref{eq:ZtZequicor}, \eqref{eq:KKTblockSc} can be
written as
 \begin{align}\label{eq:kkt_P1equicor}
 \begin{split}
 &\frac{Z_k^{\T} \eps}{n} - (1 - \rho) \wh{\alpha}_k = \gamma(s, \rho)\,
 \bm{1}^{\T} \wh{\alpha}, \; \; k \in G, \\
 &\frac{Z_k^{\T} \eps}{n} \leq \gamma(s, \rho) \bm{1}^{\T} \wh{\alpha}, \;
 \; k \notin G. 
\end{split}
\end{align}
Set $z = Z^{\T} \eps / (\sigma \sqrt{n})$ so that $z$ is a zero-mean Gaussian 
random vector with covariance 
\begin{equation*}
\frac{1}{n} Z^{\T} Z = \frac{1}{n} X_{S^c}^{\T} \Pi_S^{\perp} X_{S^c} = \Sigma_{S^c S^c} - \Sigma_{S^c S} \Sigma_{SS}^{-1} \Sigma_{S S^c}.
\end{equation*}
In view of \eqref{eq:ZtZequicor}, $z$ has the distribution as claimed in 
Proposition \ref{prop:sparsity_equicor}. From \eqref{eq:kkt_P1equicor}, we
conclude that 
\begin{equation}\label{eq:ZtZequicor_conclusion}
k \in G \Rightarrow z_k > 0 \; \; \text{and} \; \; z_k \leq z_l \; \text{for}
\; l \notin G \Rightarrow k \notin G.
\end{equation}
In particular, recalling that $z_{(1)} \geq z_{(2)} \geq \ldots \geq z_{(p-s)}$ denotes the arrangement of the components of $z$ in
decreasing order, if $z_{(1)} \leq 0$, then \eqref{eq:kkt_P1equicor} implies
that $\wh{\alpha} = 0$, $G = \emptyset$ and $|F| = s$ as stated in the
proposition. Let us henceforth assume that $z_{(1)} > 0$, in which case
\eqref{eq:kkt_P1equicor} implies that $G$ is non-empty. We may then resolve the first set of
equations in  \eqref{eq:kkt_P1equicor} with respect to $\wh{\alpha}_G$, which yields
\begin{equation*}
 \wh{\alpha}_{G} = \left( (1 - \rho) I + \gamma(s, \rho) \bm{1}
   \bm{1}^{\T} \right)^{-1} \frac{Z_{G}^{\T} \eps}{n} = \frac{1}{1 - \rho}
 \left(\frac{Z_{G}^{\T} \eps}{n} -  \frac{\gamma(s, \rho) \bm{1}
     \bm{1}^{\T} (Z_{G}^{\T} \eps/n)}{(1 - \rho) +
     \gamma(s, \rho) |G|} \right),
\end{equation*}
where the second equality is an application of the Sherman-Woodbury-Morrison
formula. This implies in turn that 
\begin{equation*}
\bm{1}^{\T} \wh{\alpha} = \bm{1}^{\T} \wh{\alpha}_{G} =
 \frac{\bm{1}^{\T} Z_{G}^{\T} \eps}{n} \frac{1}{(1 - \rho) + |G|
     \gamma(s, \rho)}.  
 \end{equation*}
Substituting this expression back into \eqref{eq:kkt_P1equicor}, we obtain
\begin{align}\label{eq:kkt_P1equicor_2}
 \begin{split}
 & \frac{\frac{Z_k^{\T} \eps}{n}}{\sum_{\ell \in G} \left(
     \frac{Z_\ell^{\T} \eps}{n}- \frac{Z_k^{\T} \eps}{n} \right)} -
 \wh{\alpha}_k \frac{(1 - \rho) + |G|
     \gamma(s, \rho)}{\sum_{\ell \in G} \left(
     \frac{Z_\ell^{\T} \eps}{n}- \frac{Z_k^{\T} \eps}{n} \right)} = \frac{\gamma(s, \rho)}{1 -
   \rho}, \; \; k \in G, \\
 & \frac{\frac{Z_k^{\T} \eps}{n}}{\sum_{\ell \in G} \left(
     \frac{Z_\ell^{\T} \eps}{n}- \frac{Z_k^{\T} \eps}{n} \right)} \leq
 \frac{\gamma(s, \rho)}{1 - \rho}, \; \; k \notin G.
\end{split}
\end{align}
Now note that for $k=1,\ldots,p-s$,
\begin{equation*}
\frac{\frac{Z_k^{\T} \eps}{n}}{\sum_{\ell \in G} \left(
     \frac{Z_\ell^{\T} \eps}{n}- \frac{Z_k^{\T} \eps}{n} \right)} = \frac{\frac{Z_k^{\T} \eps}{\sqrt{n}}}{\sum_{\ell \in G} \left(
     \frac{Z_\ell^{\T} \eps}{\sqrt{n}}- \frac{Z_k^{\T} \eps}{\sqrt{n}}
   \right)} = \frac{z_k}{\sum_{\ell \in G} \left(z_{\ell} - z_k \right)}
\end{equation*}
From
\begin{equation*}
\underbrace{\frac{z_{(2)}}{(z_{(1)} - z_{(2)})}}_{\zeta_1} \geq \underbrace{\frac{z_{(3)}}{(z_{(1)} -
    z_{(3)}) + (z_{(2)} -
    z_{(3)})}}_{\zeta_2}  \geq \ldots \geq
\underbrace{\frac{z_{(p-s)}}{\sum_{k=1}^{p-s-1} (z_{(k)} - z_{(p-s)})}}_{\zeta_{p-s-1}},  
\end{equation*}
\eqref{eq:ZtZequicor_conclusion} and the inequalities in
\eqref{eq:kkt_P1equicor_2}, it then follows that 
\begin{align*}
G &=  \{j:\,z_{j} = z_{(1)} \} \cup \left\{k \neq j:\, \frac{z_k}{\sum_{\ell: z_{\ell} \geq z_k} (z_{\ell} - z_k)}  >
  \frac{\gamma(s, \rho)}{1 - \rho} \right\}  \\
  &= \{j:\,z_{j} = z_{(1)} \} \cup \left\{k \neq j:\, \frac{z_k}{\sum_{{\ell}:
      z_{\ell} \geq z_k} (z_{\ell} - z_k)}  \geq \zeta_m \right \},
\end{align*}
where $m$ is the largest integer so that $\zeta_{m} > \gamma(s, \rho)/(1 -
\rho) = \theta(s, \rho)$ with $\theta(s, \rho)$ as defined in
\eqref{eq:distn_activeset}, which finishes the proof for $s > 0$. Turning to
the case $s = 0$, a similar scheme can be used, starting from the system 
of inequalities $\Sigma \wh{\beta} \lec \frac{X^{\T} \eps}{n} = \frac{Z^{\T}
  \eps}{n} = 
\sigma z / \sqrt{n}$. The expressions used above remain valid with $\gamma(0, \rho) = \rho$.

\section{Proof of Propositions \ref{prop:re_and_selfreg} and \ref{prop:scalingtauS}}

Proofs of Proposition \ref{prop:re_and_selfreg} and Proposition
\ref{prop:scalingtauS} are contained in the supplement.

\bibliographystyle{plainnat}
\bibliography{/warehouse/Proteomics/nips2010/references}

\clearpage

\begin{center}
{\LARGE{\textbf{Supplementary material to\\ 'Non-negative least squares for high-dimensional linear
  models: consistency and sparse recovery\\ without regularization'}}}  
\end{center} 

\setcounter{lemmachen}{0}

\section*{Proofs of Propositions 3 and 4}

We here provide proofs of Propositions 3 and 4 concerning
random equi-correlation-like matrices. These proofs rely
on a series of lemmas that are stated first. 

\subsection*{Additional lemmas}
We recall from Appendix A of the paper that a zero-mean random variable is called sub-Gaussian
if there exists $\sigma > 0$ (referred to as sub-Gaussian parameter) so that the
moment-generating function obeys the bound $\E[\exp(t Z)] \leq \exp(\sigma^2
t^2/2) \; \forall t \in \R$. If $Z_1,\ldots,Z_n$ are i.i.d.~copies of $Z$ and $v_j \in \R^n$,
$j=1,\ldots,p$, are fixed vectors, then 
\begin{equation}\label{eq:subgaussianmaximal}
\p \left(\max_{1 \leq j \leq p} |v_j^{\T} \mathbf{Z}| > \sigma \max_{1 \leq j
    \leq p} \norm{v_j}_2 \left( \sqrt{2 \log
    p} + z \right)\right) \leq 2 \exp \left(-\frac{1}{2} z^2 \right), \; \; z
\geq 0.
\end{equation}

\subsection*{Bernstein-type inequality for squared sub-Gaussian random
  variables}

The following exponential inequality combines Lemma 14, Proposition 16 and
Remark 18 in [4].
\begin{lemma}\label{lemma:subexponential} Let $Z_1,\ldots,Z_m$ be i.i.d.~zero-mean sub-Gaussian
random variables with parameter $\sigma$ and the property that $\E[Z_1^2] \leq
1$. Then for any $z \geq 0$, one has  
\begin{equation}\label{eq:concentration_normsubgaussian}
\p \left(\sum_{i = 1}^m Z_i^2 > m + z m \right) \leq \exp \left(-c \min \left\{
\frac{z^2}{\sigma^4}, \frac{z}{\sigma^2} \right \} m \right),
\end{equation}
where $c > 0$ is an absolute constant.
\end{lemma}

\subsection*{Concentration of extreme singular values of sub-Gaussian random
  matrices}

Let $s_{\min}(A)$ and $s_{\max}(A)$ denote the minimum and maximum singular
value of a matrix $A$. The following lemma is a special case of Theorem 39 in [4]. 
\begin{lemma}\label{lemma:singularvalues_general} Let $A$ be an $n \times s$ matrix with
  i.i.d.~zero-mean sub-Gaussian entries with sub-Gaussian
  parameter $\sigma$ and unit variance. Then for every $z
  \geq 0$, with probability at least $1 - 2 \exp(-c z^2)$, one has 
\begin{equation}\label{eq:deviationcovariance}
s_{\max} \left(\frac{1}{n} A^{\T} A - I \right) \leq \max(\delta, \delta^2), \;
\; \; \text{where} \; \; \delta = C \sqrt{\frac{s}{n}} + \frac{z}{\sqrt{n}}, 
\end{equation}
with $C$, $c$ depending only on $\sigma$.  
\end{lemma}

 

\subsection*{Entry-wise concentration of the Gram matrix associated with a
  sub-Gaussian random matrix}

The next lemma results from Lemma 1 in [2] and the union bound. 
\begin{lemma}\label{lemma:entrywise_general} Let $X$ be an $n \times p$ random matrix of i.i.d.~zero-mean, unit variance
sub-Gaussian entries with parameter $\sigma$. Then 
\begin{equation}\label{eq:ellinf_gram}
\p \left( \max_{1 \leq j,k \leq p} \left| \left(\frac{1}{n} X^{\T} X
      - I \right)_{jk} \right| > z
\right) \leq 4 p^2 \exp \left( -\frac{n z^2}{128(1 + 4 \sigma^2)^2}\right)
\end{equation}
for all $z \in \left(0, 8(1 + 4 \sigma^2) \right)$. 
\end{lemma}

\subsection*{Application to $\text{Ens}_+$}
Recall that the class $\text{Ens}_+$ is given by 
\begin{equation}\label{eq:Ensplus}
\text{Ens}_+: \;  X = (x_{ij})_{\substack{1 \leq i \leq n \\ 1 \leq j \leq
p}}, \;  \{ x_{ij} \} \; \text{i.i.d. from a sub-Gaussian distribution on} \; \R_+.
\end{equation}
We shall make use of the following decomposition valid for any $X$ from
\eqref{eq:Ensplus}.
\begin{equation}\label{eq:decomposition_Ensplus_general}
X = \wt{X} + \mu \one,
\end{equation}
where the entries $\{\wt{x}_{ij} \}$ of $\wt{X}$ are zero mean sub-Gaussian
random variables with parameter $\sigma$, say, $\mu =
\E[x_{11}]$ and $\one$ is an $n \times p$-matrix of ones. In the sequel,
we specialize to the case where the entries of $X$ are scaled such that 
\begin{equation}\label{eq:Sigmastar}
\Sigma^* = \E \left[\frac{1}{n}  X^{\T} X \right] = (1 - \rho) I + \rho \bm{1}
\bm{1}^{\T} 
\end{equation}
for $\rho \in (0,1)$, i.e.~the population Gram matrix has equi-correlation
structure. Then, decomposition \eqref{eq:decomposition_Ensplus_general} becomes 
\begin{equation}\label{eq:decomposition_Ensplus_specific}
X = \wt{X} + \sqrt{\rho} \one,  \quad \text{and} \; \E[\wt{x}_{11}^2] = (1-\rho). 
\end{equation}
Accordingly, we have the following expansion of $\Sigma = \frac{1}{n} X^{\T} X$. 
\begin{equation}\label{eq:decomposition_gram}
\Sigma  = \frac{1}{n} \wt{X}^{\T} \wt{X} + \sqrt{\rho} \left( \frac{1}{n} \wt{X}^{\T}
\one  + \frac{1}{n} \one^{\T} \wt{X}  \right) + \rho
\bm{1} \bm{1}^{\T}, \quad \text{where} \; \; \E \left[ \frac{1}{n} \wt{X}^{\T}
  \wt{X}  \right] = (1 - \rho) I. 
\end{equation}
Observe that 
\begin{equation}\label{eq:crossterm}
n^{-1} \wt{X}^{\T} \one = D \bm{1} \bm{1}^{\T}, \quad \text{and}
\; \; n^{-1} \one^{\T} \wt{X} = \bm{1} \bm{1}^{\T} D,  
\end{equation}
where $D \in \R^{p \times p}$ is a diagonal matrix with diagonal entries 
$d_{jj} = n^{-1} \sum_{i = 1}^n \wt{x}_{ij}$, $j=1,\ldots,p$. It hence follows
from \eqref{eq:subgaussianmaximal} that
\begin{equation}\label{eq:control_crossterm}
\p \left( \max_{j,k} \left| n^{-1} \wt{X}^{\T} \one \right|_{jk} > 2
  \sigma \sqrt{\frac{2 \log(p \vee n)}{n}} \right) \leq \frac{2}{p \vee n},
\end{equation}  
Combining \eqref{eq:Sigmastar}, \eqref{eq:decomposition_gram}, \eqref{eq:control_crossterm} and
Lemma \ref{lemma:entrywise_general}, it follows that there exists a constant
$C > 0$ depending only on $\sigma$ such that
\begin{equation}\label{deviation:gram}
\p \left( \max_{j,k} \left|\left(\frac{X^{\T} X}{n} - \Sigma^*
    \right)_{jk}\right| > C \sqrt{\frac{\log(p \vee n)}{n}} \right) \leq \frac{6}{p \vee n}. 
\end{equation}
Let now $S \subset \{1,\ldots,p\}$, $|S|=s < n$ be given. Without loss
of generality, let us assume that $S = \{1,\ldots,s \}$. In the sequel, we
control $s_{\max}(\Sigma_{SS}^* - \Sigma_{SS})$. From decomposition
\eqref{eq:decomposition_gram}, we obtain that 
\begin{equation}\label{eq:smax_decomposition}
s_{\max}(\Sigma_{SS}^* - \Sigma_{SS}) \leq (1 - \rho) s_{\max}
\left(\frac{1}{1 -  \rho} \frac{\wt{X}_S^{\T} \wt{X}_S}{n} - I \right) + 2 \sqrt{\rho}
s_{\max} \left(\frac{\wt{X}_S^{\T} \one_S}{n} \right) 
\end{equation}
Introduce $w = \left(\sum_{i=1}^n  \wt{x}_{i1}/n, \ldots, \sum_{i=1}^n \wt{x}_{is}/n \right)^{\T}$
as the vector of column means of $\wt{X}_S$. We have that
\begin{equation}\label{eq:smax_crosspart}
s_{\max}\left( \frac{\wt{X}_S^{\T} \one_S}{n} \right) = \sup_{\norm{u}_2 = 1}
\sup_{\norm{v}_2 = 1} u^{\T} \frac{\wt{X}_S^{\T} \one_S}{n} v =
 \sup_{\norm{u}_2 = 1}
\sup_{\norm{v}_2 = 1} u^{\T} w \bm{1}^{\T} v =
\sqrt{s} \nnorm{w}_2.  
\end{equation}
Moreover, 
\begin{equation}\label{eq:expandnorm}
\norm{w}_2^2 = \sum_{j = 1}^s \left(\frac{\su \wt{x}_{ij}}{n} \right)^2 =
\frac{1}{n} \sum_{j = 1}^s z_j^2, \quad  \text{where} \; \; z_j = n^{-1/2} \su
\wt{x}_{ij}.  
\end{equation}
Noting that the $\{ z_j \}_{j = 1}^s$ are i.i.d.~zero-mean sub-Gaussian
random variables with parameter $\sigma$ and variance no larger than one, we
are in position to apply Lemma \ref{lemma:subexponential}, which yields that
for any $t \geq 0$
\begin{equation}\label{eq:expandnorm_probability}
\p \left( \nnorm{w}_2^2 > \frac{s}{n} (1 + t)  \right) \leq \exp \left(-c
  \min \left(\frac{t^2}{\sigma^4}, \frac{t}{\sigma^2} \right) s \right).
\end{equation}
Combining \eqref{eq:smax_decomposition}, \eqref{eq:smax_crosspart} and     
\eqref{eq:expandnorm_probability} and using Lemma
\ref{lemma:singularvalues_general} to control the term \\$s_{\max}
\left(\frac{1}{1 -  \rho} \frac{\wt{X}_S^{\T} \wt{X}_S}{n} - I \right)$, we obtain that
for any $t \geq 0$ and any $z \geq 0$
\begin{align}\label{eq:control_smax_final}
\begin{split}
&\p \left( s_{\max}(\Sigma_{SS}^* - \Sigma_{SS}) > \max \left\{
     C \sqrt{\frac{s}{n}} + \frac{z}{\sqrt{n}}, \left(  C \sqrt{\frac{s}{n}} +
      \frac{z}{\sqrt{n}} \right)^2   \right \}  + 2
  \sqrt{\frac{s^2 (1 + t)}{n}} 
\right) \\
&\leq \exp(-c_1 \min\{t, t^2\} s) - 2 \exp(-c_2 z^2), 
\end{split}
\end{align}
where $C, c_1, c_2 > 0$ only depend on the sub-Gaussian parameter $\sigma$.
Equipped with these auxiliary results, we now turn to the proofs of 
Proposition 3 and 4. 

\subsection*{Proof of Proposition 3}

Let us first recall the restricted eigenvalue condition. 
\paragraph{Condition 2.}
Let $\mc{J}(s) = \{ J \subseteq \{1,\ldots,p\}: 1 \leq |J| \leq s \}$ and for $J
\in \mc{J}(s)$ and $\alpha \geq 1$,
\begin{equation*}
\mc{R}(J, \alpha) = \{\delta \in \R^p:\, \nnorm{\delta_{J^c}}_1 \leq  \alpha
\nnorm{\delta_J}_1 \}. 
\end{equation*}
We say that the design satisfies the $(\alpha, s)$-\textbf{restricted
  eigenvalue condition} if there exists a constant $\phi(\alpha, s)$ so
that 
\begin{equation}\label{eq:recondition}
\min_{J \in \mc{J}(s)} \, \min_{\delta \in \mc{R}(J, \alpha) \setminus
  \mathbf{0}} \, \, \frac{\delta^{\T} \Sigma \delta}{\nnorm{\delta_J}_2^2} \geq
\phi(\alpha, s) > 0. 
\end{equation}
The proof of Proposition 3 relies on a recent result in
[3]. In order to state that result, we need the following
preliminaries concerning $\psi_2$-random variables taken from [1] (see Definition 1.1.1 and
Theorem 1.1.5 therein). 
\begin{defn}
A random variable $Z$ is said to be $\psi_2$ with parameter $\theta > 0$ if 
\begin{equation}\label{eq:psi2}
\inf \left\{a > 0:\;  \E \left[\exp(Z^2 / a^2) \right] \leq e \right\} \leq
\theta. 
\end{equation}
\end{defn}
\begin{lemma}\label{lem:psi2fromtail} If a random variable $Z$ has the
  property that there exist positive constants $C, C'$ so that $\forall z \geq C'$
\begin{equation*}
\p \left(|Z| \geq z \right) \leq \exp \left(-z^2 / C^2 \right),
\end{equation*}
then $Z$ is $\psi_2$ with parameter no more than $2 \max(C, C')$. 
\end{lemma}
The following statement is essentially a special case of Theorem
1.6 in [3]. We state it in simplified form that is sufficient for
our purpose here.  
\begin{lemma}\label{lem:rudelsonzhou} Let $\Psi \in \R^{n \times p}$ be a matrix
whose rows $\Psi^1,\ldots,\Psi^n$, are independent random vectors that are
\begin{enumerate}
\item \textbf{isotropic}, i.e. $\E[\scp{\Psi^i}{u}^2] = 1$, $i=1,\ldots,n$,
\item $\bm{\psi_2}$, i.e. there exists $\theta > 0$ such that for every unit vector $u \in \R^p$
      \begin{equation}\label{eq:psi2marginals}   
      \inf \left\{a > 0:\;  \E \left[\exp(\scp{\Psi^i}{u}^2 / a^2) \right] \leq e \right\}
      \leq \theta, \; \; \; i=1,\ldots,n.
    \end{equation}
  \end{enumerate}
Let further $R \in \R^{p \times p}$ be a positive definite matrix with minimum eigenvalue $\vartheta > 0$ and set
$\Gamma = \frac{1}{n} R^{\T} \Psi^{\T} \Psi R$. Then, for any $\delta \in
(0,1)$ and any $\alpha \in [1, \infty)$, there exist positive constants $C_{\theta}$, $c > 0$ (the first
depending on the $\psi_2$ parameter $\theta$) so that if 
\begin{equation*}
n \geq \frac{C_{\theta}}{\delta^2} s \left(1 +  \frac{16 (3 \alpha^2)(3 \alpha
    + 1)}{\vartheta^2 \delta^2} \right) \log \left(c \frac{p}{s \delta} \right),
\end{equation*}
with probability at least $1 - 2 \exp(-\delta^2 n / C_{\theta})$, $\Gamma$
satisfies the $(\alpha, s)$-restricted eigenvalue condition with $\phi(\alpha, s) = \vartheta^2 (1 - \delta)^2$ .   
\end{lemma}

We now state and prove Proposition 3.
\paragraph{Proposition 3.} Let $X$ be a random matrix from $\text{Ens}_+$ \eqref{eq:Ensplus} scaled such that 
$\Sigma^* = \E[\frac{1}{n} X^{\T} X] = (1 - \rho) I
+ \rho \bm{1} \bm{1}^{\T}$ for some $\rho \in (0,1)$. Set $\delta \in (0,1)$. There exists
constants $C, c > 0$ depending only on $\delta$, $\rho$ and the sub-Gaussian parameter of
the centered entries of $X$ so that if $n \geq C \, s \log (p \vee n)$, then, with probability at least 
$1 - \exp(-c \delta^2 n ) - 6/(p \vee n)$,
$\Sigma = X^{\T} X / n$  has the self-regularizing property with $\tau^2 =
\rho/2$ and satisfies the $(3/\tau^2, s)$ restricted eigenvalue condition of Theorem 2 with $\phi(3/\tau^2, s) = (1 - \rho)(1 - \delta)^2$. 

\begin{bew} We first show that $\Sigma$ satisfies the self-regularizing property with $\tau^2 \geq \rho/2$ with probability at least $1 - 6/(p \vee
  n)$. According to Eq.(6.4) in the paper, we have
\begin{equation*}
 \tau_0^2 = \min_{\lambda \in T^{p-1}} \lambda^{\T} \Sigma \lambda \geq
\min_{\lambda \in T^{p-1}} \lambda^{\T} \Sigma^* \lambda - \max_{\lambda \in T^{p-1}} \lambda^{\T} \left( \Sigma^* -
\Sigma \right) \lambda \geq \rho - \max_{j,k} \left|\left(\Sigma - \Sigma^* \right)_{jk}\right|.
\end{equation*}
Consequently, in virtue of \eqref{deviation:gram}, there exists a numerical constant $C'$
depending on $\sigma$ and $\rho$ only so that if $n \geq C' \log(p \vee n)$, $\tau_0^2 \geq \frac{1}{2}
\rho$ with the probability as claimed. In the sequel, it will be shown that
conditional on the event $\{\tau_0^2 \geq \rho/2 \}$, Lemma \ref{lem:rudelsonzhou} can be applied with 
\begin{equation*}
\Gamma = \Sigma, \; \; \; \, R = (\Sigma^*)^{1/2}, \; \; \; \, \Psi = X
(\Sigma^*)^{-1/2}, \; \; \; \, \vartheta^2 = 1 - \rho, \; \; \; \, \alpha =
\frac{3}{\tau^2} \leq \frac{6}{\rho}, \; \; \; \, \theta = C_{\sigma, \rho},
\end{equation*}
where $(\Sigma^*)^{1/2}$ is the root of $\Sigma^*$ and $C_{\sigma, \rho}$ is a
constant depending only on $\sigma$ and $\rho$. By construction, $\Psi = X (\Sigma^*)^{-1/2}$ has independent isotropic rows. 
It remains to establish that the rows satisfy condition \eqref{eq:psi2marginals} of Lemma \ref{lem:rudelsonzhou}. Since the rows of $\Psi$ are i.i.d., it
suffices to consider a single row. Let us write $X^1$ for the transpose of the first row of $X$ and accordingly
$\Psi^1 = (\Sigma^*)^{-1/2} \,X^1$ for the transpose of the first row of
$\Psi$. Furthermore, we make use of the decomposition $X^1 = \wt{X}^1 + \sqrt{\rho} \bm{1}$, where
the entries of $\wt{X}^1$ are i.i.d~zero-mean sub-Gaussian
random variables with parameter $\sigma$ (cf.~\eqref{eq:decomposition_Ensplus_specific}). We then have for
any unit vector $u$ 
\begin{align*}
\scp{\Psi^1}{u} = \scp{(\Sigma^*)^{-1/2} \,X^1}{u} &= \scp{(\Sigma^*)^{-1/2}\,(\wt{X}^1 + \sqrt{\rho} \bm{1})}{u} \\
&=\scp{\wt{X}^1}{(\Sigma^*)^{-1/2} \, u} + \sqrt{\frac{\rho}{(1 -
    \rho) +  p \rho}} \scp{\bm{1}}{u} \\
&\leq \scp{\wt{X}^1}{(\Sigma^*)^{-1/2} \, u} + 1.             
\end{align*}
For the second equality, we have used that $\bm{1}$ is an eigenvector of
$\Sigma^*$ with eigenvalue $1 + (p -1) \rho$, while the inequality results
from Cauchy-Schwarz. We now estimate the moment-generating function of the
random variable $\scp{\Psi^1}{u}$ as follows. For any $t \geq 0$, we have
\begin{align*}
\E[\exp(t \scp{\Psi^1}{u})] &\leq \exp(t) \E \left[\exp \left(t
\scp{\wt{X}^1}{(\Sigma^*)^{-1/2} \, u} \right) \right]\\
 &\leq \exp(t) \E \left[\exp \left(\frac{\sigma^2
       t^2}{2}\nnorm{(\Sigma^*)^{-1/2} \, u}_2^2 \right) \right] \\
&\leq \exp(t) \exp \left(\frac{\sigma^2
       t^2}{2 (1 - \rho)} \right) \\
&\leq e  \exp \left(\frac{(\sigma^2 +2)
       t^2}{2 (1 - \rho)} \right)  = e \exp \left(\frac{\wt{\sigma}^2
       t^2}{2} \right),
\end{align*} 
where $\wt{\sigma} = \sqrt{(\sigma^2 + 2)/(1 - \rho)}$. For the third
equality, we have used that the maximum eigenvalue of $(\Sigma^*)^{-1}$ equals $(1 -
\rho)^{-1}$. Analogously, we obtain that 
\begin{equation*}
-\scp{\Psi^1}{u} \leq \scp{-\wt{X}^1}{(\Sigma^*)^{-1/2} \, u} + 1, \quad \; 
\text{and} \; \; \; \E[\exp(t \scp{-\Psi^1}{u})] \leq e \exp \left(\frac{\wt{\sigma}^2
       t^2}{2} \right) \; \forall t \geq 0.
\end{equation*}
From the Chernov bound, we hence obtain that for any $z \geq 0$
\begin{equation*}
\p(|\scp{\Psi^1}{u}| > z) \leq 2 e \exp \left(-\frac{z^2}{2 \wt{\sigma}^2} \right).
\end{equation*}
Invoking Lemma \ref{lem:psi2fromtail} with $C' = \wt{\sigma} \sqrt{3 \log(2
  e)}$ and $C = \sqrt{6} \wt{\sigma}$, it follows that the random
variable $\scp{\Psi^1}{u}$ is $\psi_2$ with parameter $2 \sqrt{6} \wt{\sigma}
\invcoloneq C_{\sigma, \rho}$, and we conclude that 
the rows of the matrix $\Psi$ indeed satisfy condition \eqref{eq:psi2marginals} with
$\theta$ equal to that value of the parameter.  
\end{bew}

\subsection*{Proof of Proposition 4}

\paragraph{Proposition 4.}
Let $X$ be a random matrix from $\text{Ens}_+$ \eqref{eq:Ensplus} scaled such that 
$\Sigma^* = \E[\frac{1}{n} X^{\T} X] = (1 - \rho) I
+ \rho \bm{1} \bm{1}^{\T}$ for some $\rho \in (0,1)$. Fix $S
\subset \{1,\ldots,p\}$, $|S| \leq s$. Then there exists constants $c, c', C,
C' > 0$ depending only on $\rho$ and the sub-Gaussian parameter of the
centered entries of $X$ such that for all $n \geq C s^2 \log (p \vee n)$,  
\begin{equation*}
  \tau^2(S) \geq c s^{-1} - C' \sqrt{\frac{\log p}{n}}
\end{equation*}
with probability no less than $1 - 6/(p \vee n) -  3\exp(-c' (s \vee \log n))$.

\begin{bew} The scaling of $\tau^2(S)$ is analyzed based on representation
\begin{equation}\label{eq:tauZ}
\tau^2(S) = \min_{\theta \in \R^s, \; \lambda \in T^{p-s-1}} \frac{1}{n}
\norm{X_S \theta -
X_{S^c} \lambda}_2^2.
\end{equation}
In the following, denote by $\mathbb{S}^{s-1} = \{ u \in
\R^s: \norm{u}_2 = 1 \}$ the unit sphere in $\R^s$. Expanding the square in 
\eqref{eq:tauZ}, we have
\begin{align}\label{eq:lowerboundtauS}
\begin{split}
\tau^2(S) &= \min_{\theta \in \R^s, \; \lambda \in T^{p - s - 1}} \; \theta^{\T} 
\Sigma_{SS} \theta - 2 \theta^{\T} \Sigma_{SS^c}
\lambda + \lambda^{\T} \Sigma_{S^c S^c} \lambda \\
&\geq \min_{r \geq 0, \; u \in \mathbb{S}^{s-1}, \; \lambda \in T^{p - s - 1}} 
\; r^2 u^{\T} \Sigma_{SS}^{\ast} u - r^2 s_{\max} \left(\Sigma_{SS} - \Sigma_{SS}^{\ast} \right) - \\
&\qquad - 2 r u^{\T}
       \Sigma_{SS^c} \lambda + \lambda^{\T}
        \Sigma_{S^c S^c} \lambda \\
&\geq  \min_{r \geq 0, \; u \in \mathbb{S}^{s-1}, \; \lambda \in T^{p - s - 1}} 
\; r^2 u^{\T} \Sigma_{SS}^{\ast} u 
    - r^2 s_{\max} \left(\Sigma_{SS} - \Sigma_{SS}^{\ast} \right) \\ 
&\qquad   - 2 \rho r u^{\T} \bm{1} - 2 r u^{\T} (\Sigma_{SS^c} - 
    \Sigma_{SS^c}^{\ast}) \lambda +  \rho + \frac{1 - \rho}{p - s}  - \\
&\qquad - \max_{\lambda \in T^{p - s -1}} \left|\lambda^{\T}
    (\Sigma_{S^c S^c} - \Sigma_{S^cS^c}^{\ast}) \lambda \right|. 
\end{split}
\end{align}
For the last inequality, we have used that $\min_{\lambda \in T^{p - s - 1}}
\lambda^{\T} \Sigma_{S^c S^c}^{\ast} \lambda = \rho + \frac{1 - \rho}{p - s}$. We further set 
\begin{align}
\Delta &= s_{\max} \left(\Sigma_{SS} - \Sigma_{SS}^{\ast} \right), \label{eq:prop4:Delta}\\ 
\delta &= \max_{u \in \mathbb{S}^{s
    - 1}, \lambda \in T^{p - s - 1}} \left|u^{\T}  \left(\Sigma_{S^cS^c} -
    \Sigma_{S^cS^c}^{\ast} \right) \lambda \right|. \label{eq:prop4:delta} 
\end{align}
The random terms $\Delta$ and $\delta$ will be controlled uniformly over $u
\in \mathbb{S}^{s - 1}$ and $\lambda \in T^{p - s - 1}$ below by invoking
\eqref{deviation:gram} and \eqref{eq:control_smax_final}. For the moment, we
treat these two terms as constants. We now minimize the lower bound in \eqref{eq:lowerboundtauS} w.r.t. $u$ and $r$ separately from $\lambda$. This minimization
problem involving $u$ and $r$ only reads
\begin{equation}\label{eq:runfixed}
\min_{r \geq 0, \; u \in \mathbb{S}^{s-1}} r^2 u^{\T} \Sigma_{SS}^{\ast}
u  - 2 \rho r u^{\T}
\bm{1} - r^2 \Delta - 2 r \delta.     
\end{equation}
We first derive an expression for 
\begin{equation}\label{eq:phir}
\phi(r) = \min_{u \in \mathbb{S}^{s-1}} r^2 u^{\T} \Sigma_{SS}^{\ast}
u  - 2 \rho r u^{\T} \bm{1}.
\end{equation} 
We decompose $u = u^{\parallel} + u^{\perp}$, where $u^{\parallel} = \scp{\frac{\bm{1}}{\sqrt{s}}}{u}
\frac{\bm{1}}{\sqrt{s}}$ is the projection of $u$ on the unit vector
$\bm{1}/\sqrt{s}$, which is an eigenvector of $\Sigma_{SS}^{\ast}$ associated with its largest
eigenvalue $1 + \rho (s - 1)$. By Parseval's identity, we have $\nnorm{u^{\parallel}}_2^2 = \gamma$, 
$\nnorm{u^{\perp}}_2^2 = (1 - \gamma)$ for some $\gamma \in [0,1]$. Inserting
this decomposition into \eqref{eq:phir} and noting that the remaining eigenvalues of
$\Sigma_{SS}^{\ast}$ are all equal to $(1 - \rho)$, we obtain that 
\begin{align}\label{eq:functionofgamma}
\begin{split}
&\phi(r) = \min_{\gamma \in [0,1]} \Phi(\gamma, r),\\ 
&\text{with} \; \, \Phi(\gamma, r) = r^2 \gamma \underbrace{(1 + (s - 1) \rho)}_{s_{\max}(\Sigma_{SS}^{\ast})} + r^2 (1 -
\gamma)\underbrace{(1 - \rho)}_{s_{\min}(\Sigma_{SS}^{\ast})}  - 2 \rho r
\sqrt{\gamma} \sqrt{s},
\end{split}
\end{align}
where we have used that $\nscp{u^{\perp}}{\bm{1}} = 0$. Let us put aside the
constraint $\gamma \in [0,1]$ for a moment. The function
$\Phi$ in \eqref{eq:functionofgamma} is a convex function of $\gamma$, hence we may find an (unconstrained) minimizer $\wt{\gamma}$ by differentiating
and setting the derivative equal to zero. This yields $\wt{\gamma} =
\frac{1}{r^2 s}$, which coincides with the constrained minimizer if and only
if $r \geq \frac{1}{\sqrt{s}}$. Otherwise, $\wt{\gamma} \in \{0,1 \}$. We can 
rule out the case $\wt{\gamma} = 0$, since for all $r < 1 / \sqrt{s}$ 
\begin{equation*}
\Phi(0, r) = r^2 (1 - \rho) >  r^2 (1 + (s - 1) \rho) - 2 \rho r \sqrt{s} = \Phi(1, r). 
\end{equation*}
We have $\Phi(\frac{1}{r^2 s}, r) = r^2(1 - \rho) - \rho$ and
$\Phi(\frac{1}{r^2 s}, \frac{1}{\sqrt{s}}) = \Phi(1,
\frac{1}{\sqrt{s}})$. Hence, the function $\phi(r)$ in \eqref{eq:phir} is given by
\begin{equation}\label{eq:phirfound}
\phi(r) = \begin{cases}
  r^2 s_{\max}(\Sigma_{SS}^{\ast}) - 2 \rho r \sqrt{s} \quad& r \leq 1/\sqrt{s}, \\
          r^2 (1 - \rho) - \rho  \quad&\text{otherwise}.  
          \end{cases}
\end{equation}
The minimization problem \eqref{eq:runfixed} to be considered eventually reads
\begin{equation}\label{eq:psir}
\min_{r \geq 0} \psi(r), \quad \text{where} \; \; \psi(r) = \phi(r) - r^2 \Delta - 2 r \delta.
\end{equation}
We argue that it suffices to consider the case $r \leq 1/\sqrt{s}$ in
\eqref{eq:phirfound} provided
\begin{equation}\label{eq:crucialcondition1}
((1 - \rho) - \Delta) > \delta \sqrt{s}, 
\end{equation}
a condition we will comment on below. If this condition is met,
differentiating shows that $\psi$ is increasing on $(\frac{1}{\sqrt{s}},
\infty)$. In fact, for all $r$ in that interval,  
\begin{equation*}
\frac{d}{d r} \psi(r) = 2 r (1 - \rho) - 2 r \Delta - 2 \delta, \; \text{and
  thus}
\end{equation*}
\begin{equation*}
\frac{d}{d r} \psi(r) > 0 \; \text{for all} \; r \in \; \left(\frac{1}{\sqrt{s}},
\infty \right) \; \Leftrightarrow \; \frac{1}{\sqrt{s}} ((1 - \rho) - \Delta) > \delta.   
\end{equation*}
Considering the case $r \leq 1/\sqrt{s}$, we observe that $\psi(r)$ is convex provided
\begin{equation}\label{eq:crucialcondition2}
s_{\max}(\Sigma_{SS}^{\ast}) > \Delta,
\end{equation}
a condition we shall comment on below as well. Provided
\eqref{eq:crucialcondition1} and \eqref{eq:crucialcondition2} hold true,
differentiating \eqref{eq:psir} and setting the result equal to zero, we
obtain that the minimizer $\wh{r}$ of \eqref{eq:psir} is given by $(\rho \sqrt{s} + \delta) /
(s_{\max}(\Sigma_{SS}^{\ast}) - \Delta)$. Substituting this result back into
\eqref{eq:psir} and in turn into the lower bound \eqref{eq:lowerboundtauS},
one obtains after collecting terms  
\begin{align}\label{eq:lowerboundtauS_2}
\begin{split}
\tau^2(S) \geq& \rho \frac{(1 - \rho) - \Delta}{(1 - \rho) + s \rho - \Delta} -
\frac{2 \rho \sqrt{s} \delta + \delta^2}{s_{\max}(\Sigma_{SS}^{\ast}) - \Delta}  + \frac{1 -
  \rho}{p - s} -\\
&- \max_{\lambda \in T^{p - s -1}} \left|\lambda^{\T}
    (\Sigma_{S^cS^c} - \Sigma_{S^cS^c}^{\ast}) \lambda
  \right|. 
\end{split}
\end{align}
In order to control $\Delta$ \eqref{eq:prop4:Delta}, we apply \eqref{eq:control_smax_final} with the choices
\begin{equation*}
z = \sqrt{s \vee \log n}, \quad \text{and} \; \; \; t = 1 \vee \frac{\log n}{s}.
\end{equation*}
Consequently, there exists a constant $C_1 > 0$ depending only on $\sigma$ so
that if \\$n > C_1 (s \vee \log n)$, we have that 
\begin{align}\label{eq:eventA}
\begin{split}
\p(\mc{A}) &\geq 1 - 3 \exp(-c' (s \vee \log n)), \\
&\quad \text{where} \; \; \mc{A} = \left \{\Delta \leq 2 \sqrt{\frac{s^2 (1 + 1 \vee
        (\log(n)/s))}{n}} + C' \sqrt{\frac{s \vee \log n}{n}}    \right\}
\end{split}
\end{align}
In order to control $\delta$ \eqref{eq:prop4:delta} and the last term in \eqref{eq:lowerboundtauS_2},
we make use of \eqref{deviation:gram}, which yields that 
\begin{align}\label{eq:eventB}
\begin{split}
&\p(\mc{B}) \geq 1 - \frac{6}{p \vee n}, \; \; \text{where} \; \; \\
&\mc{B} = \left\{ \delta
\leq C \sqrt{\frac{s \log (p \vee n)}{n}} \right\} \cap  \left\{ \sup_{\lambda \in T^{p - s -1}} \left|\lambda^{\T}
    (\Sigma_{S^cS^c} - \Sigma_{S^cS^c}^{\ast}) \lambda
  \right| \leq C \sqrt{\frac{ \log (p \vee n)}{n}} \right \}.
\end{split}
\end{align}
For the remainder of the proof, we work conditional on the two events $\mc{A}$
and $\mc{B}$. In view of \eqref{eq:eventA} and \eqref{eq:eventB}, we first
note that there exists $C_2 > 0$ depending only on $\sigma$ and $\rho$ such that if $n \geq
C_2 s^2 \log(p \vee n)$ the two conditions \eqref{eq:crucialcondition1} and \eqref{eq:crucialcondition2} supposed to be fulfilled
previously indeed hold. To conclude the proof, we re-write \eqref{eq:lowerboundtauS_2} as
\begin{align}\label{eq:lowerboundtauS_3}
\begin{split}
  &\tau^2(S) \geq  \frac{\rho (1 - \Delta/(1 - \rho))}{(1 - \Delta/(1 - \rho)) + s
  \frac{\rho}{1 - \rho}} + \frac{2 \rho \frac{\sqrt{s}}{1 + (s - 1) \rho}
  \delta}{1 - \Delta/(1 + (s - 1) \rho)} - \frac{\delta^2/(1 + (s-1) \rho)}{1
  - \Delta/(1 + (s - 1) \rho)} -  \\
&\qquad -\max_{\lambda \in T^{p - s -1}} \left|\lambda^{\T}
   (\Sigma_{S^c S^c} - \Sigma_{S^cS^c}^{\ast}) \lambda \right|.
\end{split}
\end{align}
Conditional on $\mc{A} \cap \mc{B}$, there exists $C_3 > 0$ depending only on
$\sigma$ and $\rho$ such that if $n \geq
 C_3 (s^2 \vee (s \log n)) $, when inserting the resulting scalings separately for each summand in
\eqref{eq:lowerboundtauS_3}, we have that 
\begin{align}\label{eq:tauSfinal}
\begin{split}
&c_1 s^{-1} - C_4 \sqrt{\frac{\log (p \vee n)}{n}} - C_5 \frac{\log (p \vee n)}{n} - C_6
\sqrt{\frac{\log (p \vee n)}{n}}\\ 
&=  c_1 s^{-1} - C_7 \sqrt{\frac{\log (p \vee
    n)}{n}}.
\end{split}
\end{align}
We conclude that if $n \geq \max\{C_1, C_2, C_3 \} s^2 \log(p \vee n) $,
\eqref{eq:tauSfinal} holds with probability no less than $1 - \frac{6}{p \vee
  n} - 3 \exp(-c' (s \vee \log n))$.


\end{bew}

\section*{Empirical scaling of $\tau^2(S)$ for $\text{Ens}_+$}

In Section 6.3, we have empirically investigated the scaling of $\tau^2(S)$
for the class \eqref{eq:Ensplus} in a high-dimensional setting for the following designs.
\begin{itemize}
\item[$E_1$:] $\{ x_{ij} \} \overset{\text{i.i.d.}}{\sim} \; a \, \text{uniform}([0,1/\sqrt{3 \cdot a}]) + (1 - a) \delta_0$, $a \in \{1, \frac{2}{3},
  \frac{1}{3}, \frac{2}{15} \} \,$ ($\rho \in \{\frac{3}{4}, \frac{1}{2},
  \frac{1}{3}, \frac{1}{10} \}$)
\item[$E_2$:] $\{ x_{ij} \} \overset{\text{i.i.d.}}{\sim}
  \; \frac{1}{\sqrt{\pi}}\, \text{Bernoulli}(\pi), \; \pi \in
  \{\frac{1}{10}, \frac{1}{4}, \frac{1}{2}, \frac{3}{4}, \frac{9}{10}  \} \,$ 
  ($\rho \in \{\frac{1}{10}, \frac{1}{4}, \frac{1}{2}, \frac{3}{4},
  \frac{9}{10} \}$)
\item[$E_3$:] $\{ x_{ij} \} \overset{\text{i.i.d.}}{\sim} \, |Z|, \, Z \sim
  a \, \text{Gaussian}(0, 1) + (1 - a) \delta_0$, $a \in
  \{1,\frac{\pi}{4},\frac{\pi}{8},\frac{\pi}{20}\} \,$ ($\rho \in \{\frac{2}{\pi}, \frac{1}{2},
  \frac{1}{4}, \frac{1}{10} \}$)
\item[$E_4$:] $\{ x_{ij} \} \overset{\text{i.i.d.}}{\sim} \, a
  \text{Poisson}(3/\sqrt{12 \,a}) + (1 - a) \delta_0$, $a \in
  \{1,\frac{2}{3},\frac{1}{3},\frac{2}{15}\} \,$ ($\rho \in \{\frac{3}{4}, \frac{1}{2},
  \frac{1}{4}, \frac{1}{10} \}$)
\end{itemize}
The results for $E_1$ are presented in the paper, and the 
results for $E_2$ to $E_4$ are displayed below.

\subsection*{$E_2$}

\begin{center}
\includegraphics[height = 0.17\textheight]{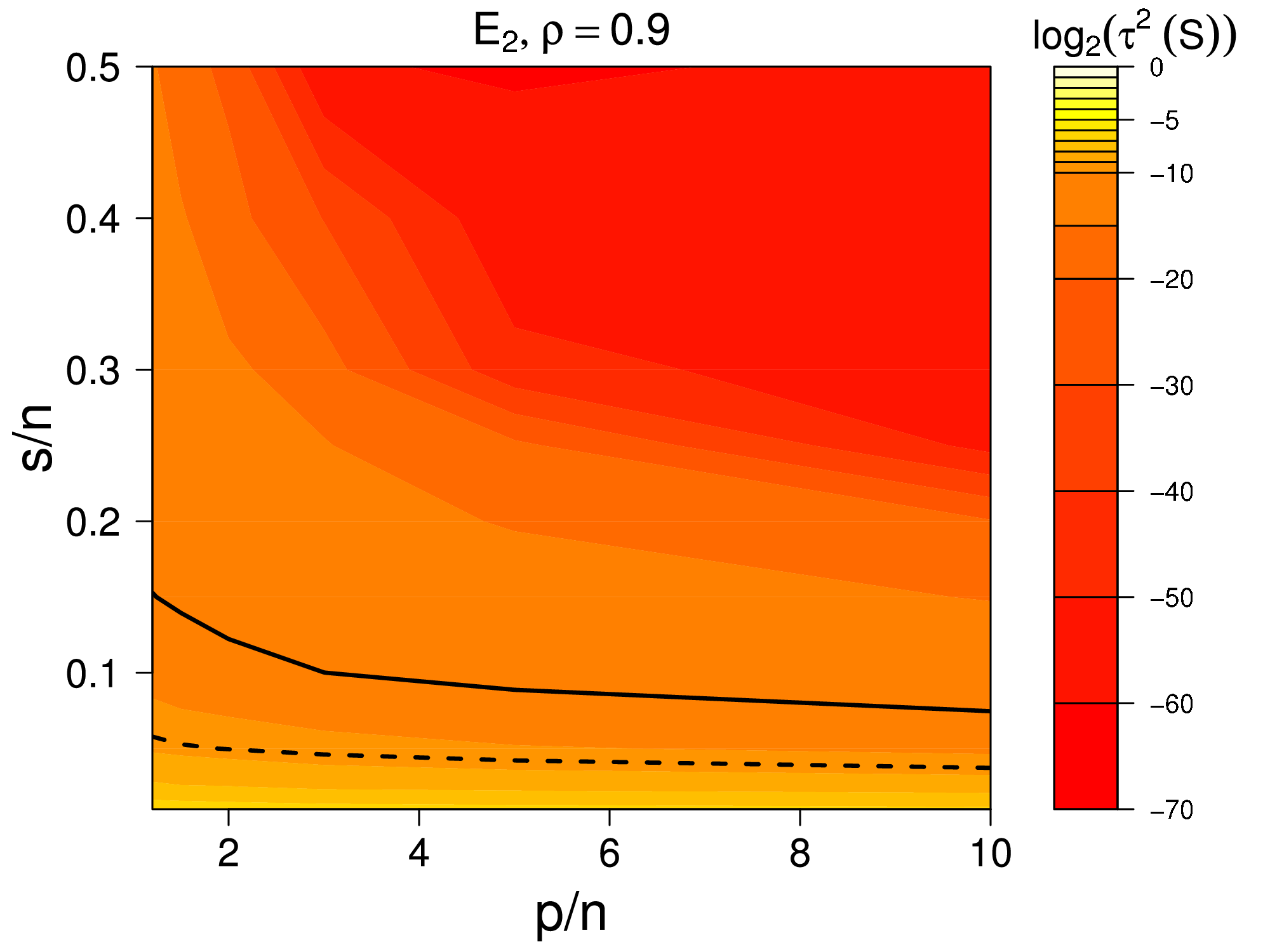}
\end{center}
\hfill\\
\hfill\\
\hfill\\
\begin{tabular}{ll}
 \hspace{.35cm}\includegraphics[height =
0.17\textheight]{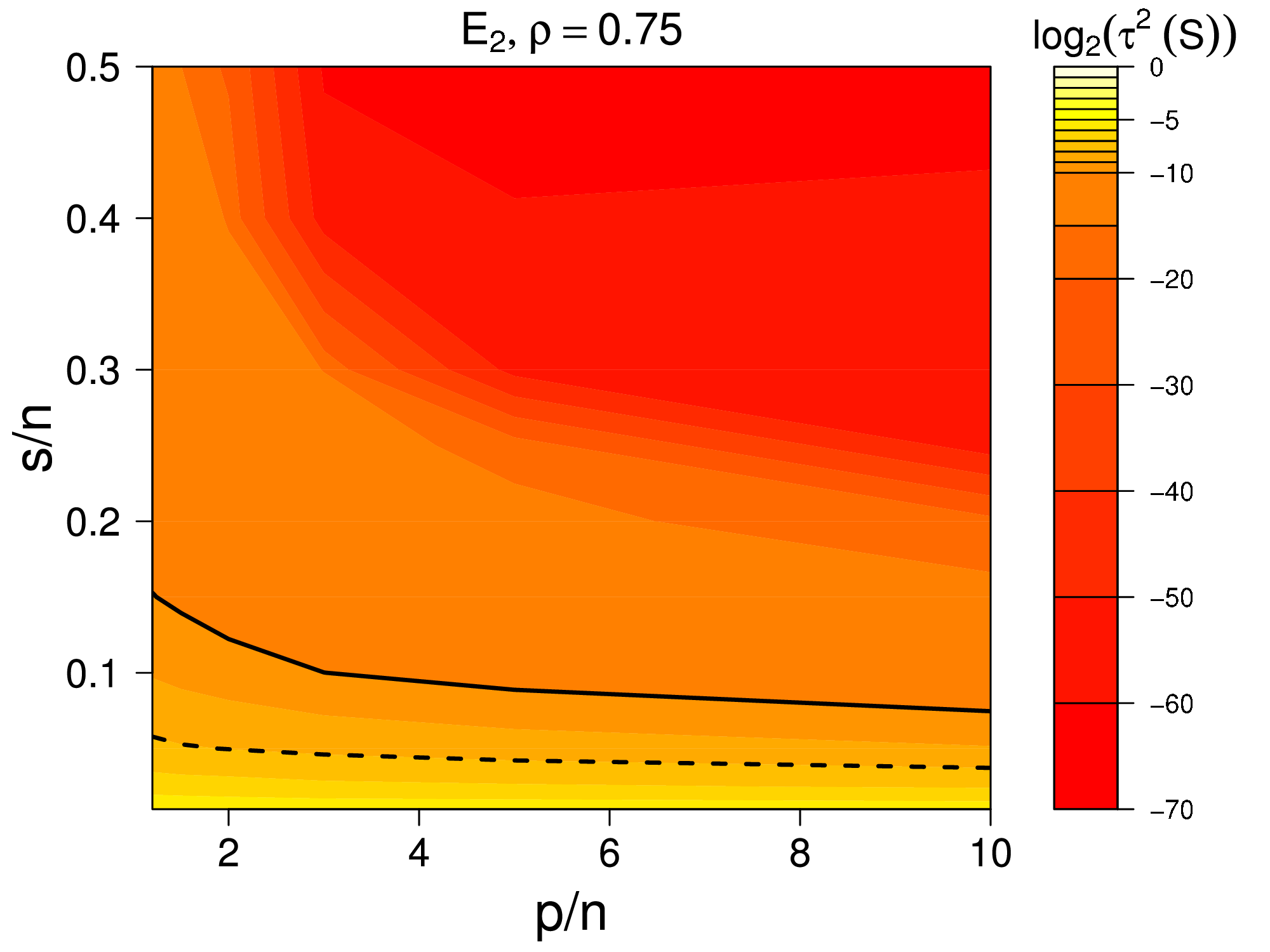}
& \hspace{2cm}\includegraphics[height =
0.17\textheight]{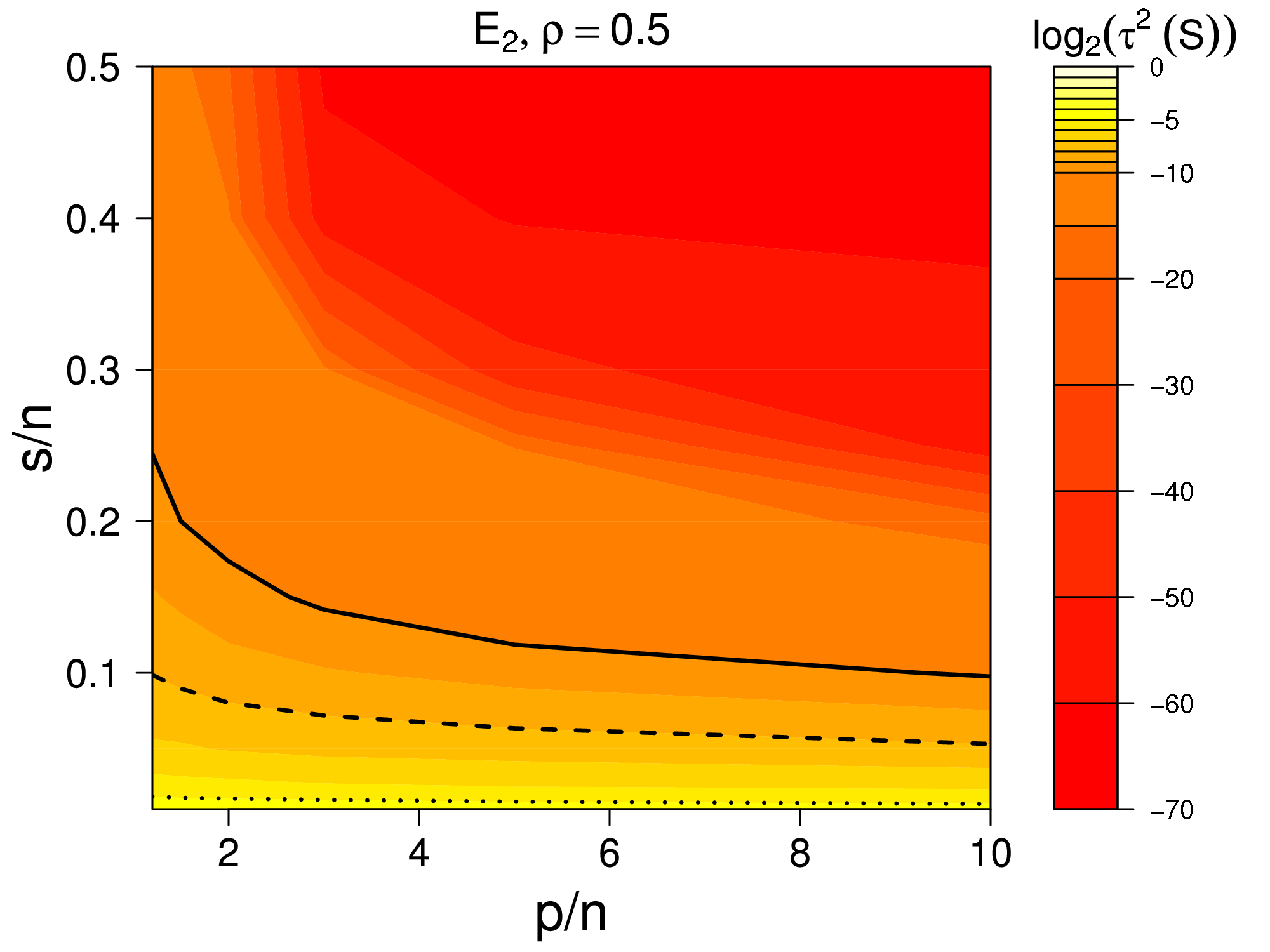} \\
  \\
\\
\\
 \hspace{.35cm}\includegraphics[height =
0.17\textheight]{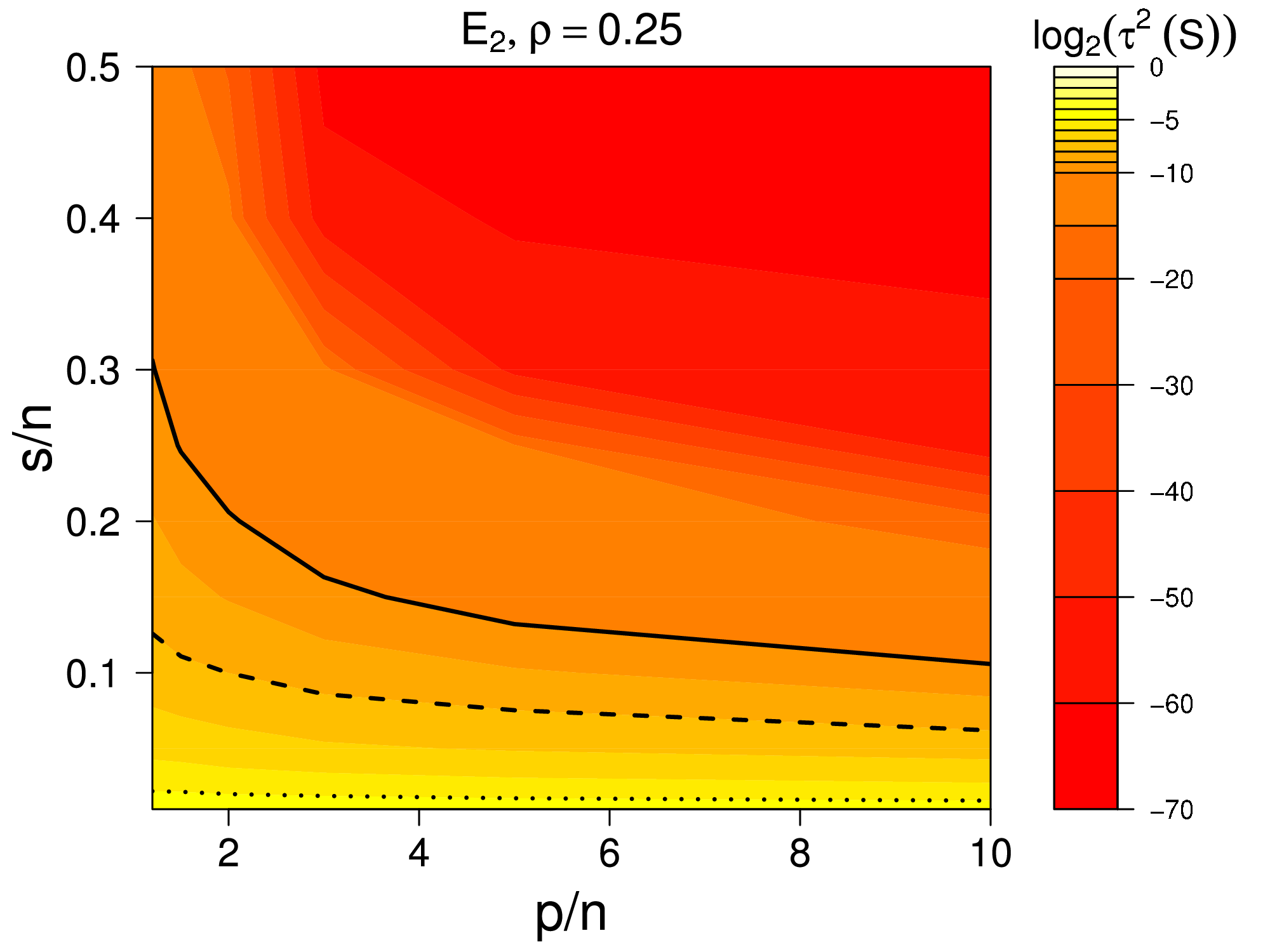} & \hspace{2cm}\includegraphics[height = 0.17\textheight]{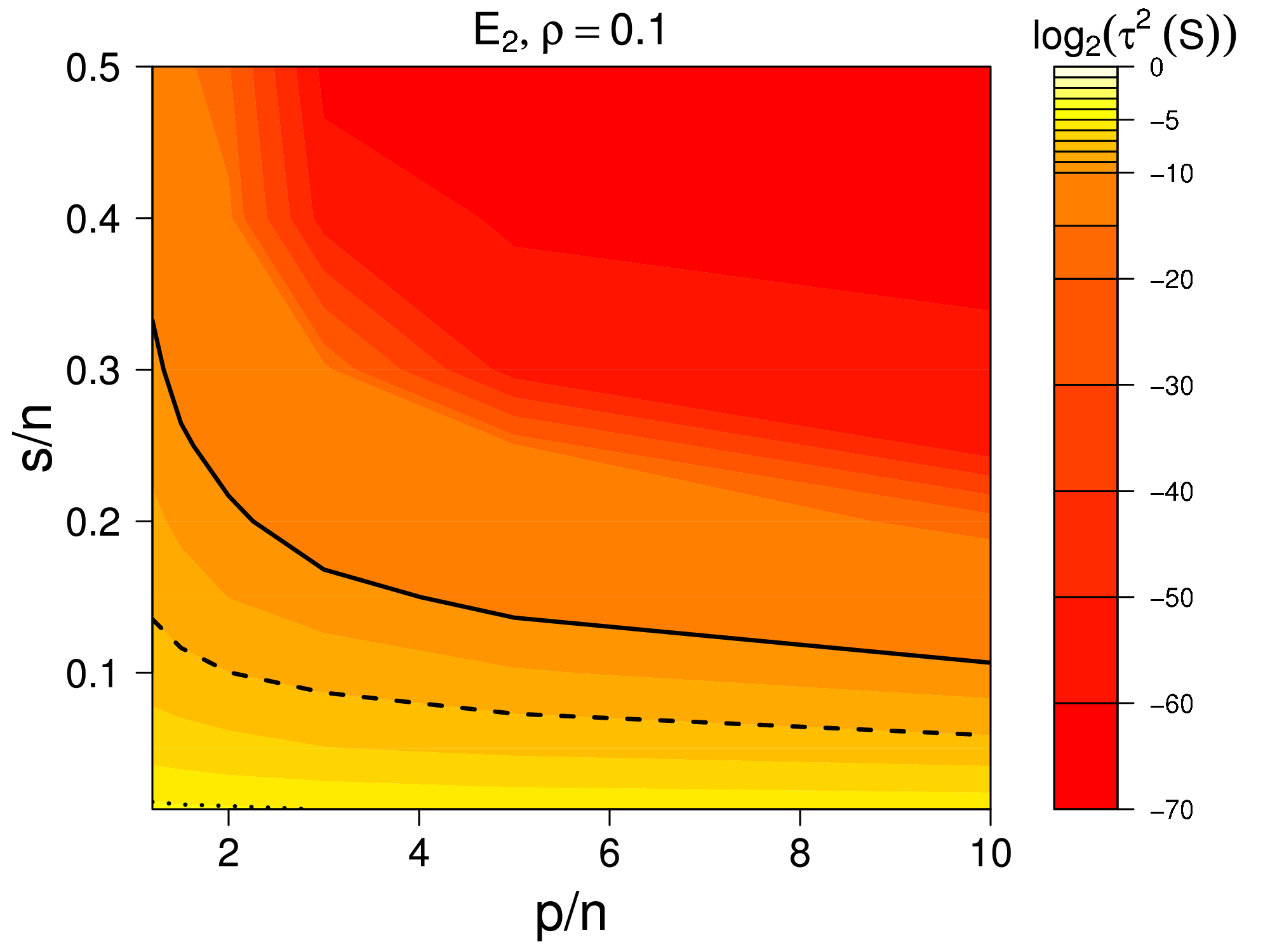}
\end{tabular}
\hfill\\
\hfill\\
\hfill\\

\subsection*{$E_3$}

\begin{tabular}{ll}
 \hspace{.35cm}\includegraphics[height =
0.17\textheight]{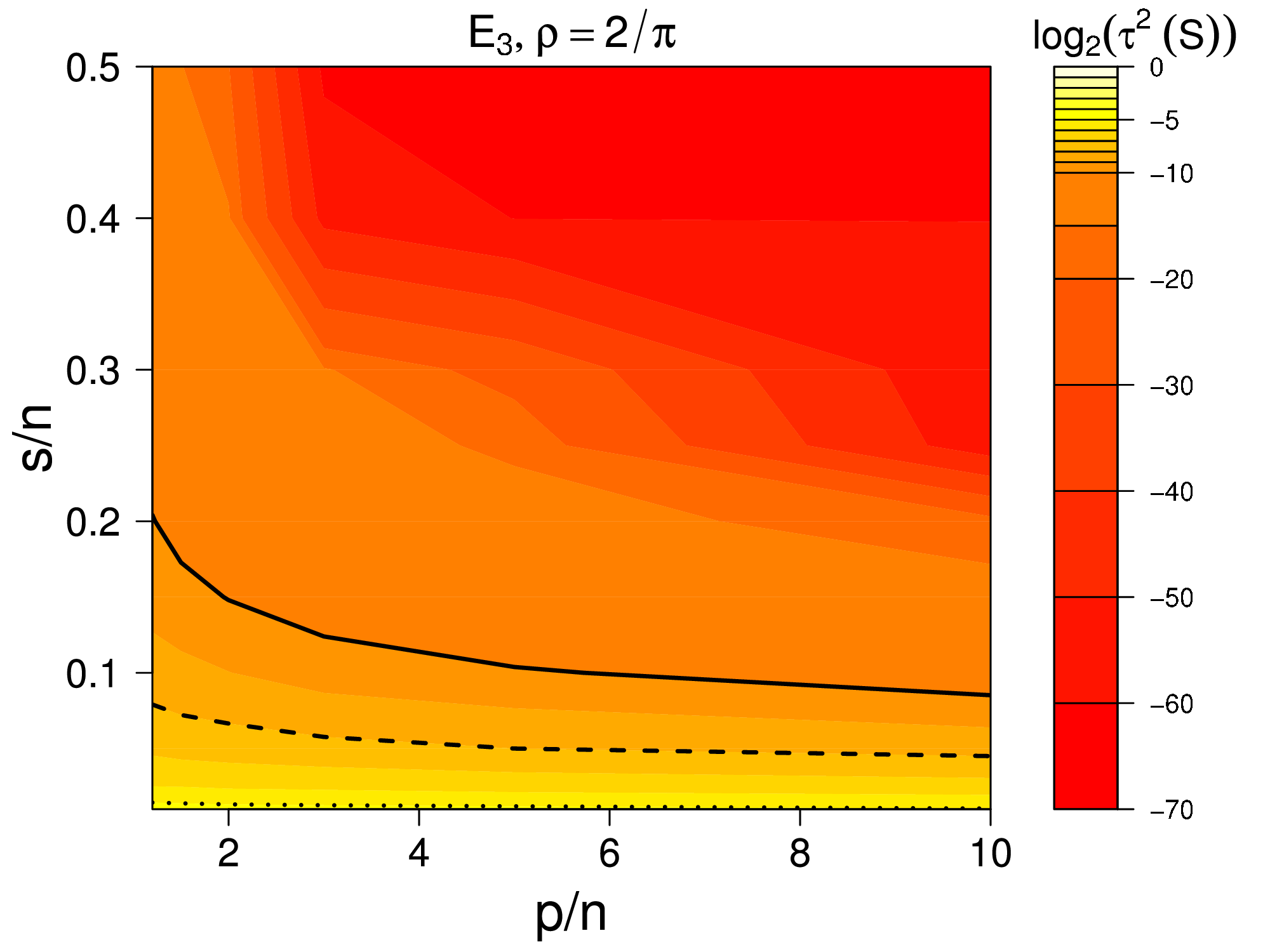}
& \hspace{2cm}\includegraphics[height =
0.17\textheight]{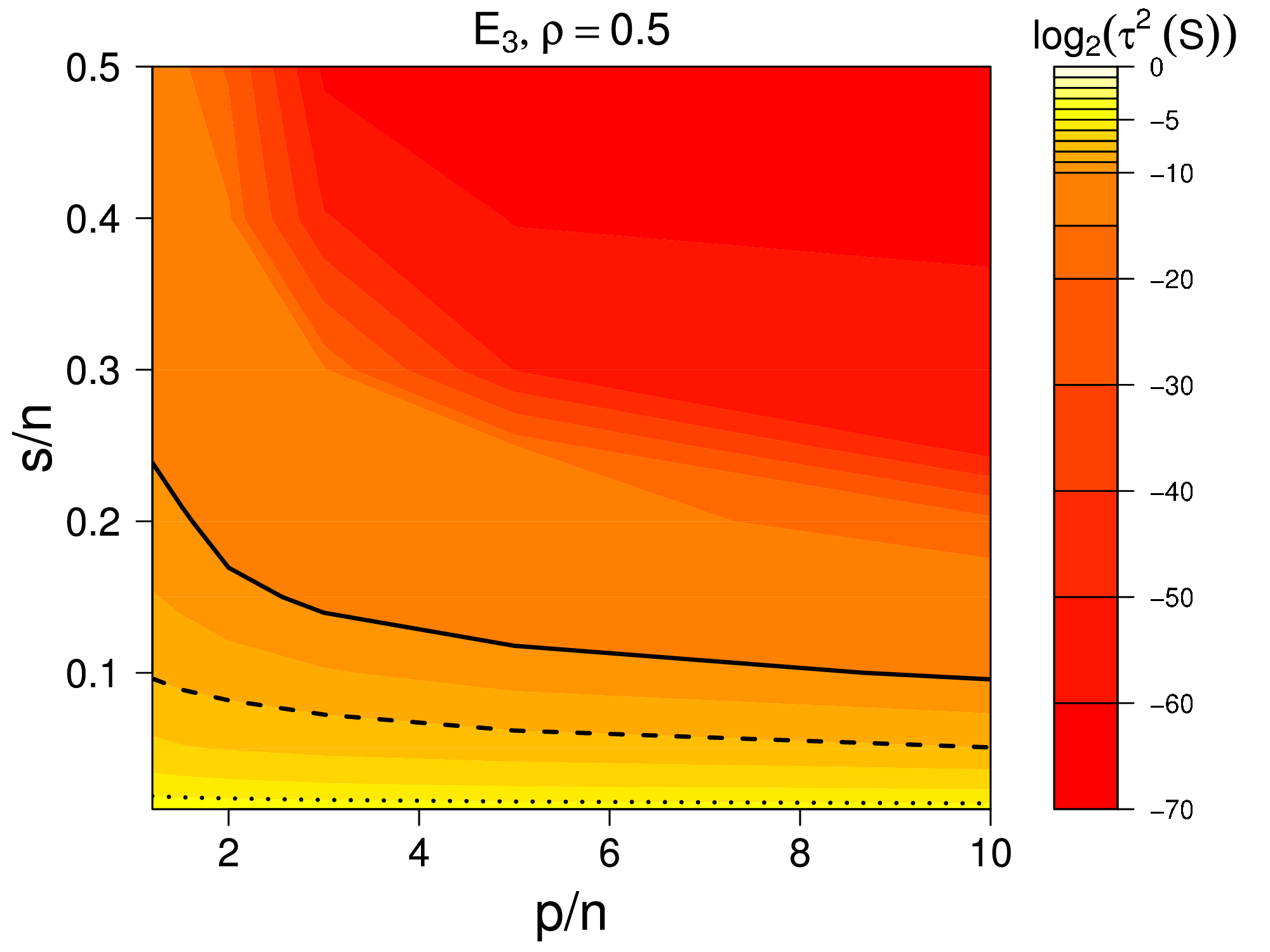} \\
  \\
\\
\\
 \hspace{.35cm}\includegraphics[height =
0.17\textheight]{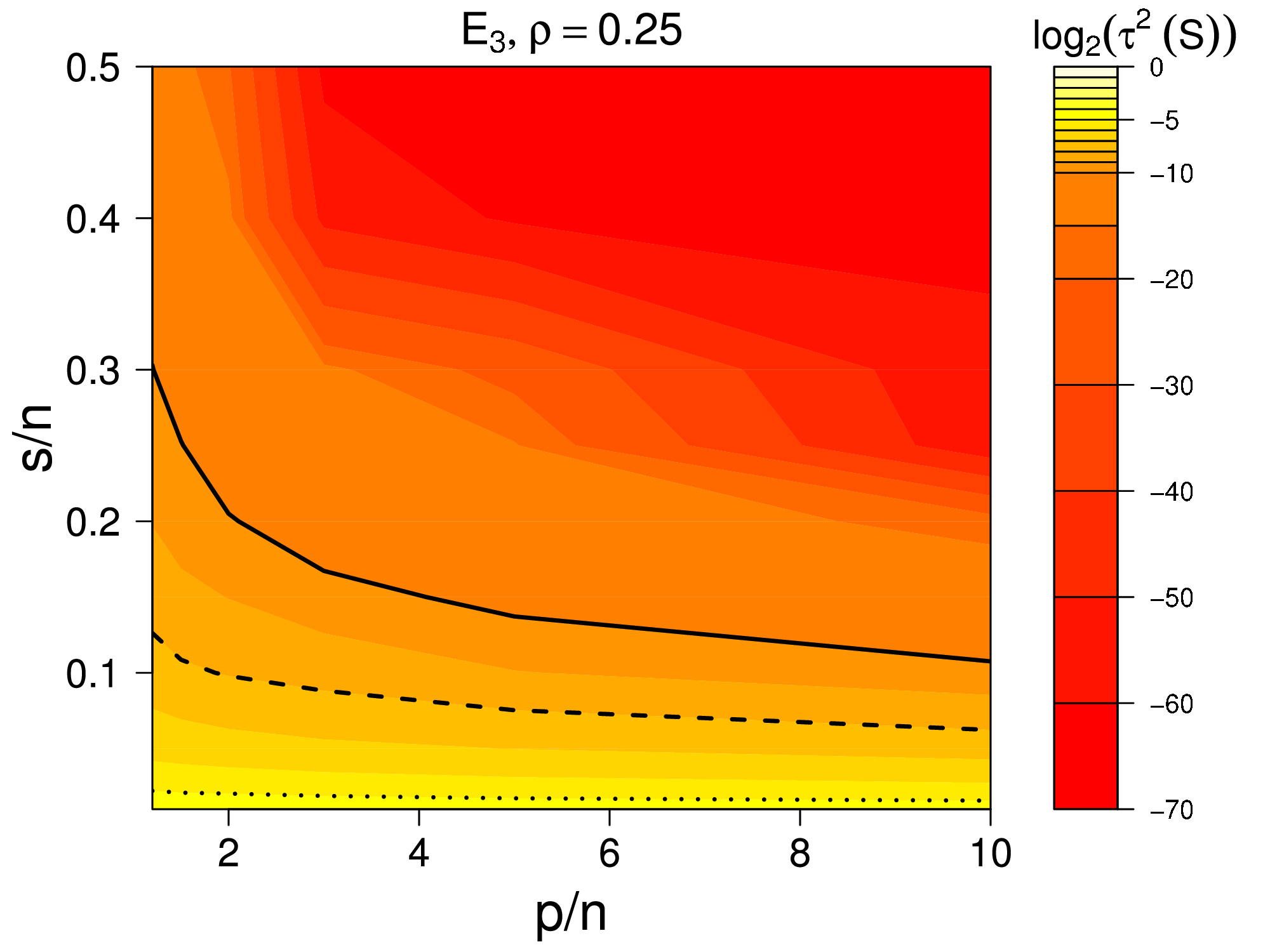} & \hspace{2cm}\includegraphics[height = 0.17\textheight]{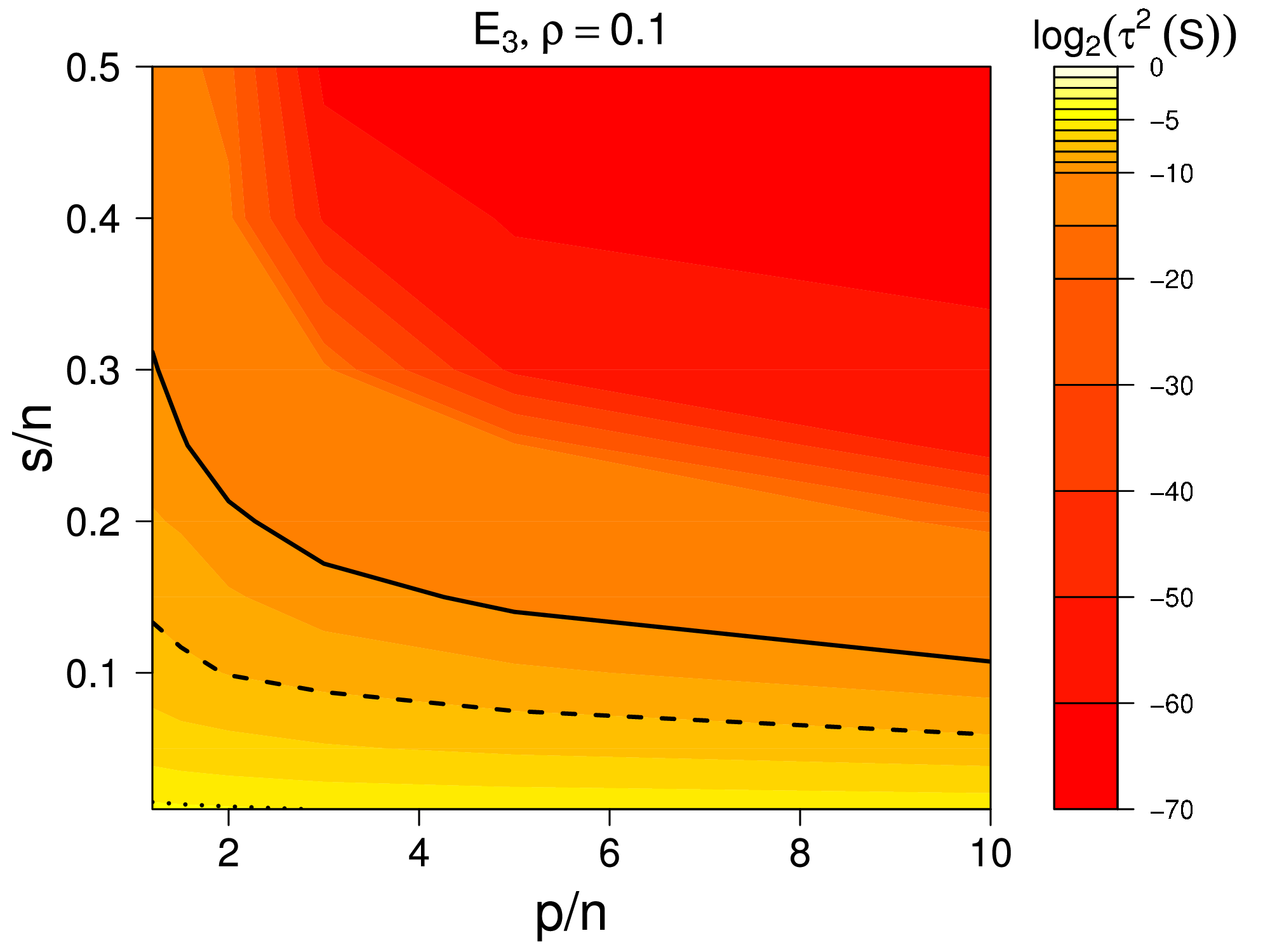}
\end{tabular}
\hfill\\
\hfill\\
\hfill\\     

\subsection*{$E_4$}

\begin{tabular}{ll}
 \hspace{.35cm}\includegraphics[height =
0.17\textheight]{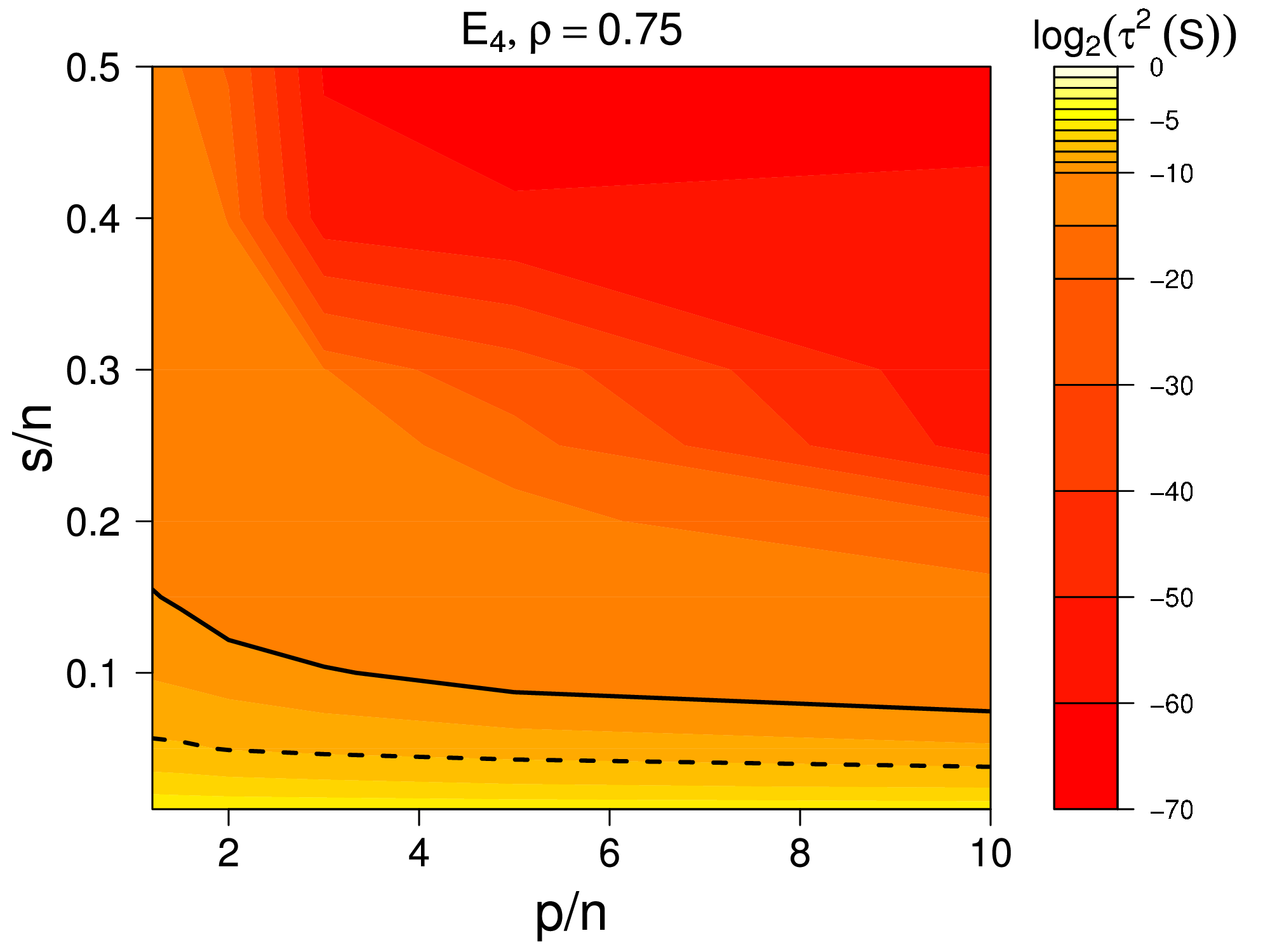}
& \hspace{2cm}\includegraphics[height =
0.17\textheight]{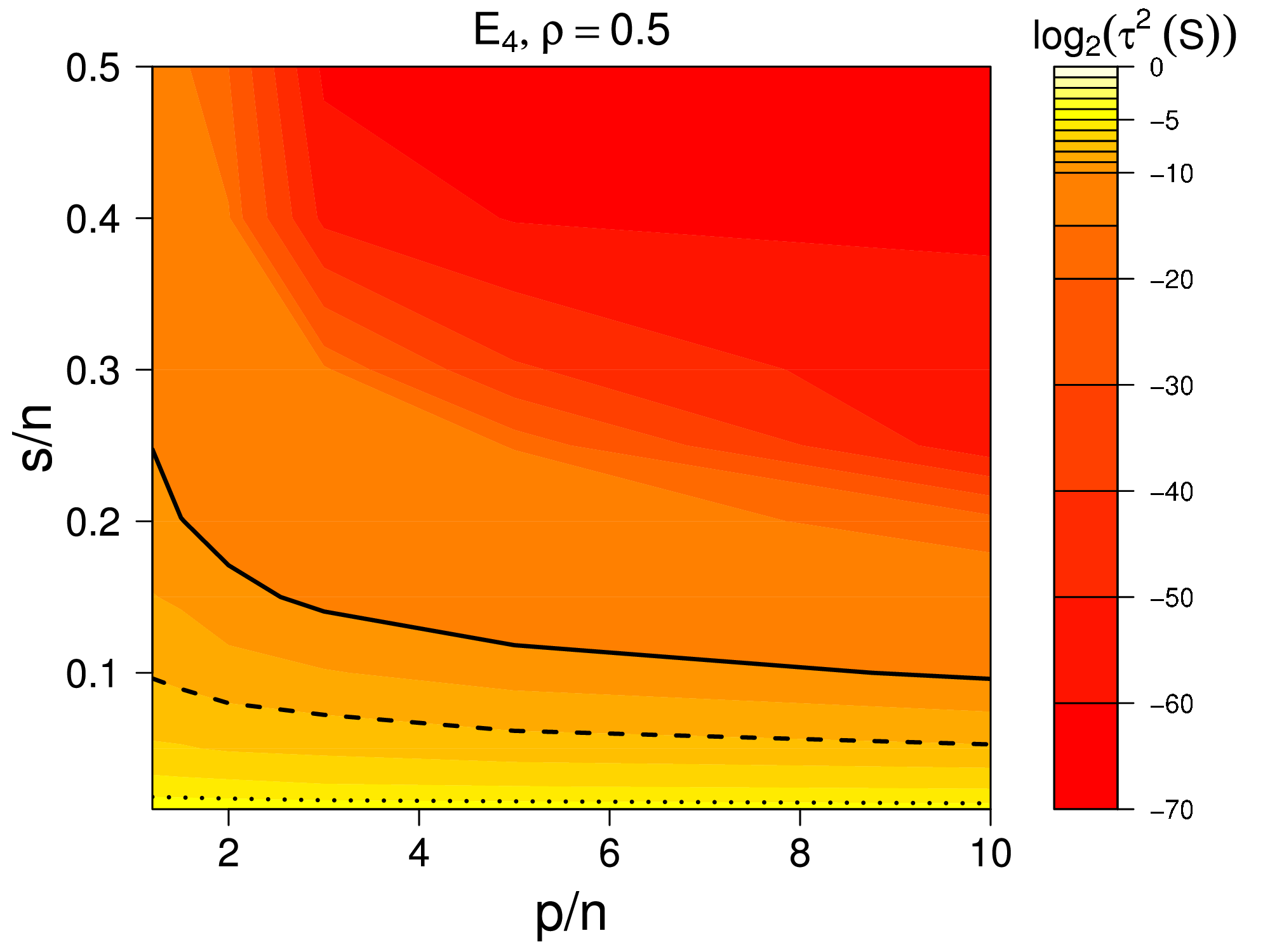} \\
  \\
\\
\\
 \hspace{.35cm}\includegraphics[height =
0.17\textheight]{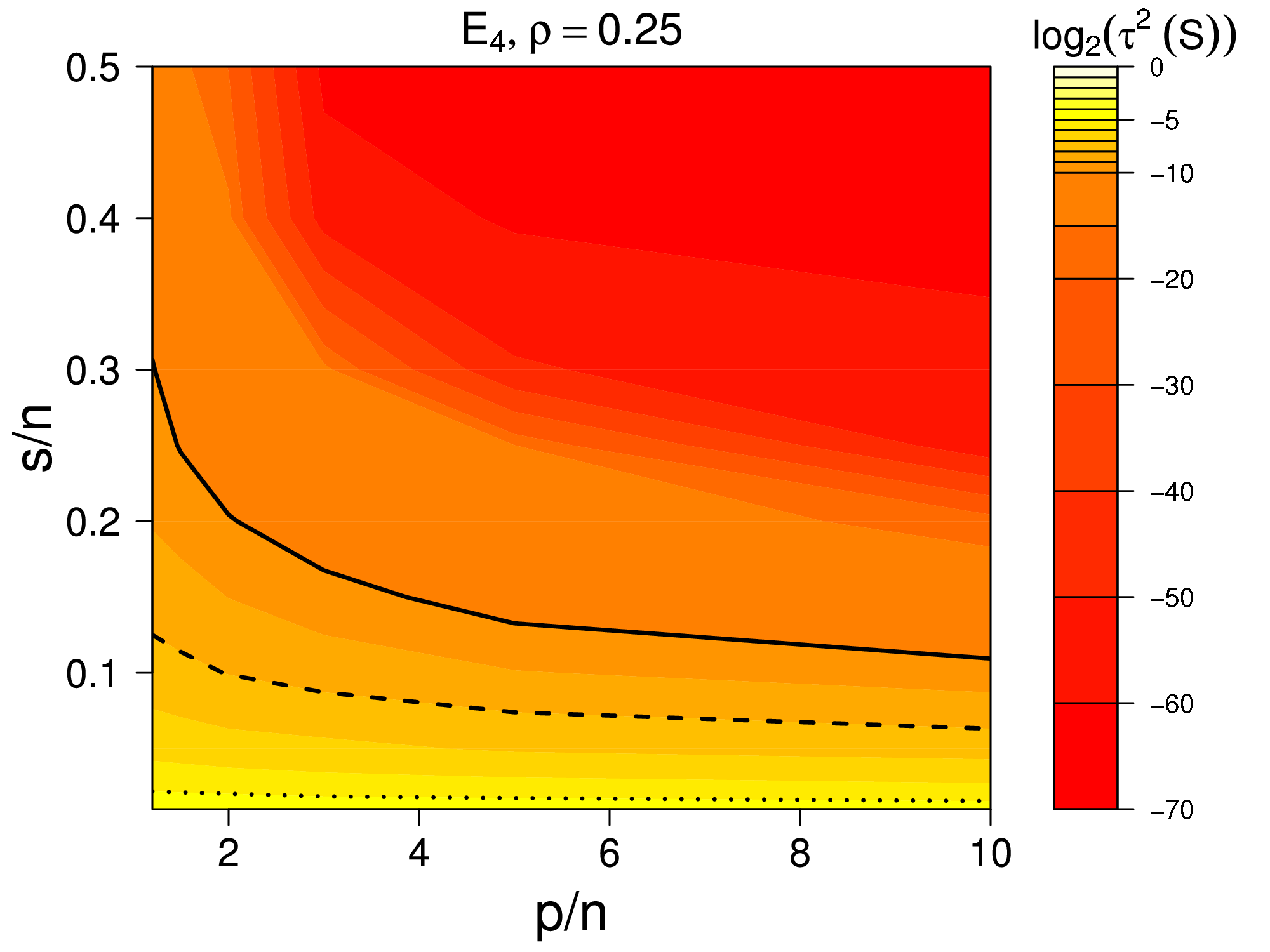} & \hspace{2cm}\includegraphics[height = 0.17\textheight]{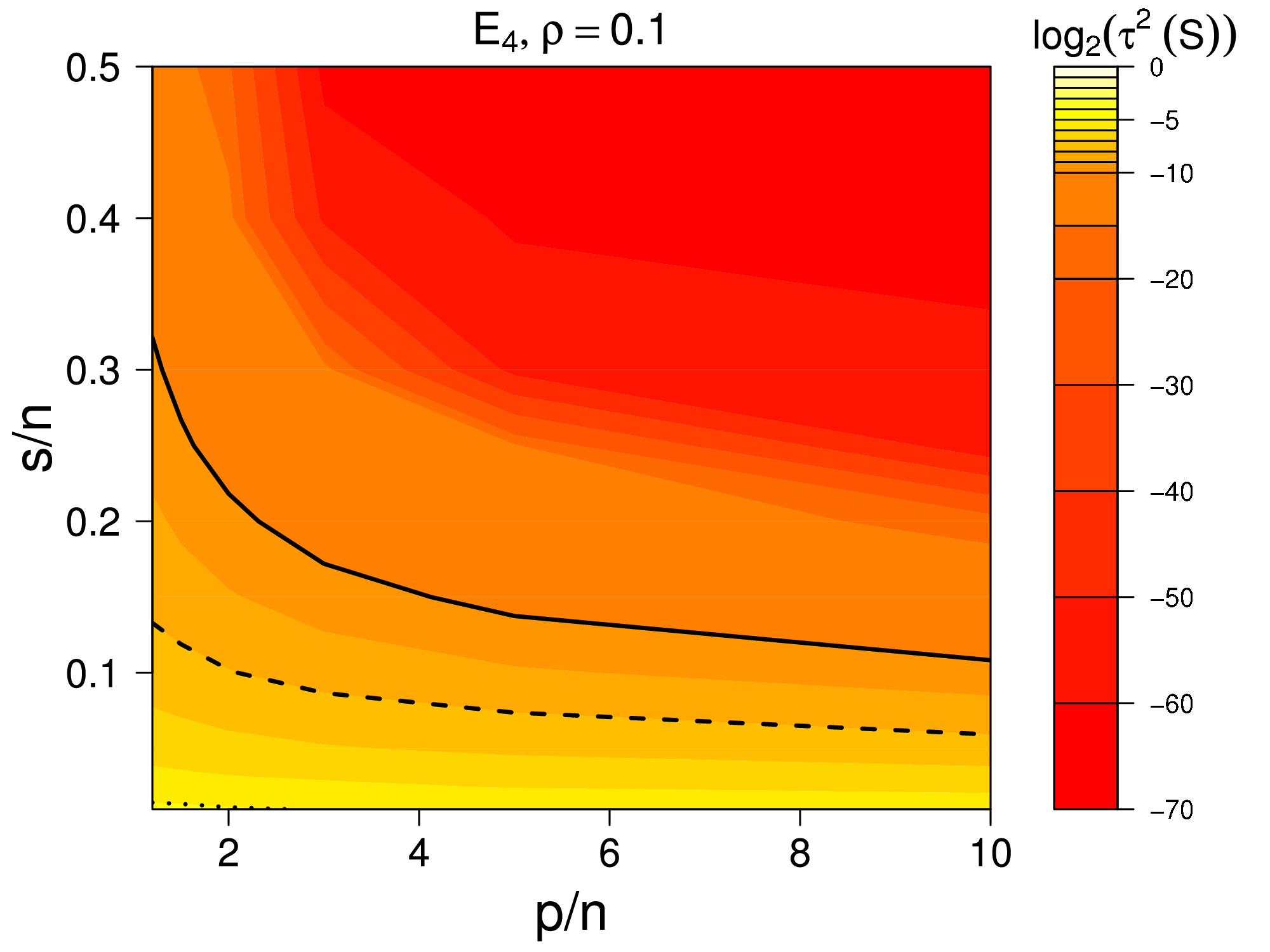}
\end{tabular}
\hfill\\
\hfill\\
\hfill\\

\section*{Additional empirical results on the $\ell_2$-error in estimating
  $\beta^*$}

The results in the two tables below are complementary to the experimental
results in Section 6.2 of the paper. We here report $\nnorm{\wh{\beta} -
  \beta^*}_2$ (NNLS) and $\nnorm{\plasso - \beta^*}_2$ (NN${\ell_1}$) in
correspondence to Tables 1 and 2 of the paper. 

\subsubsection*{Design I}

{\footnotesize
\begin{tabular}{l|l|l||l|l||l|l||l|l||}
\cline{2-9}         
& \multicolumn{8}{c|}{$p/n$} \\
\cline{2-9}
& \multicolumn{2}{c||}{$2$} & \multicolumn{2}{c||}{$3$} &
\multicolumn{2}{c||}{$5$} & \multicolumn{2}{c||}{$10$} \\ 
\cline{1-9}
\multicolumn{1}{|l|}{$s/n$} & \multicolumn{1}{|l|}{nnls}& nn${\ell_1}$ & \multicolumn{1}{|l|}{nnls}  &
nn${\ell_1}$ & \multicolumn{1}{|l|}{nnls} & nn$\ell_1$ &
\multicolumn{1}{|l|}{nnls} & nn$\ell_1$ \\
\cline{1-9}
\multicolumn{1}{|l|}{\multirow{1}{*}{$0.05$}} &  1.0\se{.01}& 1.0\se{.01} & 1.1\se{.01} & 1.1\se{.01} &
1.2\se{.01} & 1.2\se{.01} & 1.3\se{.01}& 1.3\se{.01} \\
\cline{1-9}
\multicolumn{1}{|l|}{\multirow{1}{*}{$0.1$}} & 1.4\se{.01}  & 1.4\se{.01} & 1.6\se{.01} & 1.6\se{.01} &
1.8\se{.02} & 1.8\se{.02}  & 2.1\se{.02} & 2.1\se{.02}  \\
\cline{1-9}
\multicolumn{1}{|l|}{\multirow{1}{*}{$0.15$}} &  1.8\se{.01} & 1.8\se{.02} & 2.0\se{.02} & 2.0\se{.02}
& 2.4\se{.02} & 2.4\se{.02} & 3.1\se{.04} & 3.4\se{.05}  \\
\cline{1-9}
\multicolumn{1}{|l|}{\multirow{1}{*}{$0.2$}} & 2.1\se{.02} & 2.2\se{.04} & 2.5\se{.02} & 2.6\se{.07}
& 3.1\se{.03} & 3.3\se{.04} & 5.4\se{.10} & 6.9\se{.19}  \\
\cline{1-9}
\multicolumn{1}{|l|}{\multirow{1}{*}{$0.25$}} & 2.5\se{.02} & 2.6\se{.04} & 3.1\se{.03} & 3.7\se{.14}
& 4.5\se{.07}  & 7.2\se{.27} & 12.0\se{.2} & 15.3\se{.2}  \\
\cline{1-9}
\multicolumn{1}{|l|}{\multirow{1}{*}{$0.3$}} & 3.0\se{.03} & 3.4\se{.11} & 4.0\se{.05} & 5.5\se{.24}
& 8.1\se{.19}  & 12.8\se{.3} & 18.6\se{.1} & 19.8\se{.1}  \\
\cline{1-9}
\end{tabular}
}

\subsubsection*{Design II}

{\footnotesize
\begin{tabular}{l|l|l||l|l||l|l||l|l||}
\cline{2-9}         
& \multicolumn{8}{c|}{$p/n$} \\
\cline{2-9}
& \multicolumn{2}{c||}{$2$} & \multicolumn{2}{c||}{$3$} &
\multicolumn{2}{c||}{$5$} & \multicolumn{2}{c||}{$10$} \\ 
\cline{1-9}
\multicolumn{1}{|l|}{$s/n$} & \multicolumn{1}{|l|}{nnls}& nn${\ell_1}$ & \multicolumn{1}{|l|}{nnls}  &
nn${\ell_1}$ & \multicolumn{1}{|l|}{nnls} & nn$\ell_1$ &
\multicolumn{1}{|l|}{nnls} & nn$\ell_1$ \\
\cline{1-9}
\multicolumn{1}{|l|}{\multirow{1}{*}{$0.02$}} &  0.6\se{.01}& 0.7\se{.01} & 0.6\se{.01} & 0.7\se{.01} &
0.6\se{.01} & 0.7\se{.01} & 0.6\se{.01}& 0.7\se{.01} \\
\cline{1-9}
\multicolumn{1}{|l|}{\multirow{1}{*}{$0.04$}} & 0.7\se{.01}  & 1.0\se{.01} & 0.7\se{.01} & 1.0\se{.01} &
0.7\se{.01} & 1.0\se{.01}  & 0.7\se{.01} & 1.0\se{.01}  \\
\cline{1-9}
\multicolumn{1}{|l|}{\multirow{1}{*}{$0.06$}} &  0.8\se{.01} & 1.2\se{.01} & 0.8\se{.01} & 1.2\se{.01}
& 0.8\se{.01} & 1.2\se{.02} & 0.9\se{.01} & 1.2\se{.01}  \\
\cline{1-9}
\multicolumn{1}{|l|}{\multirow{1}{*}{$0.08$}} & 0.9\se{.01} & 1.3\se{.02} & 0.9\se{.01} & 1.3\se{.02}
& 0.9\se{.01} & 1.3\se{.02} & 1.0\se{.01} & 1.4\se{.01}  \\
\cline{1-9}
\multicolumn{1}{|l|}{\multirow{1}{*}{$0.1$}} & 1.0\se{.01}  & 1.4\se{.02} &
1.0\se{.01} & 1.5\se{.02}
& 1.0\se{.01}  & 1.5\se{.02} & 1.1\se{.01} & 1.5\se{.02}  \\
\cline{1-9}
\end{tabular}
}


\end{document}